\documentclass[letter, twoside, 11pt]{article}  
\usepackage[latin1]{inputenc}
\usepackage[T1]{fontenc}
\usepackage{amssymb, amsmath, theorem, amsfonts, amscd, geometry,graphics,epsf,epsfig,psfrag,multicol,path,
changebar, array, delarray, enumerate, mathrsfs, setspace, fancyhdr, amssymb}

\usepackage[all]{xy}

\makeindex

\numberwithin{equation}{section}

\newcommand\qed{\hfill$\square$}

\newtheorem{The}{Theorem}[section]
\newtheorem{Pro}[The]{Proposition}
\newtheorem{Lem}[The]{Lemma}
\newtheorem{Def}[The]{Definition}
\newtheorem{Cor}[The]{Corollary}

\newtheorem{Rem}[The]{Remark}

\DeclareMathOperator{\End}{End}
\DeclareMathOperator{\Hom}{Hom}

\DeclareMathOperator{\Lie}{Lie}
\DeclareMathOperator{\Spec}{Spec}
\DeclareMathOperator{\Spf}{Spf}

\DeclareMathOperator{\ord}{ord}

\DeclareMathOperator{\diag}{diag}

\DeclareMathOperator{\univ}{univ}

\DeclareMathOperator{\vol}{vol}

\DeclareMathOperator{\tr}{tr}

\DeclareMathOperator{\Herm}{Herm}
\DeclareMathOperator{\charpol}{charpol}
\DeclareMathOperator{\codim}{codim}
\DeclareMathOperator{\Diff}{Diff}
\DeclareMathOperator{\supp}{supp}
\DeclareMathOperator{\red}{red}

\newcommand\V{\mathbb{V}}
\newcommand\N{\mathbb{N}}
\newcommand\Q{\mathbb{Q}}
\newcommand\R{\mathbb{R}}
\newcommand\C{\mathbb{C}}
\newcommand\A{\mathbb{A}}
\newcommand\F{\mathbb{F}}
\newcommand\Y{\mathbb{Y}}
\newcommand\X{\mathbb{X}}
\newcommand\Z{\mathbb{Z}}
\newcommand\Fp{\mathbb{F}_p}
\newcommand\Fpq{\mathbb{F}_{p^2}}

\newcommand\Zp{\mathbb{Z}_p}

\newcommand\Zpq{\mathbb{Z}_{p^2}}
\newcommand\Qp{\mathbb{Q}_p}
\newcommand\Qpq{\mathbb{Q}_{p^2}}

\newcommand\M{\mathcal{M}}
\newcommand\bk{\boldsymbol{k}}

\begin{document}
\title{Intersections of special cycles on the Shimura variety for $GU(1,2)$}
\author{Ulrich Terstiege }
\date{}
\maketitle
\thispagestyle{empty}
\section*{Introduction}
In the paper \cite{KR2} by Kudla and Rapoport, a close connection between intersection multiplicities of special cycles on arithmetic models of the  Shimura variety for $GU(1,n-1)$ and Fourier coefficients of derivatives of certain Eisenstein series is established in the case that this intersection is non-degenerate, i.e. the intersection of the cycles has dimension $0$. 
Kudla and Rapoport conjecture in that paper that their result also holds in the case of degenerate intersections. 
They show that the case of non-degenerate intersections  can be reduced to the case $n=2$. The aim of this paper is to prove this conjecture in the case $n=3$, which is the first case in which degenerate intersections of special cycles occur.
The  arithmetic models of the Shimura variety for $GU(1,n-1)$ are moduli spaces of abelian schemes with some additional structures. 
 We can regard them  as arithmetic models of moduli spaces of Picard type. In particular, for $n=3$, we are dealing with special cycles on arithmetic models of Picard surfaces. As discussed in \cite{KR2}, the intersection multiplicity of such {\em global}  special cycles can be expressed in terms of the intersection multiplicity of some {\em local} special cycles which are defined on a certain moduli space of $p$-divisible groups.  The conjecture on the intersections of the global special cycles translates into a conjectured connection between the intersection multiplicity of local special cycles and derivatives of certain hermitian representation densities.  
These local special cycles in turn were introduced by Kudla and Rapoport in \cite{KR1}, and in loc. cit. they also prove the local conjecture in the case of a non-degenerate intersection.
\newline

This paper is divided into two parts. The first part is the bulk of the paper and treats the local intersection problem for $n=3$, i.e. proves the corresponding local conjecture stated in  \cite{KR1}. In the second part we apply the main result of the first part to deduce the corresponding global result.
\newline

We  now describe the local intersection problem for $GU(1, n-1)$ as introduced in \cite{KR1}.

Let $p > 2 $ be a prime. 
 Let $\F= \overline{\F}_p$ be a fixed algebraic closure of $\Fp$ and let $W=W(\F)$ be its ring of Witt vectors.

Let $(\mathbb{X}, \iota, \lambda_{\mathbb{X}})$ be a fixed supersingular $p$-divisible group  of  dimension $n$ and height $2n$ over $\F$ which is equipped with an action $\iota: \Zpq \rightarrow \End(\mathbb{X})$ satisfying the signature condition $(1, n-1)$ (see below) and with a $p$-principal polarization $\lambda_{\mathbb{X}}$ of $\mathbb{X}$ for which the Rosati involution satisfies $\iota(a)^*=\iota(a^{\sigma}).$ The triple $(\mathbb{X}, \iota, \lambda_{\mathbb{X}})$ is unique up to isogeny. 

We consider the following functor on the category Nilp of $W$-schemes $S$ such that $p$ is locally nilpotent in $\mathcal{O}_S.$ 

It associates to a scheme $S \in $ Nilp the set of isomorphism classes of tuples $(X,\iota_X, \lambda_X, \varrho_X )$, where $X$ is a $p$-divisible group over $S$ which is equipped with an $\Zpq$-action $\iota_X$ satisfying the   
 signature condition $(1, n-1),$  where $\lambda_X$ is a $p$-principal polarization on $X$ for which the Rosati involution induces the Galois involution $\sigma $ on $\Zpq$, and where $\varrho_X$ is an $\Zpq$-linear quasi-isogeny of height zero 
$$
\varrho_X: X\times_S \overline{S} \rightarrow \mathbb{X}\times_{\F} \overline{S}.
$$
Here $\overline{S}= S\times_W \F$. Further, we require that up to a scalar in $\Zp^{\times}$ we have the identity $\hat{\varrho_X}\circ\lambda_{\mathbb{X}}\circ \varrho_{X} = \lambda_X$.
Recall that the signature condition $(1,n-1)$ means the following. Let $\phi_0$ and $\phi_1$  be the two embeddings of $\Zpq$ into $W$. Then we require that the characteristic polynomial of $\iota(a)$ is of the form
$$
\charpol (\iota(a), \Lie X)(T)=(T-\phi_0(a))(T-\phi_1(a))^{n-1} \in \mathcal{O}_S [ T ] .
$$
This functor is representable by a separated formal scheme $\mathcal{N}$ over $\Spf W$ which is locally formally of finite type over $W$ and formally smooth of dimension $n-1$ over $W.$ (More generally, considering the same functor but with  signature condition $(n-r,r)$ for any $0\leq r \leq n$, this functor   is representable by a separated formal scheme over $\Spf W$ which is locally formally of finite type over $W$ and formally smooth of dimension $(n-r)\cdot r$ over $W.$)

To define the notion of  a {\em special cycle},  we first need to know what a {\em special homomorphism} is. 
Let $(\mathbb{Y}, \iota, \lambda_{\mathbb{Y}})$ be the basic triple over $\F$ used for the definition of $\mathcal{N}$ in the case $n=1.$ Let $(\overline{\mathbb{Y}}, \iota, \lambda_{\overline{\mathbb{Y}}})$ be the triple obtained from  $(\mathbb{Y}, \iota, \lambda_{\mathbb{Y}})$ by changing $\iota$ to $\iota \circ \sigma.$ 
 The pair  $(\overline{\mathbb{Y}}, \iota)$ has a canonical lift  $(\overline{Y}, \iota)$ over $W.$ 
The space of special homomorphisms is defined to be the $\Qpq$-vector space $$\mathbb{V}:=\Hom_{\Zpq}(\overline{\mathbb{Y}}, \mathbb{X})\otimes_{\Z}\Q.$$
It is equipped with an hermitian form given by $$h(x,y)= \lambda^{-1}_{\overline{\Y}}\circ \hat{y}\circ \lambda_{\X}\circ x \in \End_{\Zpq}(\overline{\Y})\otimes \Q \cong \Qpq,$$
where the last isomorphism is via $\iota^{-1}$, and where $\hat{y}$ is the dual of $y.$
For a special homomorphism $x\in \V$, we define the {\em special cycle} $\mathcal{Z}(x)$ to be the closed formal subscheme of $\mathcal{N}$ with the property that, for $S \in$ Nilp, the set  $\mathcal{Z}(x)(S)$ consists of all points $(X,\iota_X, \lambda_X, \varrho_X )$ in $\mathcal{N}(S)$ where the quasi-homomorphism
$$
\overline{\Y}\times_{\F} \overline{S}\stackrel{x}{\longrightarrow}{\X}\times_{\F} \overline{S}\stackrel{\varrho_X^{-1}}{\longrightarrow}X\times_S \overline{S}
$$
extends to an $\Zpq$-linear homomorphism $$\overline{Y}\times_W S \longrightarrow X.$$
Then $\mathcal{Z}(x)$ is (for $x \neq 0$) a relative divisor in $\mathcal{N}$ (or empty).

For a collection $j_1,...,j_m$ of special homomorphisms, the {\em fundamental matrix} $T(j_1,..,j_m)$ is defined to be
$$T(j_1,...,j_m)=(h(j_k, j_l))\in \Herm_m(\Qpq).$$

Suppose we are given  special homomorphisms  $j_1,...,j_n$ such that  the intersection $\mathcal{Z}(j_1)\cap...\cap \mathcal{Z}(j_n)$ is non-empty and the fundamental matrix  is non-singular. We define the  intersection multiplicity of the  special cycles $\mathcal{Z}_i= \mathcal{Z}(j_i)$
to be 
$$
(\mathcal{Z}_1,..., \mathcal{Z}_n)= \chi(\mathcal{O}_{\mathcal{Z}_1} \otimes^{\mathbb{L}}...\otimes^{\mathbb{L}}\mathcal{O}_{\mathcal{Z}_n} ),
$$
where $\otimes^{\mathbb{L}}$ is the derived tensor product of $\mathcal{O}_{\mathcal{N}}$-modules and $\chi$ is the Euler-Poincar\'e characteristic.
It follows from our assumptions on $j_1,...,j_n$ that the intersection $\bigcap_i Z(j_i)$ has support in the supersingular locus and that this support is proper over $\F$. Therefore $(\mathcal{Z}_1,..., \mathcal{Z}_n)$ is a finite number. 

We want to relate in the case 
 $n=3$ the intersection multiplicity $(\mathcal{Z}_1, \mathcal{Z}_2, \mathcal{Z}_3)$ to certain hermitian representation densities which we now recall. 

For non-singular hermitian matrices $S\in \Herm_m(\Zpq)$ and $T\in \Herm_l(\Zpq),$ where $m \geq l$, the representation density $\alpha_p(S,T)$ is defined as $$\alpha_p(S,T)= \operatorname*{lim}_{k\rightarrow\infty} p^{-kl(2m-l)} \mid \{x \in M_{m,l}(\Zpq/p^k \Zpq);
\ S[x] \equiv T \mod p^k \}\mid,$$
where $S[x]= {^t} xS\sigma(x).$
Note that  $\alpha_p(S,T)$ depends only on the $GL_m(\Zpq)$- resp. $GL_l(\Zpq)$-equivalence classes of  $S$ and $T$.
For $r \geq 0$ let $S_r= \diag(S, 1_r).$ Then one can show that there is a polynomial $F_p(S,T;X)\in \Q[X]$ such that $\alpha_p(S_r, T)=F_p(S,T; (-p)^{-r}).$ The derivative   $\alpha^{'}_p(S,T)$ is defined to be
$$\alpha^{'}_p(S,T):= -\frac{\partial}{\partial X}F_p(S,T;X)\mid_{X=1}. $$

The following theorem is the main result of the local theory in this paper and confirms the local conjecture of Kudla and Rapoport in the case $n=3$. 
Suppose  for the remainder of the Introduction that $n=3$ (unless otherwise mentioned) and let us fix  special homomorphisms  $j_1, j_2, j_3$ such that  the intersection $\mathcal{Z}(j_1)\cap \mathcal{Z}(j_2) \cap \mathcal{Z}(j_3)$ is non-empty and the fundamental matrix  is non-singular.
\begin{The}\label{HSeinf}
Let $S=1_3$ be the $3\times 3$ unit matrix and let $T$ be the fundamental matrix of the fixed special homomorphisms $j_1,j_2,j_3$. Then there the following identity between intersection multiplicities of special cycles and representation densities, 
 $$
(\mathcal{Z}_1, \mathcal{Z}_2, \mathcal{Z}_3)=\frac{\alpha^{'}_p(S,T)}{\alpha_p(S,S)}.
 $$  
\end{The}
 In the case that the intersection has dimension $0$ or, equivalently, that $T$ is $GL_3(\Zpq)$-equivalent to a diagonal matrix of the form $\diag(1,p^a,p^b)$, this theorem was proved  in \cite{KR1}. (More generally, in \cite{KR1} it is shown for any $n$ that the intersection is non-degenerate  if and only if   $T$ is $GL_n(\Zpq)$-equivalent to a diagonal matrix of the form $\diag(1_{n-2},p^a,p^b)$ and that the  analogue of the above formula holds in that case.)

Sections 1-5 deal with the proof of this theorem. The proof is by explicit calculation of both sides of the equation. In section 5 an explicit recursive formula for the relevant representation densities is derived (see below). It works for general $n$. Sections 1 to 4 deal with the geometry of special cycles and the calculation of their intersection multiplicities. 
 A main ingredient  is an explicit description of the special fiber of special cycles.
Before we formulate  it, we first briefly summarize some results of Vollaard  on the underlying reduced subscheme of $\mathcal{N}$: To the isocrystal $N$ of $\mathbb{X}$ one associates an hermitian space $C$ of dimension $n=3$ over $\Qpq$. Then there is a bijection $\Lambda \mapsto \mathcal{V}(\Lambda)$ between  the set of self-dual $\Zpq$-lattices in $C$ and the set irreducible components (of the underlying reduced subscheme) of  $\mathcal{N}$. Each $\mathcal{V}(\Lambda)$ is isomorphic to the Fermat curve in $\mathbb{P}^2_{\F}$ given by $X^{p+1}+Y^{p+1}+Z^{p+1}=0.$ The underlying reduced subscheme of  $\mathcal{N}$ can be written as union of the  $\mathcal{V}(\Lambda)$ in such a way that $\mathcal{N}_{\red}$ becomes a tree whose vertices are the intersections of the curves (see \cite{V} for details and proofs of these statements).
In section 2, we will first give a description of the underlying reduced subscheme of a special cycle. In particular, we will see that  it is connected. Let now $S(j)$ be the set of lattices $\Lambda$ such that $\mathcal{V}(\Lambda)$ is contained in $Z(j)$  (see Proposition \ref{redloc}), and denote by $b_j(\Lambda)$ the maximal integer $b$ such that $\Lambda \in S(j/p^b)$.  We now give a description of the special fiber of a special cycle as a divisor in the special fiber  $\mathcal{N}_p$ of $\mathcal{N}$. Note that $\mathcal{N}_p$ is formally smooth of dimension $2$ over $\F$. 
\begin{The} \label{spezfintro}
Let $j\in \mathbb{V}$ such that  $r:=\ord_p (h(j,j))\geq 0$. The special  fiber of $\mathcal{Z}(j)$ as a divisor in $\mathcal{N}_p$ can be described as follows,
$$
\mathcal{Z}(j)_p= p^{r}\cdot t + \sum_{\Lambda \in S(j)}(1+p^2+...+ p^{2b_j(\Lambda)})\cdot \mathcal{V}(\Lambda).
$$
Here $t$ is zero for $r$ odd. For $r$ even, it is an irreducible divisor in $\mathcal{N}_p$ which passes through the unique $\F$-valued point of $\mathcal{Z}(j/p^{r/2}).$ Its intersection multiplicity with each $\mathcal{V}(\Lambda)$ which contains this point is $1$.
\end{The}
To prove this theorem, we will use the theory of displays and windows, see \cite{Z1}, \cite{Z2}. 
This theorem is the key to determining equations for the intersections of two special cycles, as done in section 3. We will proceed as follows.
 To $\mathcal{Z}(j)$ we associate (as in \cite{T}) the "difference divisor" $\mathcal{D}(j)=\mathcal{Z}(j)-\mathcal{Z}(j/p)$. We will show that $\mathcal{D}(j)$ is regular. Further, we will see that, given $j_1,j_2,j_3$ as in Theorem \ref{HSeinf}, we have $$
(\mathcal{Z}(j_1),\mathcal{Z}(j_2),\mathcal{Z}(j_3))=\sum_{l_1,l_2,l_3}(\mathcal{D}(j_1/p^{l_1}),\mathcal{D}(j_2/p^{l_2}),\mathcal{D}(j_3/p^{l_3})).
$$ Therefore, by induction, it is enough to compute the intersection multiplicity $ (\mathcal{D}(j_1),\mathcal{D}(j_2),\mathcal{D}(j_3))$ which can be written as the intersection multiplicity of $\mathcal{D}(j_1)\cap \mathcal{D}(j_2)$ and $\mathcal{D}(j_1)\cap \mathcal{D}(j_3)$ as divisors in the regular (formal) scheme  $\mathcal{D}(j_1)$. Thus we will focus in section 3 on determining equations for the intersections $\mathcal{D}(j_k)\cap \mathcal{D}(j_l)$. The methods for this will be a combination of  Theorem \ref{spezfintro}, of Grothendieck-Messing theory, and of the results of \cite{KR1} in the case of a non-degenerate intersection. For example, we will use Grothendieck-Messing theory to show that, if locally around a point $x$ we have  $\mathcal{D}(j_1)_p\subset \mathcal{D}(j_2/p)_p$, then locally around $x$ we have  $\mathcal{D}(j_1)\cap \mathcal{D}(j_2)=\mathcal{D}(j_1)_p$,  whose structure in turn is known by Theorem \ref{spezfintro} (compare Lemma \ref{pmultlem}). 

Another important property of the intersection multiplicity $(\mathcal{Z}(j_1),\mathcal{Z}(j_2),\mathcal{Z}(j_3))$ 
contained in the claim of Theorem \ref{HSeinf} is that it only depends on the $\Zpq$-span of $j_1,j_2,j_3$ in $\mathbb{V}$. We will see (in Proposition \ref{diag}) that even  $ \mathcal{O}_{\mathcal{Z}(j_1)}\otimes^{\mathbb{L}}\mathcal{O}_{\mathcal{Z}(j_2)}\otimes^{\mathbb{L}}\mathcal{O}_{\mathcal{Z}(j_3)}$ depends only on the $\Zpq$-span  of $j_1,j_2,j_3$ in $\mathbb{V}$. 
This allows us,  for example, to assume that $j_1,j_2,j_3$ are perpendicular to each other, which simplifies the combinatorics. 

Having determined equations for the intersections $\mathcal{D}(j_k)\cap \mathcal{D}(j_l)$ the calculation of the intersection multiplicities will be a combinatorial problem solved in section 4.  The result is the following. 
 \begin{The}\label{intmulteinf}
 Suppose that the fundamental matrix of the fixed special homomorphisms $j_1,j_2,j_3$ is $GL_3(\Zpq)$-equivalent to the diagonal matrix $diag(p^{a_1},p^{a_2},p^{a_3})$, where $0\leq a_1 \leq a_2 \leq a_3$. Then the intersection multiplicity $(\mathcal{Z}(j_1),\mathcal{Z}(j_2),\mathcal{Z}(j_3))$ is finite and is given by the formula
$$
(\mathcal{Z}(j_1),\mathcal{Z}(j_2),\mathcal{Z}(j_3))= -\frac{1}{2}\sum_{k=0}^{a_1}\sum_{l=0}^{a_1+a_2-2k}(-1)^k((k+l)p^{2k+l}-(k+l+a_3+1)p^{a_1+a_2-l}).
$$
\end{The}

The restriction to the case $n=3$ simplifies the problem in several aspects. For example, the underlying reduced subscheme of $\mathcal{N}$ is in this case a tree, which simplifies the combinatorics. Also, the fact that  $\mathcal{N}_{\red}$ is one-dimensional and the fact that we can reduce the problem to the calculation of the intersection multiplicity of two divisors on a (formal) arithmetic curve (namely $\mathcal{D}(j_1)$) simplifies the situation considerably.

Finally, we need to calculate the relevant hermitian representation densities. 
For any hermitian $m \times m$ matrix $T$  over $\Zpq$,  let $F_p(X;T)$ be the polynomial such that $\alpha_p(1_s,T)=F_p((-p)^{-s};T).$ In section 5 we will prove the following inductive statement.
\begin{The}\label{repeinf}
Suppose that  the matrix   $T\in \Herm_n(\Zpq)$ is $GL_m(\Zpq)$-equivalent to the diagonal matrix $diag(p^{a_1},...,p^{a_m})$, where  $0 \leq a_1 \leq ... \leq a_m.$ 
\begin{enumerate}
\item If $m=1$, then 
$$F_p(X;T)=(1-X)\sum_{l=0}^{a_1}(pX)^l. $$
\item If  $m\geq 2$, let $T^{-}= diag(p^{a_1},...,p^{a_{m-1}}).$ Further, let $H_p(X)=(1-X)(1+pX)$ and let $$A_p(X;T)=(1-(-p)^{m-1}X)(1-(-p)^{-m}X^{-1})(-p)^{a_1+...+a_{m-1}}(-1)^{(m+1)a_m}(p^mX)^{a_m+2}.$$ Then there is the following inductive formula, 
$$
F_p(X;T)=\frac{H_p(X)F_p(p^2X;T^{-})-A_p(X;T)F_p(X;T^{-})}{1-p^{2m}X^2}.
$$
\end{enumerate}
\end{The}

A similar statement for quadratic representation densities was proved by Katsurada in \cite{Ka}. The proof here is also similar to that in \cite{Ka}. More precisely, we will combine a weaker recursion formula for the polynomials $F_p(X;T)$ with a functional equation for  those polynomials obtained from the paper \cite{I} by Ikeda. Combining Theorem 
\ref{intmulteinf} and Theorem \ref{repeinf}, we obtain Theorem \ref{HSeinf}.
\newline

In the last two sections, we will apply the local result of Theorem \ref{HSeinf} to the global intersection problem introduced in \cite{KR2} which we describe now (see loc. cit. and section 6 for further details).
For the moment, let again $n$ be arbitrary.
Let $\bk$ be an imaginary quadratic field with ring of integers $O_{\bk}$. Let $n\geq 1$ and let $0 \leq r \leq n$. We consider the Deligne-Mumford-stack $\M(n-r,r)$ over  $O_{\bk}$ whose functor of points 
associates to any 
   locally noetherian $O_{\bk}$-scheme $S$ the 
   groupoid of tuples $(A,\iota, \lambda)$, where $A$ is an abelian scheme over $S$, where $\iota:  O_{\bk} \rightarrow \End_S(A)$ is an  $O_{\bk}$-action on $A$ satisfying the signature condition $(n-r,r)$ (plus another technical condition, see section 6), and where $\lambda: A \rightarrow A^{\vee}$ is a principal polarization for which the Rosati involution satisfies $\iota(a)^*=\iota(a^{\sigma}).$ See \cite{KR2}, section 4 (or section 6 of this paper) for the relation of $\M(n-r,r)$ to arithmetic models of the Shimura variety for $GU(n-r,r)$.  We consider 
 the fiber product
$$\M=\M(n-r,r)\times_{\Spec O_{\bk}} \M(1,0).$$
Given a point in $\M(S)$, i.e. a pair $(A,\iota, \lambda) \in \M(n-r,r)(S)$ and $(E,\iota_0, \lambda_0) \in \M(1,0)(S)$, we consider the free $O_{\bk}$-module 
$$V^{'}(E,A)=\Hom_{O_{\bk}}(E,A).$$
It is equipped with an $O_{\bk}$-valued  hermitian form $h^{'}$ given by $$h^{'}(x,y)=\iota_0^{-1}(\lambda_0^{-1}\circ y^{\vee}\circ \lambda \circ x),$$
where $y^{\vee}$ is the dual of $y.$ Given elements $x_1,...,x_m \in V^{'}(E,A)$,  the {\em fundamental matrix} of the collection $x_1,...,x_m$ is defined to be the hermitian $m \times m$ matrix over $ O_{\bk}$ with entries $h^{'}(x_i,x_j).$
\newline
For an hermitian $m \times m$ matrix $T$ with entries in $O_{\bk}$, the  (global) {\em special cycle} $Z(T)$ is the stack of collections $(A,\iota, \lambda, E,\iota_0, \lambda_0, x_1,...,x_m)$, where $(A,\iota, \lambda) \in \M(n-r,r)(S),$ $(E,\iota_0, \lambda_0) \in \M(1,0)(S)$ and  $x_1,...,x_m \in V^{'}(E,A)$ such that the fundamental matrix of the $x_i$ is $T.$ 
Then $Z(T)$ is representable by a Deligne-Mumford-stack which is finite and unramified over $\M.$ We suppose from now on that $r=1.$

Let $T\in \Herm_n(O_{\bk})_{>0}$ with $Z(T)\neq\emptyset$. Then the support of  $Z(T)$ is contained in a union over finitely many $p$  of the supersingular locus of the fiber at $p$ of $\M$. One defines  
$$
\widehat{\deg}\ Z(T)=
\sum_p \chi(Z(T)_p, \mathcal{O}_{Z(t_1)}\otimes^{\mathbb{L}}...\otimes^{\mathbb{L}}\mathcal{O}_{Z(t_n)})\cdot \log p,
$$ where $\chi$ is the Euler-Poincar\'e characteristic, $Z(T)_p$ is the part of $Z(T)$ with support in the fiber at $p$, and $t_1,...,t_n$ are the diagonal entries of $T$.

Suppose now further that  $p$ is a prime which is inert in $\bk$ and for which $\ord_p ( \det T)$ is odd. As shown in \cite{KR2}, it follows that $p$ is the only prime with that property (since otherwise $Z(T)$ is empty) and  $Z(T)$ then  has support in the supersingular locus of the fiber at $p$ of $\M$. 
(If there is no such prime, then $Z(T)$ has support in the union of the fibers at the ramified primes, see \cite{KR2}, Proposition 2.22). Thus we
 have
$$
\widehat{\deg}\ Z(T)=
\chi(Z(T)_p, \mathcal{O}_{Z(t_1)}\otimes^{\mathbb{L}}...\otimes^{\mathbb{L}}\mathcal{O}_{Z(t_n)})\cdot \log p.
$$ 
Now let  $n=3$  and $p>2$. In the global part of the paper we relate  $\widehat{\deg}\ Z(T)$  to the $T$-th Fourier coefficient of the derivative of 
 a certain (incoherent) Eisenstein series $E(z,s)$ for $U(3,3)$ as introduced in \cite{KR2}. We will use Theorem \ref{HSeinf} to show the following Theorem.

\begin{The}\label{globHSeinf} Let $n=3$, 
let $T\in \Herm_3(O_{\bk})_{>0}$ with $Z(T)\neq\emptyset$, and let $p>2$ be a prime which is inert in $\bk$ and for which $\ord_p ( \det T)$ is odd. 
 Then
$$
E^{'}_T(z,0)=C_1\cdot \widehat{\deg}\ Z(T) \cdot q^T,
$$ for some explicit constant $C_1.$
\end{The}

Similarly as Theorem \ref{HSeinf}, this theorem was proved in the case that $Z(T)$ has dimension $0$ in \cite{KR2}. The proof of Theorem \ref{globHSeinf}, given in section 7, shows that, in order to prove the corresponding statement for any $n$, it is enough to show the corresponding statement of Theorem \ref{HSeinf}. For this in turn one can use Theorem \ref{repeinf} for the representation densities.
\smallskip
\smallskip

{\em Acknowledgements.} I want to thank all people who helped me to write this paper. In particular, I thank M. Rapoport for  suggesting  the problem  and for his continuous interest in in my work. I also thank  B. Gross, B. Howard, H. Katsurada, S. Kudla and T. Wedhorn for helpful comments and discussions. This paper has been written during my stay at the Harvard University, and I  want to thank the Department of Mathematics for providing such a great working environment. This work was supported by a fellowship within the Postdoc-Programme of the German Academic Exchange Service (DAAD).  I  thank  the DAAD and the Hausdorff Center for Mathematics for their support of my work.
\tableofcontents
\newpage
\centerline
{\bf\Large Part I:  Local Theory }

\section{On the structure  of  $\mathcal{N}$}

In this section we gather together several facts on the structure of the moduli space $\mathcal{N}$ described in the introduction, in particular about the reduced locus of $\mathcal{N}$. Proofs can be found in \cite{V} resp. \cite{VW}. 

 Let $\mathbb{M}$ be the Dieudonn\'e module of $\X$, and  let $N$ be its isocrystal. Then $\mathbb{M}$ and hence also $N$ are equipped with $\sigma$- resp. $\sigma^{-1}$-linear automorphisms $F$ resp. $V$. From the $\Zpq$-action we get gradings $\mathbb{M}=\mathbb{M}_0\oplus \mathbb{M}_1$ resp. $N=N_0\oplus N_1$. 
From the $p$-principal polarization we get a perfect alternating pairing 
$$
\langle \ , \  \rangle: N \times N \rightarrow W_{\Q}
$$
such that $ \langle Fx, y \rangle =\langle x, Vy \rangle^{\sigma}$, and, for $a \in \Zpq$, we have  $\langle \iota(a)x, y \rangle =\langle x, \iota(a)^{\sigma} y \rangle$.
Let $\tau=V^{-1}F$ and let $C=N_0^{\tau}$. It is a  $\Qpq$-vector space of dimension $n$. We define the skew hermitian form  $\{  , \}$ on $C$  given by $\{ x ,y \}=\langle x, Fy \rangle$. The same formula gives a form on $C_{W_{Q}}$ (which is not skew hermitian any more).
For a $W$-lattice $A\subset C_{W_{\Q}}$,  let $A^{\vee}$ denote the dual lattice.
There is an bijection between $\mathcal{N}(\F)$ and the set 
$$
\{  
W\text{-lattices } A\subset C_{W_{\Q} } | \ pA^{\vee} \overset{1}{\subset} A  \overset{n-1}{\subset} A^{\vee}
\}. $$
If  $x\in \mathcal{N}(\F)$ and if $M$ is its Dieudonn\'e module, then the corresponding lattice in $ C_{W_{\Q} } $ is given by $M_0$.

Let $\mathcal{L}$ be the set of $\Zpq$-lattices $\Lambda$ in C which satisfy $$p\Lambda^{\vee} \subsetneqq \Lambda  \subset  \Lambda^{\vee}.$$ For $\Lambda \in \mathcal{L}$, we set $$\mathcal{V}(\Lambda)(\F)=\{ \text{$W$-lattices }   A\subset \Lambda_{W} | \ pA^{\vee} \overset{1}{\subset} A      \}.$$
The 
type of  $\Lambda \in \mathcal{L}$ is defined to  be the index of $p\Lambda^{\vee} $ in $\Lambda.$ We denote by $\mathcal{L}^{(l)}$ the set of elements in   $\mathcal{L}$ which have type $l$. The type of  $\Lambda \in \mathcal{L}$ is always an odd integer between $1$ and $n$. 

If $x\in \mathcal{N}(\F)$ and if $M$ is its Dieudonn\'e module, we say that $x$ is superspecial if $\tau M=M.$ We then also say that  $M$, resp. the corresponding lattice $M_0$ in $ C_{W_{\Q} }$ are superspecial. 
The following facts  can be found in \cite{V} and \cite{VW}.

\begin{enumerate}
\item $ \mathcal{N}(\F)= \underset{\Lambda \in \mathcal{L}}\bigcup \mathcal{V}(\Lambda)(\F)$. 
\item $ \mathcal{V}(\Lambda)(\F) \subset  \mathcal{V}( \Lambda^{'})(\F) $ if and only if $\Lambda \subset \Lambda^{'}.$ In particular  $ \mathcal{V}(\Lambda)(\F) =  \mathcal{V}( \Lambda^{'})(\F) $ if and only if $\Lambda = \Lambda^{'}.$
\item $ \mathcal{V}(\Lambda)(\F)\cap  \mathcal{V}( \Lambda^{'})(\F)
=\begin{cases}
 \mathcal{V}(\Lambda \cap \Lambda^{'})(\F) ,  & \text{if  } \Lambda \cap \Lambda^{'} \in \mathcal{L},\\
 \emptyset, & \text{otherwise.}\
\end{cases}$
\item The set $\mathcal{V}(\Lambda)(\F)$ always contains a superspecial point.
\item The cardinality of $\mathcal{V}(\Lambda)(\F)$ is $1$ if and only if the type of $\Lambda$ is $1$. 
\item For  any $\Lambda \in \mathcal{L}^{(l)}$, the set  $\mathcal{V}(\Lambda)(\F)$ is the set of $\F$-valued points of an irreducible smooth subvariety $\mathcal{V}(\Lambda)$ of  $\mathcal{N}_{\red}$ over $\F$. Its  dimension is $\frac{1}{2}(l-1)$. 
\item The inclusions  $ \mathcal{V}(\Lambda)(\F) \subset  \mathcal{V}( \Lambda^{'})(\F)$, for  $\Lambda \subset \Lambda^{'}$, are induced by inclusions of algebraic varieties $ \mathcal{V}(\Lambda) \subset  \mathcal{V}( \Lambda^{'}) $ (as subvarieties of $\mathcal{N}_{\red}$).  
\item All irreducible components of $\mathcal{N}_{\red}$ are of the form  $ \mathcal{V}(\Lambda)$ for some $\Lambda$ of maximal type.
\item The set $\mathcal{L}$ can be identified with the set of vertices in the building $\mathcal{B}(SU(N_0, \{,\}), \Qp)$, and this identification is $SU(N_0, \{,\})(\Qp)$ invariant.
\item In the case $n=3$, for each $\Lambda \in \mathcal{L}^{(3)}$, the   variety $ \mathcal{V}(\Lambda)$ is isomorphic to the Fermat curve in $\mathbb{P}^2_{\F}$ given by $X^{p+1}+Y^{p+1}+Z^{p+1}=0.$
\item Suppose we are given $\Lambda \in \mathcal{L}^{(l)}$, let $V=\Lambda/p\Lambda^{\vee}$, and let $V^{'}=\Lambda^{\vee}/\Lambda$. These are $\Fpq$-vector spaces of dimension $l$ resp. $n-l$. The skew hermitian form $\{ ,\}$ on $C$ induces a skew hermitian form $(,)$ on $V$ by setting 
$$
(\overline{x}, \overline{y})=\overline{\{x,y\}}\in \Fpq,
$$
 for $\overline{x}, \overline{y}\in V$ and lifts $x,y\in \Lambda.$ We can extend this form to $V \otimes \F$ via $(x \otimes a, y\otimes b )=ab^{\sigma}(x,y).$
 
 Similarly, $V^{'}$ is equipped with a skew hermitian form $(\ ,\ )^{'}$ given by
  $$
(\overline{x}, \overline{y})^{'}=\overline{(p\{x,y\})}\in \Fpq,
$$
 for $\overline{x}, \overline{y}\in V^{'}$ and lifts $x,y\in \Lambda^{\vee}.$ The following assertions hold.
\begin{itemize}
\item The set of lattices $\Lambda_1\in \mathcal{L}^{(l_1)}$ with $\Lambda_1 \subset \Lambda$ can be identified with the set of $\Fpq$-subspaces $U \subset V$ of dimension $(l+l_1)/2$ with $U^{\perp}\subset U.$
The subspace in $V$ belonging to $\Lambda_1$ is $\Lambda_1/p\Lambda^{\vee}.$ The corresponding description  for the set $\mathcal{V}(\Lambda)(\F)$ in terms of the subspaces $U \subset V \otimes \F$ which have dimension  
 $(l+1)/2$  and for which $U^{\perp}\subset U$ also holds.
\item Similarly,  the set of lattices $\Lambda_1\in \mathcal{L}^{(l_1)}$ with $\Lambda \subset \Lambda_1$ can be identified with the set of $\Fpq$-subspaces $U \subset V^{'}$ of dimension $n-(l+l_1)/2$ with $U^{\perp^{'}}\subset U.$
\end{itemize}

\end{enumerate}
In particular, for $n=3$, we have the following description of  $\mathcal{N}_{\red}$. 
For each $\Lambda \in \mathcal{L}^{(3)}$, there is a curve  $ \mathcal{V}(\Lambda) \subset \mathcal{N}$ which is isomorphic to the Fermat curve and $ \mathcal{N}_{\red}= \underset{\Lambda \in \mathcal{L}^{(3)}} \bigcup \mathcal{V}(\Lambda).$ On each such curve $ \mathcal{V}(\Lambda)$ are precisley $p^3+1$ superspecial points and each superspecial point is the intersection point of precisely $p+1$ curves $ \mathcal{V}(\Lambda^{'}) $ for some  $\Lambda^{'} \in \mathcal{L}^{(3)}$ (including  $\Lambda$).  Further, the graph whose vertices are the  superspecial points in $\mathcal{N}$ and in which two vertices are connected by an edge if and only if the corresponding superspecial points are contained in a curve  $ \mathcal{V}(\Lambda) $ for some  $\Lambda \in \mathcal{L}^{(3)}$ is a tree (since the  building in point 9 has, for $n=3$, a tree as underlying simplicial complex).   
\section{On the structure  of special cycles}
Recall that the space of special homomorphisms is defined to be the $\Qpq$-vector space $$\mathbb{V}:=\Hom_{\Zpq}(\overline{\mathbb{Y}}, \mathbb{X})\otimes_{\Z}\Q.$$
It is equipped with an hermitian form given by $$h(x,y)= \lambda^{-1}_{\overline{\Y}}\circ \hat{y}\circ \lambda_{\X}\circ x \in \End_{\Zpq}(\overline{\Y})\otimes \Q \cong \Qpq,$$
where the last isomorphism is via $\iota^{-1}$ and where $\hat{y}$ is the dual of $y.$
Also, recall that, for a special homomorphism $x\in \V$, the {\em special cycle} $\mathcal{Z}(x)$ is defined  to be the closed formal subscheme of $\mathcal{N}$ with the property that, for $S \in$ Nilp, the set  $\mathcal{Z}(x)(S)$ consists of all points $(X,\iota_X, \lambda_X, \varrho_X )$ in $\mathcal{N}(S)$ where the quasi-homomorphism
$$
\overline{\Y}\times_{\F} \overline{S}\stackrel{x}{\longrightarrow}{\X}\times_{\F} \overline{S}\stackrel{\varrho_X^{-1}}{\longrightarrow}X\times_S \overline{S}
$$
extends to an $\Zpq$-linear homomorphism $$\overline{Y}\times_W S \longrightarrow X.$$
By \cite{KR1}, Proposition 3.5,  $\mathcal{Z}(x)$ is (for $x \neq 0$) a relative divisor in $\mathcal{N}$ (or empty). We denote by $\overline{\mathbb{M}}^0=\overline{\mathbb{M}}^0_0 \oplus \overline{\mathbb{M}}^0_1=W\overline{1}_0\oplus W\overline{1}_1$ the Dieudonn\'e module of $\overline{\mathbb{Y}}.$ We denote by $\overline{N}^0=\overline{N}^0_0\oplus \overline{N}^0_1$ its isocrystal. A special homomorphism $j\in \mathbb{V}$ gives rise to a map $j_0:\ \overline{N}^0_0 \rightarrow N_0.$ By \cite{KR1}, Proposition 3.10,  a point $x\in \mathcal{N}(\F)$ with corresponding lattice $A\subset C_{W_{\Q}}=N_0$ belongs to the special cycle $\mathcal{Z}(j)$ if and  only if $j_0(\overline{1}_0)\in pA^{\vee}.$ Also, for $\Lambda \in \mathcal{L}$, the subscheme $\mathcal{V}(\Lambda)$ belongs to $\mathcal{Z}(j)$ if and only if 
$j_0(\overline{1}_0)\in p\Lambda^{\vee}.$ 
Further, by \cite{KR1}, Lemma 3.9, the map $j \mapsto j_0(\overline{1}_0)$ defines an isomorphism $\mathbb{V} \rightarrow C$ and we have $\{j_0(\overline{1}_0), j_0(\overline{1}_0)\}=p\cdot h(j,j).$ We say that $j$ is even resp. odd if the $p$-adic valuation of $h(j,j)$ is even resp. odd. By the valuation of a special homomorphism $j$ we mean the $p$-adic valuation of $h(j,j)$. We also write $\nu_p(j)$  for the valuation of $j$. We say that two special homomorphisms $j_1,j_2$ are perpendicular to each other (and write $j_1 \perp j_2$), if $h(j_1,j_2)=0.$ By a special homomorphism  of non-negative valuation we mean a special homomorphism $j$ with $h(j,j)\in \Zp \setminus \{ 0 \}.$
\newline

{\Large
 From now on we assume for the remainder of part I that  $n=3$.}

\subsection{The reduced locus of special cycles}
We introduce the following notation. Let $T \subset \mathcal{L}^{(3)}$ and let $\Lambda \in \mathcal{L}^{(3)}.$ Then we denote by $d(\Lambda, T)$ the minimal integer $d \geq 0$ such that there exist $\Lambda_0,...,\Lambda_d \in \mathcal{L}^{(3)}$ where $\Lambda_0=\Lambda$ and $\Lambda_d \in T$ and, for any $i \in \{0,...,d-1\}$, the curve $\mathcal{V}(\Lambda_i)$ intersects $\mathcal{V}(\Lambda_{i+1}).$
Analogously one defines the distance $d(\Lambda, \Lambda^{'})$ for $\Lambda, \Lambda^{'}\in  \mathcal{L}^{(3)}.$
If $\Lambda \in  \mathcal{L}^{(1)} $ and $ \Lambda^{'}\in  \mathcal{L}^{(3)}$, then we define the distance  $d(\Lambda^{'}, \Lambda)$ to be  the minimal integer $d \geq 1$ such  that there exist $\Lambda_1,...,\Lambda_d \in \mathcal{L}^{(3)}$ where $\Lambda \subset \Lambda_1$ and $\Lambda_d = \Lambda^{'}$ and, for any $i \in \{1,...,d-1\}$, the curve $\mathcal{V}(\Lambda_i)$ intersects $\mathcal{V}(\Lambda_{i+1}).$

\begin{Pro} \label{redloc}
 Let $j$ be a special homomorphism.
\begin{enumerate}
\item If the valuation of $j$ is $0$, then the underlying reduced subscheme of $\mathcal{Z}(j)$ consists of precisely one $\F$-valued point. This point is superspecial.
\item Assume that the valuation of $j$ is of the form $2r>0$ and let $A\in \mathcal{L}^{(1)}$ be the lattice corresponding to the unique superspecial point in $\mathcal{Z}(j/p^r).$ Then the underlying reduced subscheme of $\mathcal{Z}(j)$ is given by
$$\mathcal{Z}(j)_{\red}= \underset{\Lambda \in \mathcal{L}^{(3)}, \ d(\Lambda, A)\leq r }\bigcup \mathcal{V}(\Lambda).$$
\item If the valuation of $j$ is $1$, then the underlying reduced subscheme of $\mathcal{Z}(j)$ is of the form
$$\mathcal{Z}(j)_{\red}= \underset{\Lambda \in T(j)}\bigcup \mathcal{V}(\Lambda) ,$$ where $T(j)\subset \mathcal{L}^{(3)}$ has the following properties: For each $\Lambda \in T(j)$, there are precisely $p+1$ superspecial points in $\mathcal{V}(\Lambda)$ such that,  for each such superspecial point, all $p+1$ curves $\mathcal{V}(\Lambda^{'})$
 passing through it also belong to $\mathcal{Z}(j),$ i.e. the lattice $\Lambda^{'}$ also belongs to $T(j)$. 
For  all other curves passing through $\mathcal{V}(\Lambda)$, the corresponding lattices do not belong to $T(j)$. Further, for any $\Lambda,\Lambda^{'}\in T(j)$, there is there is a finite sequence $\Lambda_0, ..., \Lambda_N \in T(j)$ such that $\Lambda=\Lambda_0$ and $\Lambda^{'}=\Lambda_N$ and $\mathcal{V}(\Lambda_i)$ intersects $\mathcal{V}(\Lambda_{i+1}).$ In other words, $\mathcal{Z}(j)_{\red}$  is connected.
\item Assume that the valuation of $j$ is of the form $2r+1>0.$ 
Then the underlying reduced subscheme of $\mathcal{Z}(j)$ is given by
$$\mathcal{Z}(j)_{\red}= \underset{\Lambda \in \mathcal{L}^{(3)}, \ d(\Lambda, T(j/p^r))\leq r }\bigcup \mathcal{V}(\Lambda), $$
where 
 $T(j/p^r)$ is as in 3.
\item If the valuation of $j$ is negative, then $\mathcal{Z}(j)$ is empty.
\end{enumerate}
\end{Pro}
{\em Proof.}
1.)  We choose $j_1,j_2 \in \mathbb{V}$ of non negative valuation and linearly independent and both perpendicular to $j$. By \cite{KR1}, Lemma 5.2, for any $N \in \N$, the intersection $\mathcal{Z}(j)_{\red}\cap \mathcal{Z}(p^N j_1)_{\red} \cap \mathcal{Z}(p^N j_2)_{\red}$ consists of precisely one $\F$-valued point. But, for any subset $S\subset  \mathcal{N}(\F)$, we can choose $N$  so large that $S \subset \mathcal{Z}(p^N j_1)_{\red} \cap \mathcal{Z}(p^N j_2)_{\red}$. From this the claim of 1.) follows.

Before passing to the proof of 2.) we prove some lemmas.

\begin{Lem} \label{lem1}
Let $j\in \mathbb{V},$ let $x \in \mathcal{Z}(j)(\F)$,  and suppose that  $x \in \mathcal{V}(\Lambda)(\F)$ for some $\Lambda \in \mathcal{L}^{(3)}$. Then $\mathcal{V}(\Lambda)\subset \mathcal{Z}(p\cdot j).$\end{Lem}
{\em Proof.}
Let $M$ be the Dieudonn\'e module of the $p$-divisible group corresponding to $x.$ Then, for the lattice $A=M_0 \subset  N_0 = C_{W_{\Q}}$, we get $j_0(\overline{1}_0)\in pA^{\vee} \subset A\subset \Lambda.$ Hence $p\cdot j_0(\overline{1}_0)\in p\Lambda=p\Lambda^{\vee}.$ (Note that $\Lambda = \Lambda^{\vee}$ since  $\Lambda \in \mathcal{L}^{(3)}$.) \qed

\begin{Lem}\label{lem2h}
Let $j$ be a special homomorphism of valuation at least $2$. Suppose that $ \mathcal{V}(\Lambda) \subset Z(j)_{\red}$ for some $\Lambda \in \mathcal{L}^{(3)}$. Then either $ \mathcal{V}(\Lambda) \subset \mathcal{Z}(j/p)_{\red}$ or there is a unique point in $  \mathcal{V}(\Lambda)(\F)$ which is contained in $\mathcal{Z}(j/p)(\F)$. Further, this point is superspecial. 
\end{Lem}
{\em Proof.} Suppose that  $ \mathcal{V}(\Lambda) \not\subset Z(j/p)_{\red}.$ Let $V=\Lambda/p\Lambda^{ \vee}$ with skew hermitian from $(.,.)$ as in point 11 of the list in section 1. Then $(j/p)_0(\overline{1}_0)\notin p\Lambda^{ \vee}$ since otherwise  $ \mathcal{V}(\Lambda)\subset Z(j/p)_{\red}.$  Let $v$ be the image of $(j/p)_0(\overline{1}_0)$ in $V$. Then $v\neq 0$ but $(v,v)=0$ (since the valuation of $j$ is at least $2$). Let $w\in V$ such that $U:=v^{\perp}$ is the span of $ v$ and $w.$
Then $U$ is the unique $2$-dimensional subspace of $V$ for which $U^{\perp} \subset U$  and $v \in U.$ Let $A$ resp. $A_W$ be the corresponding lattice in $C$ resp.  $C_{W_{\Q}}.$ Since $(v,v)=0$ and $(v,w)=0$, it follows that, for any $a \in A_W$, we have $\{(j/p)_0(\overline{1}_0), a\}\in pW.$ Thus $(j/p)_0(\overline{1}_0) \in pA_W^{\vee}$ and hence the $\F$-valued point in $\mathcal{V}(\Lambda)$ corresponding to $A_W$ belongs to $\mathcal{Z}(j/p)(\F).$ The construction shows that the point is superspecial. The uniqueness follows from the uniqueness of $U$ resp. the corresponding subspace $U_{\F}$ in $V_{\F}.$  \qed
\newline

We proceed in  the proof of the proposition. First we observe that, for any special homomorphism of valuation at least $1$, the reduced locus $\mathcal{Z}(j)_{\red}$ is purely of dimension $1$. This follows from \cite{KR1}, Corollary 4.3 (reasoning analogously to the proof of point 1).

Now we prove point 2 of the proposition by induction on $r$. For $r=0$, the claim is given by point 1. Let now $r >0$ and assume the claim  for $r-1.$ Then  $\mathcal{Z}(j)_{\red}$  is a union of curves of the form  $ \mathcal{V}(\Lambda)$ for  some $\Lambda \in \mathcal{L}^{(3)}$.  Lemma \ref{lem1} and the induction hypothesis show that  $\mathcal{Z}(j)_{\red}$ contains the set $ \underset{\Lambda \in \mathcal{L}^{(3)}, \ d(\Lambda, A)\leq r }\bigcup \mathcal{V}(\Lambda).$ Let  $\mathcal{V}(\Lambda) \subset Z(j)_{\red}$ for some $\Lambda \in \mathcal{L}^{(3)}$. By Lemma \ref{lem2h},  there is at least one  point in $  \mathcal{V}(\Lambda)(\F)$ which is contained in $\mathcal{Z}(j/p)(\F).$ This shows together with the induction hypothesis that the set $ \underset{\Lambda \in \mathcal{L}^{(3)}, \ d(\Lambda, A)\leq r }\bigcup \mathcal{V}(\Lambda)$ contains $\mathcal{Z}(j)_{\red}$. This ends the proof of point 2 of the proposition.

\begin{Lem} \label{lem4}
Let $j\in \mathbb{V}$ with $h(j,j)\neq 0.$ Then $\mathcal{Z}(j)_{\red}$ is connected.
\end{Lem}
{\em Proof.} We may assume that the valuation of $j$ is at least $1$ and hence   $\mathcal{Z}(j)_{\red}$  is a union of curves of the form  $ \mathcal{V}(\Lambda)$ for  some $\Lambda \in \mathcal{L}^{(3)}$.  Assume  $\mathcal{Z}(j)_{\red}$ was not  connected. Then we could find a superspecial  point $x$ of $\mathcal{N}$ not lying in $\mathcal{Z}(j)$ but lying between two components $C$ and $C^{'}$ of $\mathcal{Z}(j)$. This means that there are pairwise distinct $\Lambda_1,..., \Lambda_N$ and pairwise distinct $\Lambda^{'}_1,..., \Lambda^{'}_{N^{'}}$ in  $\mathcal{L}^{(3)}$ with the following properties.  The point $x$ lies in the intersection  $\mathcal{V}(\Lambda_1)(\F)\cap \mathcal{V}(\Lambda^{'}_1)(\F) $  and the intersections $ \mathcal{V}(\Lambda_i)(\F) \cap \mathcal{V}(\Lambda_{i+1})(\F)$ resp. $ \mathcal{V}(\Lambda^{'}_i)(\F) \cap \mathcal{V}(\Lambda^{'}_{i+1})(\F)$ are not empty for any $i<N$ resp $i<N^{'}$. Further,  $\mathcal{V}(\Lambda_i)\cap C$ is empty for $i < N$ and the intersection $\mathcal{V}(\Lambda_N)\cap C$  consists of precisely one $\F$-valued point. Analogously,  $\mathcal{V}(\Lambda^{'}_i)\cap C^{'}$ is empty for $i < N^{'}$, and the intersection $\mathcal{V}(\Lambda^{'}_{N^{'}})\cap C^{'}$  consists of precisely one $\F$-valued point. We may assume that $N=N^{'}.$

 Now we choose a special homomorphism  $j_1$ of valuation $0$ such that $x$ is the (unique) $\F$-valued point of $\mathcal{Z}(j_1)$ and such that the fundamental matrix of $j$ and $j_1$ is non-singular. (It is easy to see that such a $j_1$ exists.) Then the intersection  $\mathcal{Z}(j)_{\red}\cap \mathcal{Z}(p^Nj_1)_{\red}$ has pure dimension (same reasoning as above). 
 Using point 2, we see that the intersection $\mathcal{Z}(j)_{\red}\cap \mathcal{Z}(p^Nj_1)_{\red}$ has dimension zero and contains at least  two distinct isolated $\F$-valued points. Now we choose a special homomorphism $j_2$ of positive valuation which is perpendicular to $j$ and $j_1$ such that $\mathcal{Z}(j_2)_{\red}$ contains these two points. Hence the intersection $\mathcal{Z}(j)_{\red} \cap \mathcal{Z}(p^N j_1)_{\red} \cap \mathcal{Z}(j_2)_{\red}$ is of dimension $0$ and contains (at least) two points. This contradicts \cite{KR1}, Corollary  4.3 and Corollary 4.7  \qed
\newline

We proceed in the proof of Proposition \ref{redloc}. 

Let us prove point 3). 
We already know  
that $\mathcal{Z}(j)_{\red}$ is of pure dimension $1$. Thus we find $\Lambda \in \mathcal{L}^{3}$ such that $\mathcal{V}(\Lambda)\subset \mathcal{Z}(j).$ We may assume that $h(j,j)=p,$ 
hence $\{j_0(\overline{1}_0),j_0(\overline{1}_0)\}=p^2.$ Then $j_0(\overline{1}_0) \in p\Lambda^{\vee}=p\Lambda.$ Writing $j_0(\overline{1}_0)=p\cdot l$ it follows that $l \in \Lambda.$ 
If now $\mathcal{V}(\Lambda^{'})$ intersects 
$\mathcal{V}(\Lambda)$ and is also contained in $\mathcal{Z}(j)$, then it follows that $l \in \Lambda \cap \Lambda^{'}=:A$. If $\mathcal{V}(\Lambda^{''})$ is a third curve also passing through the intersection point of $\mathcal{V}(\Lambda)$ and $\mathcal{V}(\Lambda^{'})$ (which is the $\F$-valued point associated to $A$), then $\Lambda \cap \Lambda^{''}=A,$ hence $pl \in p\Lambda^{'' \vee},$ hence $\mathcal{V}(\Lambda^{''})\subset \mathcal{Z}(j).$
Thus, if two curves $\mathcal{V}(\Lambda), \mathcal{V}(\Lambda^{'}) \subset \mathcal{Z}(j)$ intersect in a point, then all curves $\mathcal{V}(\Lambda^{''})$ passing through this point belong to $\mathcal{Z}(j).$

Next we want to count those superspecial points on $\mathcal{V}(\Lambda)$ for which there is another curve $\mathcal{V}(\Lambda^{'})\subset \mathcal{Z}(j)$ passing through that point. In view of the above reasoning this means we count the superspecial points on $\mathcal{V}(\Lambda)$  such that, for the corresponding lattice $A$, we have $l \in A.$ Let $\bar{l}$ denote the image of $l$ in $V=\Lambda/p\Lambda.$ Since $\{l,l\}=1$ we have  $\bar{l}\neq 0.$
In view of the last point  of the list in section 1 this means we count the subspaces $U$ in $V$ of dimension $2$ which contain $l$ and for which $U^{\perp}\subset U.$ Denoting by $W$ the orthogonal complement of the span of $l$  this means that we count the number of subspaces $U$ in $W$ of dimension $1$ with $U^{\perp}\subset U.$ By \cite{VW}, Example 4.6, this number is $p+1.$

Since we already know from Lemma \ref{lem4} that $\mathcal{Z}(j)$ is connected, it follows that $\mathcal{Z}(j)$ has the desired description.

The claim of 4.) follows from 3.) by induction on $r$ and with the help of Lemma \ref{lem1} and Lemma \ref{lem2h} in the same way as point 2 follows from point 1.

Point 5 is clear (comp. also \cite{KR1}, Corollary 3.11).
\qed

\begin{Pro}\label{redschnitt}
Let $j_1,j_2$ be perpendicular special homomorphisms of nonnegative valuations $a_1 \leq a_2.$
\begin{enumerate}
\item The intersection $\mathcal{Z}(j_1)\cap \mathcal{Z}(j_2)$ is non-empty and connected.
\item  If  $a_1=0$, then this intersection has precisely one $\F$-valued point. It is superspecial.
\item If  $a_1=a_2=1$, then  $(\mathcal{Z}(j_1)\cap \mathcal{Z}(j_2))_{\red}$ consists of precisley one curve $\mathcal{V}(\Lambda)$ for some $\Lambda \in \mathcal{L}^{(3)}$.
\item If $j_3$ is a third special homomorphism which is perpendicular to $j_1$ and to $j_2$ and has valuation $a_3\geq a_1,a_2$, then   $(\mathcal{Z}(j_1)\cap \mathcal{Z}(j_2)\cap \mathcal{Z}(j_3))_{\red}= (\mathcal{Z}(j_1)\cap \mathcal{Z}(j_2))_{\red}.$
\end{enumerate}
\end{Pro}
{\em Proof.} 1.) To show that $\mathcal{Z}(j_1)\cap \mathcal{Z}(j_2)$ is non-empty, we may restrict ourselves to the case that $a_1,a_2 \in \{0,1\}.$
In the case $a_1=0$ it follows from \cite{KR1}, Lemma 5.2  that  $\mathcal{Z}(j_1)\cap \mathcal{Z}(j_2)$ is non-empty. Suppose now that $a_1=a_2=1.$ We write $j_i=j_{i,0}\oplus j_{i,1}$ for the $0$ and $1$ component of $j_i$. Let $l_i=p^{-1}j_{i,0}(\overline{1}_0) \in C.$
Then $\{l_0,l_1\}=0$ and $\{l_0,l_0\}$ and $\{l_1,l_1\}$ both have valuation $0.$ Thus we see that there is precisely one lattice $\Lambda \in  \mathcal{L}^{(3)}$ with $l_0,l_1 \in  \mathcal{L}^{(3)}.$ From this the claim follows. The fact that $\mathcal{Z}(j_1) \cap \mathcal{Z}(j_2)$ is connected can be proved like Lemma \ref{lem4}

2.) follows from \cite{KR1}, Lemma 5.2.

3.) Follows from 1.) and its proof.

4.) For $i=1,2,3$, we write $j_i=p^{r_i}y_i$ where $y_i$ has valuation $0$ or $1.$ One shows as in the preceding points that $(\mathcal{Z}(y_1)\cap \mathcal{Z}(y_2) \cap \mathcal{Z}(y_3))_{\red}$ is not empty and that it consists of one single $\F$-valued point (and this point is superspecial) if one of the $j_i$ is even. If all $j_i$ are odd, one shows as above that $(\mathcal{Z}(y_1)\cap \mathcal{Z}(y_2) \cap \mathcal{Z}(y_3))_{\red}$  consists of precisley one curve $\mathcal{V}(\Lambda)$ for some $\Lambda \in \mathcal{L}^{(3)}$. Now one uses the description of the $\mathcal{Z}(j_i)$ given in Proposition \ref{redloc} to conclude the claim. \qed
\newline 

Finally, we note the following (easy to prove) statement.
\begin{Lem}\label{lem5}
Suppose we are given two perpendicular special homomorphisms $j_0$ and $j_1$, and suppose that the valuation of $j_0$ is $0$ and that the valuation of $j_1$ is $1$. Let $x$ be the (unique) $\F$-valued point of $\mathcal{Z}(j_0)$. Then $\mathcal{Z}(j_1)$ contains all $p+1$ curves $\mathcal{V}(\Lambda_1),...,\mathcal{V}(\Lambda_{p+1})$ passing through $x$. \qed
\end{Lem}

\subsection{The special fiber of special cycles}
For any (formal) subscheme $X$ of $\mathcal{N}$, we denote by $X_p$ its special fiber.
For  $j\in \mathbb{V}$ of non-negative valuation, we denote 
by  $S(j)\subset \mathcal{L}^{(3)}$ the subset of such $\Lambda$ for which $\mathcal{V}(\Lambda) \subset \mathcal{Z}(j)$, see Proposition \ref{redloc}.
Further, for $\Lambda \in S(j)$, we denote by $b_j(\Lambda)$ the maximal integer $a$ such that $\Lambda \in S(j/p^a),$ see also Proposition \ref{redloc}.

\begin{The} \label{spezf}
Let $j\in \mathbb{V}$ be of valuation $r \geq 0$. The special  fiber of $\mathcal{Z}(j)$ as a divisor in $\mathcal{N}_p$ can be written as follows,
$$
\mathcal{Z}(j)_p= p^{r}\cdot t + \sum_{\Lambda \in S(j)}(1+p^2+...+ p^{2b_j(\Lambda)})\cdot \mathcal{V}(\Lambda). 
$$
Here $t$ is zero for $r$. For $r$ even, it is an irreducible divisor in $\mathcal{N}_p$ which passes through the unique supersingular point of $\mathcal{Z}(j/p^{r/2}).$ Its intersection multiplicity with each $\mathcal{V}(\Lambda)$ which contains this point is $1$.
\end{The}
{\em Proof.}
Let  $\mathcal{V}(\Lambda) \subset \mathcal{Z}(j)$ 
 Let $x\in \mathcal{V}(\Lambda)(\F)$ be a non superspecial point. We want to determine the equation of $\mathcal{Z}(j)$ in $\mathcal{O}_{\mathcal{N}_p,x}$ and in particular the  multiplicity of  $\mathcal{V}(\Lambda)$ in $\mathcal{Z}(j)_p.$ Let $X$ be the corresponding $p$-divisible group and let $M$ be its Dieudonn\'e module. Let $M_p$ be the reduction mod $p$ of $M.$ We denote the  operators on $M_p$ induced by $V$ and $F$ also be these letters. Since $M$ is endowed with a $\Zpq$-action, we get a grading $M=M_0 \oplus M_1$ and likewise $M_p=M_{p,0}\oplus M_{p,1}.$ By \cite{VW}, Theorem 2.1, the reduction $M_p$ is isomorphic to $\mathbb{B}(3).$ This means that we find  bases $\bar{e_1},\bar{e_2},\bar{e_3}$ of $M_{p,0}$ and  $\bar{f_1},\bar{f_2},\bar{f_3}$ of $M_{p,1}$ such that $V( \bar{f_1})=\bar{e_2},$ 
 $V( \bar{f_2})=\bar{e_3},$   $V( \bar{e_3})=\bar{f_1}$ and $F( \bar{e_1})=\bar{f_3},$ 
 $F (\bar{f_2})=\bar{e_1},$   $F( \bar{f_3})=-\bar{e_2}$, and for the induced alternating form we have $\langle \bar{e_i}, \bar{f_j} \rangle= \pm\delta_{ij}$. 

 Let $f_2\in M_1$ be a lift of $\bar{f_2}.$ We define $e_3=V(f_2), f_1=V^2(f_2), e_2=V^3(f_2), e_1=F(f_2), f_3=F^2(f_2).$ Then the $e_i$ are lifts of the $\bar{e_i}$ and form a basis of $M_0$ and  the $f_i$ are lifts of the $\bar{f_i}$ and form a basis of $M_1$. Further, there are $x,y,z \in W$ such that $F(f_3)= px\cdot e_1 + (-1+py)\cdot e_2 + pz \cdot e_3$.

 Denote by  $T$ the $W$-span of $e_1,f_2,f_3$ and by $L$ the $W$-span of $e_2,e_3,f_1$. Then
\[
M=L\oplus T, \ \ \ VM=L\oplus pT.
\]
Let $h_1=e_1,h_2=f_2,h_3=f_3,h_4=e_2,h_5=e_3,h_6=f_1.$ 
Define the matrix $(\alpha_{ij})$ by 
\[
Fh_j=\sum_i \alpha_{ij}h_i \text{ for } j=1,2,3,  \\
 \]
 \[
 V^{-1}h_j=\sum_i \alpha_{ij}h_i \text{ for } j=4,5,6.
\]
Hence 
\[
(\alpha_{ij})=
\begin{pmatrix}
&1&px&& \\
&&&&1& \\
1&&&&&& \\
&&-1+py&& \\
&&pz&&&1\\
&&&1&&\\
\end{pmatrix}.
\]
It follows (see \cite{Z1}, p. 48) that the universal deformation of $X$ over $\F[\![ t_1,...,t_9]\!]$ corresponds to the display $(L\oplus T)\otimes W(\F[\![t_1,...,t_9 ]\!]) $ with matrix $(\alpha_{ij})^{\univ}$ (wrt. the basis $h_1,...h_6$ and with entries in $W(\F[\![t_1,...,t_9 ]\!])$ given by 
\[
(\alpha_{ij})^{\univ}=
\begin{pmatrix}
1&&&[t_1]&[t_2]&[t_3] \\
&1&&[t_4]&[t_5]&[t_6] \\
&&1&[t_7]&[t_8]&[t_9] \\
&&&1&& \\
&&&&1& \\
&&&&&1 \\
\end{pmatrix} \cdot
(\alpha_{ij})
=
\begin{pmatrix}
&1&px+[t_1](-1+py)+[t_2]pz&[t_3]&&[t_2]\\
&&[t_4](-1+py)+[t_5]pz&[t_6]&1&[t_5]\\
1&&[t_7](-1+py)+[t_8]pz&[t_9]&&[t_8]\\
&&-1+py&&&\\
&&pz&&&1\\
&&&1&&\\
\end{pmatrix}.
\]
Here the $[t_i]$ denote the Teichm\"uller representatives of the $t_i$.

 Now let $A^{'}=W[\![t_1,...,t_9]\!]$,  let  $R^{'}=\F[\![t_1,...,t_9]\!]$. 

   We extend the Frobenius  $\sigma $ on $W$ to $A^{'}$ 
  putting $\sigma(t_i)=t_i^p.$ 
 Let $R$ be the completed universal deformation ring of $X$ together with the $\Zpq$-action. Then $R$ is a quotient of $R^{'}$ by an ideal $J$. 
We observe that $t_3,t_4,t_5,t_7,t_8\in J.$ This follows from the fact that the corresponding map  $\alpha_{ij}$ of the corresponding display is homogenous of degree $1$. On the other hand, we see similarly that indeed $J=(t_3,t_4,t_5,t_7,t_8).$ 
Denoting the image of $t_i$ in $R$ for $i=1,2,6,9$ also by $t_i$, we write $R=\F[\![t_1,t_2,t_6,t_9]\!]$. We also define $A=W[\![t_1,t_2,t_6,t_9]\!].$
 Further, for any $n \in \N$,  denote by $\mathfrak{a}_n$ resp. $\mathfrak{r}_n$ the ideal in $A$ resp. in $R$ generated by the monomials $t_1^{a_1} t_2^{a_2} t_6^{a_6} t_9^{a_9}$, where  $a_i\geq 0$ and $\sum a_i = n$. Let $A_n=A/\mathfrak{a}_n$ and $R_n=R/\mathfrak{r}_n$.  Then $A^{'}$ is a frame for $R^{'}$ resp. $A$ is a frame for $R$ resp. $A_n$ is a frame for $R_n$.

For an  $A^{'}$ - $R^{'}$-window $(M^{'}, M^{'}_1, \Phi^{'},\Phi^{'}_1)$, let  $M_1^{'^{\sigma}}=A^{'}\otimes_{A^{'}, \sigma}M_1^{'}$ and denote by $\Psi^{'}:  M_1^{'^{\sigma}} \rightarrow M^{'}$ the linearization of $\Phi^{'}_1$. It is an isomorphism of  $A^{'}$-modules.
 Denote by  $\alpha^{'}:M_1^{'} \rightarrow M_1^{'^{\sigma}}$ the composition of the inclusion map $M_1^{'}  \hookrightarrow M^{'} $ followed by $\Psi^{'^{-1}}$. In this way, the category of pairs $(M_1^{'}, \alpha^{'})$ consisting of a free $A^{'}$-module of finite rank and an $A^{'}$-linear injective homomorphism $\alpha^{'}:M_1^{'} \rightarrow M_1^{'^{\sigma}}$ such that Coker  $\alpha^{'}$ is a free $R^{'}$-module becomes equivalent to the category of formal $p$-divisible groups over $R^{'}.$  (Note that the so called nilpotence condition is fulfilled automatically here, comp. also \cite{KR1}, section 8.) A corresponding description holds for the category of  formal $p$-divisible groups over $R$ resp. $R_n$.

 We consider the $A^{'}$ - $R^{'}$ window $(M^{'},M^{'}_1, \Phi^{'},\Phi^{'}_1)$ given by 
\[ 
M^{'}=M\otimes A^{'},  \ M_1^{'}=VM\otimes A^{'}, \ 
\Phi^{'}=
\begin{pmatrix}
&1&px+[t_1](-1+py)+[t_2]pz&p[t_3]&&p[t_2]\\
&&[t_4](-1+py)+[t_5]pz&p[t_6]&p&p[t_5]\\
1&&[t_7](-1+py)+[t_8]pz&p[t_9]&&p[t_8]\\
&&-1+py&&&\\
&&pz&&&p\\
&&&p&&\\
\end{pmatrix} \sigma,  \   
\] and $\Phi^{'}_1=\frac{1}{p}\cdot \Phi^{'}$, where 
the matrix of $\Phi^{'}$ is described in the basis $h_1,...,h_6$. The corresponding display is the universal display described above (easy to see using the procedure described on p.2 of \cite{Z2}). Hence $(M^{'},M^{'}_1, \Phi^{'},\Phi^{'}_1)$ is the universal window.

The corresponding matrix of $\alpha^{'}$ wrt. the basis $ph_1,ph_2,ph_3,h_4,h_5,h_6$ of $ M_1^{'}$ resp. the basis $p(1\otimes h_1),p(1\otimes h_2),p(1\otimes h_3),1\otimes h_4,1\otimes h_5,1\otimes h_6$ of $ M_1^{'^{\sigma}}$ is given by 
\[\alpha^{'}=
\begin{pmatrix}
&&p&-t_7&-t_8&-t_9 \\
p&&&-\frac{px}{-1+py}-t_1&-t_2&-t_3 \\
&&&\frac{1}{-1+py}& \\
&&&&&1\\
&p&&-t_4&-t_5&-t_6\\
&&&-\frac{pz}{-1+py}&1
\end{pmatrix}.
\]

We conclude that the map  $\alpha: M_1 \rightarrow M_1^{\sigma}$ corresponding to the $A$-$R$ window of the universal defomation of $X, \iota$  (which is the base change from $(M_1^{'}, \alpha^{'})$) can be written as follows (using the bases $pe_1,e_2,e_3, f_1,pf_2,pf_3$ resp. $p(1\otimes e_1),1\otimes e_2,1\otimes e_3,1\otimes  f_1,p(1\otimes f_2),p(1\otimes f_3)$ ),

\[\alpha=
\begin{pmatrix}
&&&-t_9&&p \\
&&&1&& \\
&&&-t_6&p& \\
&-\frac{pz}{-1+py}&1&& \\
p&-\frac{px}{-1+py}-t_1&-t_2&\\
&\frac{1}{-1+py}\\
\end{pmatrix}=\begin{pmatrix}
&U\\
\tilde{U}&\\
 \end{pmatrix},
\]
where
\[U=
\begin{pmatrix}
-t_9&&p \\
1&& \\
-t_6&p& \\

\end{pmatrix}
\ \ \ \text{and } \ \ \
\tilde{U}=
\begin{pmatrix}
&-\frac{pz}{-1+py}&1\\
p&-\frac{px}{-1+py}-t_1&-t_2\\
&\frac{1}{-1+py}\\
\end{pmatrix}.
\]

The corresponding universal $p$-divisible groups over 
  $R_n$ correspond to the pairs 
  $(M_1(n), \alpha(n))$ obtained by base change from $(M_1, \alpha)$.

Consider the $p$-divisible group $\overline{\mathbb{Y}}$ with its Dieudonn\'e module  $\overline{\mathbb{M}}^0=W\overline{1}_0\oplus W\overline{1}_1.$ Let $\overline{M}^0_{,1}=W\overline{1}_0\oplus p\cdot W \overline{1}_1$ and let $n_0=\overline{1}_0$ and  $n_1=p\overline{1}_1$. Then  $\overline{\mathbb{Y}}$ corresponds to the pair $(\overline{M}^0_{,1}, \beta)$ where $\beta(n_0)=1\otimes n_1$ and $\beta(n_1)=-p\otimes n_0$. By base change $W \rightarrow A_n$ resp.  $W \rightarrow A$ we obtain pairs $({\overline{M}_{,1}^0}_{R}, \beta)$ resp. $({\overline{M}_{,1}^0}_{R_n}, \beta)$ corresponding   to the constant $p$-divisible group  $\overline{\mathbb{Y}}$ over $R$ resp. $R_n$. We denote the matrix of $\beta$  by $S$, hence
$$S=
\begin{pmatrix}
0&-p\\
1&0\\
\end{pmatrix}.
$$

Next we want to investigate the ideal in $R$ describing the maximal deformation of the homomorphism $j$. To this end we determine the image of this ideal in the rings $R_n.$

The map $j$ corresponds to a map $ j(1):\overline{M}_{,1}^0 \rightarrow  M_1(1) $ such that the following diagram commutes
\[
\xymatrix{ \overline{M}_{,1}^0 \ar[d]_{j(1)} \ar[r]^{\beta}&  {\overline{M}_{,1}^0}^{\sigma} 
\ar[d]^{\sigma(j(1))}  \\M_1(1)\ar[r]_{\alpha(1)} & M_1(1)^{\sigma}  .}
\]
Then $j$ lifts over $R_n$ if and only if there is a lift $j(n)$ of $j(1)$ such that 

\[
\xymatrix{ {\overline{M}_{,1}^0}_{R_n} \ar[d]_{j(n)} \ar[r]^{\beta}&  {{\overline{M}_{,1}^0}_{R_n}}^{\sigma} 
\ar[d]^{\sigma(j(n))}  \\M_1(n)\ar[r]_{\alpha(n)} & M_1(n)^{\sigma}  .}
\]

We write $j=j_0 \oplus j_1$ where $j_i:\overline{M}^0_i \rightarrow M_i.$ We write $j_0(\overline{1}_0)=a \cdot e_1+b \cdot e_2+c\cdot e_3.$ We also write $j(1)=(X(1), Y(1)).$ Then $X(1)$ can be written in the above basis as 
$$ X(1)=
\begin{pmatrix}
a/p&0 \\
b&0\\
c&0\\
\end{pmatrix}.
$$
Using $Fj=jF$ and hence $j(\overline{1}_1)=FjF^{-1}(\overline{1}_1)$  one sees that 
$$ Y(1)=
\begin{pmatrix}
0&-pb^{\sigma} \\
0&-c^{\sigma}\\
0&-a^{\sigma}/p\\
\end{pmatrix}.$$
Using  $j(\overline{1}_1)=FjF^{-1}(\overline{1}_1)$ one easily checks that for the $p$-adic valuations we have $\nu_p(a)=\nu_p(b)+1$ and that $c$ has at least valuation $\nu_p(b)+1.$

We are looking for liftings $X(n)$ of $X(1)$ and  $Y(n)$ of $Y(1)$ over $A_n$ such that
 \[ UY(n)= \sigma (X(n))S \  \text{ and } \   \tilde{U}X(n)= \sigma (Y(n))S. \tag{$*$}\]
Suppose $n=p^l$, where $l\geq1$, and suppose we have found liftings  $X(p^{l-1})$ and  $Y(p^{l-1})$  satifying $(*)$. For any choice of liftings   $X(p^{l})$ and  $Y(p^{l})$  of $X(p^{l-1})$ and  $Y(p^{l-1})$, 
the matrices $\sigma(X(p^{l}))$ resp. $\sigma(Y(p^{l}))$ are equal  to $\sigma(X(p^{l-1})$ resp. $\sigma(Y(p^{l-1}))$ interpreted as matrix over $A_{p^l}$. Hence there are liftings  $X(p^{l})$ and  $Y(p^{l})$  satifying (*) if and only if 
\[ U^{-1}\sigma (X(p^{l-1}))S \  \text{ and } \  \tilde{U}^{-1} \sigma (Y(p^{l-1}))S\] are integral, and in this case 
\[ X(p^l)= \tilde{U}^{-1} \sigma (Y(p^{l-1}))S  \  \text{ and } \ Y(p^l)= U^{-1}\sigma (X(p^{l-1}))S.\] 
Define now inductively matrices $X_{\Q}(p^l)$ and $Y_{\Q}(p^l)$ over $A_{p^l}\otimes_{\Z}\Q$ as follows: $X_{\Q}(1)=X(1)$ and $Y_{\Q}(1)=Y(1)$ and 
\[ X_{\Q}(p^{l+1})=   \tilde{U}^{-1} \sigma (Y(p^{l}))S  \  \text{ and } \  Y_{\Q}(p^{l+1})=U^{-1}\sigma (X(p^{l}))S.\] (Again $\sigma (X_{\Q}(p^l))$ and $\sigma (Y_{\Q}(p^l))$ are well defined over $A_{p^{l+1}}\otimes_{\Z}\Q$.)

Using the above form of $X(1)$ and $Y(1)$, the following lemma can easily be proved by induction.
\begin{Lem}
For any integer $k\geq 0$, the matrix 
$X(p^{2k+1})$ is of the form 
$$ X(p^{2k+1})= 
\begin{pmatrix}
t_1^{1+p^2+...+p^{2k}}\cdot p^{\nu(b)-1-k} \cdot \varepsilon&0\\
0&0\\
0&0\\
\end{pmatrix}
+A(p^{2k+1}), 
$$
where $\varepsilon$ is a unit in $W$ and $A(p^{2k+1})$ is a matrix with entries in $p^{\nu(b)-k} A_{p^{2k+1}}.$
 The matrix 
$Y(p^{2k+1})$ has entries in $p^{\nu(b)-k} A_{p^{2k+1}}.$
For any integer $k\geq 0$, the matrices
$X(p^{2k})$ and $Y(p^{2k})$ both have entries in $p^{\nu(b)-k} A_{p^{2k}}.$
\end{Lem}
From this lemma it follows that that the ideal in $R$ describing the deformation of $j$ is $(t_1^{1+p^2+...+p^{2\nu(b)}}).$

The (completed) local ring  of $\mathcal{N}_p$ in $x$ is a quotient $\overline{R}$ of $R$  and the equation for $\mathcal{Z}(j)$ in $x$ is $\overline{t_1}^{1+p^2+...+p^{2\nu(b)}}=0,$ where $\overline{t_1}$ is the image of $t_1$ in $\overline{R}.$
\newline

Next we show that $\overline{t_1}=0$ is the equation of $\mathcal{V}(\Lambda)$ in  (the completed  local ring of) $x.$
First we observe that, for any special homomorphism $y$ such that $x \in \mathcal{Z}(y)(\F)$, 
 the equation of $\mathcal{Z}(y)_p$ is locally around $x$ given by $\overline{t_1}^q=0$ for some $q.$ Now we choose three linearly independent special homomorphisms $y_1,y_2,y_3$ which come from homomorphisms between the corresponding abelian varieties of $\overline{\mathbb{Y}}$ and $X$ such that  $x \in (\mathcal{Z}(y_1) \cap \mathcal{Z}(y_2) \cap \mathcal{Z}(y_3))(\F).$  Then it follows from  \cite{KR2}, Lemma 2.21 that  the intersection $\mathcal{Z}(y_1) \cap \mathcal{Z}(y_2) \cap \mathcal{Z}(y_3)$ has support in the supersingular locus. But locally around $x$ the locus defined by $\overline{t_1}=0$ in $\mathcal{N}_p$ is contained in  $\mathcal{Z}(y_1) \cap \mathcal{Z}(y_2) \cap \mathcal{Z}(y_3).$ Thus (locally around $x$) the locus defined by $\overline{t_1}=0$ also has support in the supersingular locus, i.e. in $\mathcal{V}(\Lambda)$.
 Thus locally around $x$ the divisor given by  $\overline{t_1}=0$ in $\mathcal{N}_p$ is a multiple of the divisor $\mathcal{V}(\Lambda).$ It remains to show that this multiplicity is $1$. 
 We choose a special homomorphism $y$ for which $\mathcal{V}(\Lambda)\subset \mathcal{Z}(y)$ and  $b_y(\Lambda)=0$, hence the equation of $\mathcal{Z}(y)_p$ in (the completed  local ring of) $x$ is given by $\overline{t_1}=0.$ 
 We choose a superspecial point $z$ on $\mathcal{V}(\Lambda)$ such that  $\mathcal{V}(\Lambda)$ is the only curve passing through $z$ which belongs to $\mathcal{Z}(y).$ We choose three special homomorphisms $y_1,y_2,y_3$ which are perpendicular to each other and such that $z\in \mathcal{Z}(y_i)(\F)$ for all $i$ and such that the valuations of $y_1$ and $y_2$ are both $0$ and the valuation of $y_3$ is $1$. We write $y=\alpha y_1+ \beta y_2 + \gamma y_3$ for some $\alpha, \beta, \gamma \in \Zpq.$ It follows that $\alpha$ and $\beta$ are units. Let $\tilde{y}= \beta y_2 + \gamma y_3$. Then $z \in \mathcal{Z}(\tilde{y})(\F)$ and $\tilde{y}$ has valuation $0$ and $\tilde{y}$ is perpendicular to $y_1$. Using that $\mathcal{Z}(y)_p\cap \mathcal{Z}(\tilde{y})_p=\mathcal{Z}(y_1)_p\cap \mathcal{Z}(\tilde{y})_p$ it follows from \cite{KR1}, Proposition 8.2 that the length of the ring $\mathcal{O}_{\mathcal{Z}(y)_p \cap \mathcal{Z}(\tilde{y})_p, z}$ is $1.$ From this it follows that the equation $\overline{t_1}=0$ defines indeed a reduced divisor (locally around $x$).
 
 Thus we have shown that,  for any $x \in  \mathcal{V}(\Lambda)(\F)$ which is not superspecial,  the divisor $\mathcal{Z}(j)$ is locally around $x$ indeed equal to $(1+p^2+...+ p^{2b_j(\Lambda)})\cdot \mathcal{V}(\Lambda)$.
\newline

Suppose now that $x$ is a superspecial point in $\mathcal{Z}(j)(\F)$ and that, in the case that $j$ is even, $x$ is not the superspecial point of $\mathcal{Z}(j/p^{r/2}).$
We also suppose that in the case that $r$ is odd $x$ is not the intersection point of two curves $\mathcal{V}(\Lambda)$ and $\mathcal{V}(\Lambda^{'})$ which belong to $\mathcal{Z}(j/p^{(r-1)/2}).$
 Again we choose three special homomorphisms $y_1,y_2,y_3$ which are perpendicular to each other and such that $x\in \mathcal{Z}(y_i)(\F)$ for all $i$ and such that the valuations of $y_1$ and $y_2$ are both $0$ and the valuation of $y_3$ is $1$. We write $j=\alpha y_1+ \beta y_2 + \gamma y_3$ for some $\alpha, \beta, \gamma \in \Zpq.$ Suppose that $x$ belongs to $\mathcal{Z}(j/p^a)(\F)$ but not to $\mathcal{Z}(j/p^{a+1})(\F).$  Then it follows that there is a curve  $\mathcal{V}(\Lambda)$ with $x \in  \mathcal{V}(\Lambda)(\F) \subset \mathcal{Z}(j)(\F)$ with $b_j(\Lambda)=a$. In the case $a>0$ there are $p$ other curves  $\mathcal{V}(\Lambda^{'})$ with $x \in  \mathcal{V}(\Lambda^{'})(\F) \subset \mathcal{Z}(j)(\F)$ with $b_j(\Lambda^{'})=a-1.$ In the case $a=0$ there is no other curve  $\mathcal{V}(\Lambda^{'})$ which contains $x$ and belongs to $\mathcal{Z}(j).$
It follows that $\nu(\alpha)=\nu(\beta)=a$ and $\nu(\gamma)\geq a.$  Now it follows from \cite{KR1}, Proposition 8.2 that the length of the ring $\mathcal{O}_{\mathcal{Z}(j)_p \cap \mathcal{Z}(y_1)_p, x}$ is $1+p+...+p^{2a}.$ We already know that $\mathcal{V}(\Lambda)$ has multiplicity $1+p^2+...+p^{2a}$ in $\mathcal{Z}(j)_p$. Each of the $p$ other curves  $\mathcal{V}(\Lambda^{'})$ has (for $a>0$) multiplicity  $1+p^2+...+p^{2a-2}.$ Thus the intersection multiplicity of $\mathcal{Z}(y_1)_p$ with the supersingular curves  of $\mathcal{Z}(j)_p$ (counted with multiplicities) is at least $1+p^2+...+p^{2a} + p(1+p^2+...+p^{2a-2})=1+p+p^2+...+p^{2a}$ which is already the length of the ring $\mathcal{O}_{\mathcal{Z}(j)_p \cap \mathcal{Z}(y_1)_p, x}.$ Thus locally around $x$ the divisor $\mathcal{Z}(j)_p$ in $\mathcal{N}_p$ consists only of the supersingular curves passing through $x$ with their multiplicities, i.e. there are no non-supersingular components of $\mathcal{Z}(j)_p$ passing through $x.$

Next suppose that $j$ is odd and that  $x$ is the intersection point of two curves $\mathcal{V}(\Lambda)$ and $\mathcal{V}(\Lambda^{'})$ which belong to $\mathcal{Z}(j/p^{(r-1)/2}).$ There is a special homomorphism $y$ of valuation $0$ such that $y\perp j$ and $x\in \mathcal{Z}(y)(\F).$  It follows from \cite{KR1}, Proposition 8.2 that the length of the ring $\mathcal{O}_{\mathcal{Z}(j)_p \cap \mathcal{Z}(y)_p, x}$ is $1+p+...+p^{r}.$ We know that each of the $p+1$ supersingular curves passing through $x$ are contained in $\mathcal{Z}(j)_p$ with multiplicity $1+p^2+...+p^{r-1}.$ Thus the intersection multiplicity of $\mathcal{Z}(y)_p$ with the supersingular curves  of $\mathcal{Z}(j)_p$ (counted with multiplicities) is at least $(p+1)(1+p^2+...+p^{r-1})=1+p+p^2+...+p^{r}$ which is already the length of the ring $\mathcal{O}_{\mathcal{Z}(j)_p \cap \mathcal{Z}(y)_p, x}.$ Thus locally around $x$ the divisor $\mathcal{Z}(j)_p$ in $\mathcal{N}_p$ consists only of the supersingular curves passing through $x$ with their multiplicities, i.e. there are no non supersingular components of $\mathcal{Z}(j)_p$ passing through $x.$

Now we assume that $j$ is even and that $x$ is the superspecial point of $\mathcal{Z}(j/p^{r/2}).$  Let $X$ be the corresponding $p$-divisible group and let $M=M_0\oplus M_1$ be its Dieudonn\'e module. We find bases $e_1,e_2,e_3$ of $M_0$ and $f_1,f_2,f_3$  of $M_1$ such that for the alternating form coming from the $p$-principal polarization we have $\langle e_i, f_j \rangle= \varepsilon_i \delta_{ij}$, where $ \varepsilon_1= \varepsilon_2=-1$ and $ \varepsilon_3=1$ and such that the operators $F$ and $V$ are given as follows:
$-V(f_1)=F(f_1)=e_1,\ -V(f_2)=F(f_2)=e_2$ and $-V(e_3)=F(e_3)=f_3.$ Let $L$ be the span of $e_1,e_2,f_3$ and let $T$ be the span of $f_1,f_2,e_3.$
Let $h_1=f_1,h_2=f_2,h_3=e_3,h_4=e_1,h_5=e_2,h_6=f_3.$ 
Define the matrix $(\alpha_{ij})$ by 
\[
Fh_j=\sum_i \alpha_{ij}h_i \text{ for } j=1,2,3,  \\
 \]
 \[
 V^{-1}h_j=\sum_i \alpha_{ij}h_i \text{ for } j=4,5,6.
\]
Hence 
\[
(\alpha_{ij})=
\begin{pmatrix}
&&&-1&& \\
&&&&-1& \\
&&&&&-1& \\
1&&&& \\
&1&&&&\\
&&1&&&\\
\end{pmatrix}.
\]
It follows (see \cite{Z1}, p. 48) that the universal deformation of $X$ over $\F[\![ t_1,...,t_9]\!]$ corresponds to the display $(L\oplus T)\otimes W(\F[\![t_1,...,t_9 ]\!]) $ with matrix $(\alpha_{ij})^{\univ}$ (wrt. the basis $h_1,...h_6$ and with entries in $W(\F[\![t_1,...,t_9 ]\!])$ given by 
\[
(\alpha_{ij})^{\univ}=
\begin{pmatrix}
1&&&[t_1]&[t_2]&[t_3] \\
&1&&[t_4]&[t_5]&[t_6] \\
&&1&[t_7]&[t_8]&[t_9] \\
&&&1&& \\
&&&&1& \\
&&&&&1 \\
\end{pmatrix} \cdot
(\alpha_{ij})
=
\begin{pmatrix}
[t_1]&[t_2]&[t_3]&-1 \\
[t_4]&[t_5]&[t_6] &&-1\\
[t_7]&[t_8]&[t_9]&&&-1\\
1&&&& \\
&1&&&&\\
&&1&&&\\
\end{pmatrix}.
\]
Here the $[t_i]$ denote again the Teichm\"uller representatives of the $t_i$.

Next it is easy to see that the ideal in  $R^{'}$ describing the deformation of $\iota$ and $\lambda$ equals $(t_1,t_2,t_4,t_5,t_9,t_3-t_7,t_6-t_8).$
Let $R=\F[\![t_1,...,t_9]\!]/(t_1,t_2,t_4,t_5,t_9,t_3-t_7,t_6-t_8)$ and let $A=W[\![t_1,...,t_9]\!]/(t_1,t_2,t_4,t_5,t_9,t_3-t_7,t_6-t_8).$ Again let us 
 for any $n \in \N$ denote by $\mathfrak{a}_n$ resp. $\mathfrak{r}_n$ the ideal in $A$ resp. in $R$ generated by the monomials $t_3^{a_3} t_6^{a_6}$ where  $a_i\geq 0$ and $\sum a_i = n$. Let $A_n=A/\mathfrak{a}_n$ and $R_n=R/\mathfrak{r}_n$.  Then $A$ is a frame for $R$ resp. $A_n$ is a frame for $R_n$.
As above we construct the corresponding window of the above universal display and the corresponding pair $(M_{,1}, \alpha).$ Here $M_{,1}$ is obtained by tensoring the $W$- submodule of $M$ spanned by $pf_1,pf_2,pe_3,e_1,e_2,f_3$ with $A$. The matrix of $\alpha$ is given by 
\[\alpha=
\begin{pmatrix}
&&&-p&&t_3 \\
&&&&-p&t_6 \\
&&&&&1 \\
1&&&& \\
&1&\\
t_3&t_6&-p&\\
\end{pmatrix}
=\begin{pmatrix}
&U\\
\tilde{U}&\\
 \end{pmatrix}.
\]
Tensoring with $A_n$ we obtain the windows corresponding to the base change from $R$ to $R_n$ and we denote them by $(M_{,1}(n), \alpha(n)).$
The window corresponding to $\overline{\mathbb{Y}}$ is described as above.
As before we write $ j(1):\overline{M}_{,1}^0 \rightarrow  M_1(1) $ and we write $j(1)=(X(1),Y(1)).$ As in \cite{KR1}, section 8 we may assume that $X(1)$ is of the form
$$ X(1)=
\begin{pmatrix}
p^{r/2}&0 \\
0&0\\
0&0\\
\end{pmatrix}.
$$
Using $Fj=jF$ and hence $j(\overline{1}_1)=F^{-1}jF(\overline{1}_1)$  one sees that 
$$ Y(1)=
\begin{pmatrix}
0&p^{r/2} \\
0&0\\
0&0\\
\end{pmatrix}.$$
Analogously as above we construct for any $l$ matrices $X_{\Q}(p^l)$ and $Y_{\Q}(p^l)$ over $A_{p^l}\otimes_{\Z}\Q$ as follows: $X_{\Q}(1)=X(1)$ and $Y_{\Q}(1)=Y(1)$ and 
\[ X_{\Q}(p^{l+1})=   \tilde{U}^{-1} \sigma (Y(p^{l}))S  \  \text{ and } \  Y_{\Q}(p^{l+1})=U^{-1}\sigma (X(p^{l}))S.\] (Again $\sigma (X_{\Q}(p^l))$ and $\sigma (Y_{\Q}(p^l))$ are well defined over $A_{p^{l+1}}\otimes_{\Z}\Q$.)
Then $j$ lifts over  $A_{p^l}$ if and only if  $X_{\Q}(p^l)$ and $Y_{\Q}(p^l)$ have entries in   $A_{p^l}$. One checks that 
$$ X(p)=
\begin{pmatrix}
p^{r/2}&0 \\
0&0\\
t_3p^{r/2-1}&0\\
\end{pmatrix}.
$$ and 
$$ Y(p)=
\begin{pmatrix}
0&p^{r/2} \\
0&0\\
0&0\\
\end{pmatrix}.$$
Hence in the case $r=0$ the equation of $\mathcal{Z}(j)$ is $t_3=0.$ This divisor (we call it $t$) contains $x$ as its only supersingular point. Inductively one checks now the following lemma.
\begin{Lem}
For any integer $k\geq 0$, the matrix 
$X(p^{2k+1})$ is of the form 
$$ X(p^{2k+1})= 
\begin{pmatrix}
0&0\\
0&0\\
t_3^{p^{2k}}p^{r/2-1-k}\cdot P&0\\
\end{pmatrix}
+A(p^{2k+1}),  
$$
where $P\in A_{p^{2k+1}}$  and $A(p^{2k+1})$ is a matrix with entries in $p^{r/2-k} A_{p^{2k+1}}.$
 The matrix 
$Y(p^{2k+1})$ has entries in $p^{r/2-k} A_{p^{2k+1}}.$
For any integer $k\geq 0$, the matrices
$X(p^{2k})$ and $Y(p^{2k})$ both have entries in $p^{r/2-k} A_{p^{2k}}.$
\end{Lem}

Thus  modulo $t_3^{p^r}$ the matrices  $X_{\Q}(p^{r+1})$ and $Y_{\Q}(p^{r+1})$ have entries in   $A_{p^r}$. Thus the divisor defined by $t_3^{p^r}=0$ belongs to $\mathcal{Z}(j)$, i.e. $t$ has (at least) multiplicity $p^r.$ There is a special homomorphism $y$ of valuation $0$ such that $y\perp j$ and $x\in \mathcal{Z}(y)(\F).$  It follows from \cite{KR1}, Proposition 8.2 that the length of the ring $\mathcal{O}_{\mathcal{Z}(j)_p \cap \mathcal{Z}(y)_p, x}$ is $1+p+...+p^{r}.$ We know that each of the $p+1$ supersingular curves passing through $x$ are contained in $\mathcal{Z}(j)_p$ with multiplicity $1+p^2+...+p^{r-2}.$ Further, as just shown, there is a subdivisor of $\mathcal{Z}(j)_p$ passing through $x$ defined by  $t_3^{p^r}=0.$ The intersection multiplicity of $\mathcal{Z}(y)_p$ with this divisor together with the supersingular part of $\mathcal{Z}(j)_p$ is thus at least $(1+p)(1+p^2+...+p^{r-2})+p^r=1+p+...+p^r$ which is already the length of $\mathcal{O}_{\mathcal{Z}(j)_p \cap \mathcal{Z}(y)_p, x}.$ Thus there are no other components of $\mathcal{Z}(j)_p$ passing through $x$ and the multiplicity of $t$ is precisely $p^r.$

It remains to check that for even $r$ the intersection multiplicity of $t$ with any of the curves $\mathcal{V}(\Lambda)$ passing through $x$ is $1.$ To this end we consider a special homomorphism $y$ of valuation $1$ such that $y\perp j$ and all $p+1$ curves $\mathcal{V}(\Lambda)$ passing through $x$ belong to $\mathcal{Z}(y).$ By  \cite{KR1}, Proposition 8.2, the length of the ring $\mathcal{O}_{\mathcal{Z}(j/p^{r/2})_p \cap \mathcal{Z}(y)_p, x}$ is $p+1$. From this the claim follows.  
\qed

\subsection{Further results}
\begin{Def}
Let $j$ be a special homomorphism. The difference divisor $\mathcal{D}(j)=\mathcal{Z}(j)-\mathcal{Z}(j/p)$ is defined as follows. Suppose $\mathcal{Z}(j)$ is locally given by the equation $f=0$ and $\mathcal{Z}(j/p)$ is locally given by $g=0.$  Then $\mathcal{D}(j)$ is the divisor locally given by $f\cdot g^{-1}=0$ (note that $g$ divides $f$). 
\end{Def}
\begin{Lem} \label{pmultlem}
Let  $j $ be a special homomorphism of nonnegative valuation,  
let
  $x \in \mathcal{Z}(j)(\F)$  and let $D\subset \mathcal{N}$ be  a divisor which is regular in $x$. Suppose that locally around $x$ we have  $D_p \subset \mathcal{Z}(j)_p$. 
  Let  $R= \mathcal{O}_{D,x}$ and let $(f)$ be the ideal of $\mathcal{Z}(j)\cap D$ in $R$ (i.e., $\mathcal{O}_{\mathcal{Z}(j)\cap D,x}=R/(f)$).
   Then  the ideal of $\mathcal{Z}(pj)\cap D$ in $R$ is  $( p \cdot f).$ In other words, (if $f \neq 0$) the equation of  $\mathcal{D}(pj) \cap D$ in $R$ is given by $p=0.$
\end{Lem}
 This is a stronger version of Lemma 3.12 in \cite{T}. The first part of its proof is analogous, but we prove it here completely.

\emph{Proof.} 
We may assume that $f \neq 0.$ 

The canonical map
\[
\varrho_0: \Spf R/(f) \rightarrow D
\]
factors  through $\mathcal{Z}(j)\cap D$. 
\smallskip

\begin{bf} Claim: \end{bf} 
\emph{The canonical map}
\[
\varrho_{1}:\Spf R/(pf) \rightarrow D
\]
\emph{factors through  $\mathcal{Z}(pj)\cap D$.} 
\smallskip

By our assumption  $D_p \subset \mathcal{Z}(j)_p$ (locally around $x$), we know that $f$ is divisible by $p$.
We have $R/(f)=(R/(pf))/I$, where $I=(f)/(pf)$. Since $f$ is divisible by $p$, the ideal  $I$ carries a nilpotent  $pd$-structure.  Hence we may apply   Grothendieck-Messing theory for the pair $R/(pf)$, $R/(f)$.
We denote by  $M$  the value of the crystal of the  $R/(pf)$-valued point $\varrho_1$ in $R/(pf)$, and by $\overline{M}$ the value of the crystal of the $R/(f)$-valued point $\varrho_0$ in $R/(f)$. Then $\overline{M}=M\otimes R/(f)$. 
The value of the crystal of $\overline{Y}\times_W R/(pf)$ resp.  $\overline{Y}\times_W R/(f)$ in $R/(pf)$ resp. $R/(f)$ is  $\overline{\mathbb{M}}^0_{R/(pf)} = \overline{\mathbb{M}}^0\otimes_W R/(pf)$ resp. $\overline{\mathbb{M}}^0_{R/(f)} = \overline{\mathbb{M}}^0\otimes_W R/(f).$ 
(Recall that $\overline{\mathbb{M}}^0$ denotes the Dieudonn\'e module of $\overline{\mathbb{Y}}$.)

Denote by  ${\mathcal{F}}\hookrightarrow M$ the Hodge filtration of  $\varrho_{1}$ and by $\overline{\mathcal{F}}\hookrightarrow \overline{M}$ the Hodge filtration of $\varrho_{0}$. Then the Hodge filtration of $\varrho_{1}$ lifts the  Hodge filtration of $\varrho_{0}$.
We write $F_{\overline{Y}_{R/(pf)}}\rightarrow \overline{\mathbb{M}}^0_{R/(pf)}$ resp. $F_{\overline{Y}_{R/(f)}}\rightarrow \overline{\mathbb{M}}^0_{R/(f)}$ for  the Hodge filtrations of $\overline{Y}_{R/(pf)}$ resp. $\overline{Y}_{R/(f)}$. Here  $F_{\overline{Y}_{R/(pf)}}$ resp. $F_{\overline{Y}_{R/(f)}}$ are generated by $\overline{1}_0.$
 We have a map $j:\overline{\mathbb{M}}^0_{R/(f)} \rightarrow \overline{M}$ with $jF_{\overline{Y}_{R/(f)}}\subset \overline{\mathcal{F}}$. It lifts to a map $j:\overline{\mathbb{M}}^0_{R/(pf)} \rightarrow {M}$, and we have to show that $pj F_{\overline{Y}_{R/(pf)}} \subset {\mathcal{F}}$.  
Since  $jF_{\overline{Y}_{R/(f)}}\subset \overline{\mathcal{F}}$, we have  $jF_{\overline{Y}_{R/(pf)}}\subset{\mathcal{F}}+IM$. Since $pI=0$ it follows that  $pjF_{\overline{Y}_{R/(pf)}}\subset{\mathcal{F}}.$ This confirms the claim.

Let $\mathcal{O}_{\mathcal{Z}(pj)\cap D,x}=R/(g)$. Since $D$ is regular in $x$ we know that $R$ is a UFD. We write $g=h\cdot s \cdot f$, where $h$ is prime to $p$ and $V(s)$ has support in the special fibre of $D$. 
We just saw that $s$ is divisible by $p$ and want to show next that (up to a unit)  $s=p$. 
Suppose that $q$ is prime divisor of $p$ in $R$ such that $g$ divisible by $qpf$. 
Let $C=R/(qpf)$, let $\overline{C}=R/(qf)$, and let $\overline{\overline{C}}=R/(f)$. Since $p>2$, the ideals $(f)/(pf)$ and ($f)/(pqf)$ are equipped  nilpotent $pd$-structures.
We consider the canonical $C$-, resp. $\overline{C}$ -, resp. $\overline{\overline{C}}$- valued points 
\[
\varrho_C: \Spf C \rightarrow D, \ \ \ 
\varrho_{\overline{C}}: \Spf \overline{C} \rightarrow D, \ \ \ 
\varrho_{\overline{\overline{C}}}: \Spf \overline{\overline{C}} \rightarrow D.
\]

\begin{bf} Claim: \end{bf} {\em$\varrho_{\overline{C}}$ factors through $\mathcal{Z}(j)\cap D$.}
\smallskip

We denote by    $M$ the value of the crystal of the  $\varrho_C$-valued point  in $C$,  analogously we define   $\overline{M}$ and  $\overline{\overline{M}}$. 
Further, let ${\mathcal{F}}\hookrightarrow M$ by the Hodge filtration of $\varrho_{C}$ and analogously we define $\overline{\mathcal{F}}\hookrightarrow \overline{M}$ and $\overline{\overline{\mathcal{F}}}\hookrightarrow \overline{\overline{M}}$. 
The filtrations  $F_{\overline{Y}_{C}}\rightarrow \overline{\mathbb{M}}^0_{C}$ resp. $F_{\overline{Y}_{\overline{C}}}\rightarrow \overline{\mathbb{M}}^0_{\overline{C}}$ resp. $F_{\overline{Y}_{\overline{\overline{C}}}}\rightarrow \overline{\mathbb{M}}^0_{\overline{\overline{C}}}$ are defined as above.
Let $f_{C}$ be a generator of  $F_{\overline{Y}_{C}}$. Its images $f_{\overline{C}}$  in  $F_{\overline{Y}_{\overline{C}}}$ resp. $f_{\overline{\overline{C}}}$ in $F_{\overline{Y}_{\overline{\overline{C}}}}$ are  also generators. Again, we have a map $j:\overline{\mathbb{M}}^0_{\overline{\overline{C}}}\rightarrow  \overline{\overline{M}}$ respecting the filtrations. It lifts to maps (which we denote by $j$, too)  $\overline{\mathbb{M}}^0_{{\overline{C}}}\rightarrow {\overline{M}}$ and  $\overline{\mathbb{M}}^0_{{{C}}}\rightarrow {{M}}$.   Suppose that $M$ is the direct sum of $F$ and the span of $b_1,b_2,b_3.$ Then we may write $j(f_C)=a+\overline{x_1}b_1+\overline{x_2}b_2+\overline{x_3}b_3,$ where $a \in F$  and $\overline{x_i} \in R/(pqf).$ Let $x_i \in R$ for $i=1,2,3$ be lifts of $\overline{x_i}.$ Then it follows that $px_i \in (pqf),$  hence $x_i \in (qf).$ Therefore the image of  $j(f_C)$ in $\overline{M}$ lies in $\overline{F}.$ This confirms the claim.

Finally, we want to show that $h$ is a unit. This follows from the following
\smallskip

\begin{bf} Claim: \end{bf} {\em If $ H \in R$ is coprime to $p$ and $ \Spf R/(H) \rightarrow D$ factors through $\mathcal{Z}(pj) \cap D$, then it also factors through $\mathcal{Z}(j) \cap D$.} 
\smallskip

For any positive integer $n$, let $R_n=R/(H,p^n).$ Since $p>2$, the ring $R_{n+2}$ is always a  pd-extension of $R_n$ (the ideal $(p^n)/(p^{n+2})$ in $R_{n+2}$ carries a nilpotent pd-structure). We show by induction that  $ \Spf R_n \rightarrow D$ factors through $\mathcal{Z}(j) \cap D$. From this the claim follows. 
For $n=1$, this is true since by assumption   $ \Spf R/(p) \rightarrow D$ factors through $\mathcal{Z}(j)\cap D$. We denote by $F_n \rightarrow M_n$ the Hodge filtration of   $ \Spf R/(p^n,H) \rightarrow D$ und by  $F_{\overline{Y}_n} \rightarrow  \overline{\mathbb{M}}^0_n$ the Hodge filtration of $\overline{Y}_{R_n}.$ 
Suppose   $ \Spf R/(p^n,H) \rightarrow D$ factors through $\mathcal{Z}(j) \cap D$. We want to show that   $ \Spf R/(p^{n+1},H) \rightarrow D$ factors through $\mathcal{Z}(j)\cap D.$
By the induction hypothesis, we have a morphism $j:\overline{\mathbb{M}}^0_n \rightarrow M_n$ which respects the filtrations. It induces morphisms  $j:\overline{\mathbb{M}}^0_{n+1} \rightarrow M_{n+1}$ and $j:\overline{\mathbb{M}}^0_{n+2} \rightarrow M_{n+2}.$ Let $f_{ n+2}$ be a generator of  $F_{\overline{Y}_{n+2}}$. Its images $f_{n+1}$  in  $F_{\overline{Y}_{n+1}}$ resp $f_n$ in $F_{\overline{Y}_n}$ are  also generators. Suppose that $M_{n+2}$ is the direct sum of $F_{n+2}$ and the span of $b_1,b_2,b_3.$ Then we may write $j(f_{n+2})=a+\overline{x_1}b_1+\overline{x_2}b_2+\overline{x_3}b_3,$ where $a \in F_{n+2}$  and $\overline{x_i} \in R_{n+2}.$ Let $x_i\in R$ for $i=1,2,3$ be lifts of the $\overline{x_i}.$ Then it follows that $px_i \in (H,p^{n+2}),$  hence $px_i$ is of the form $px_i= y_iH +z_ip^{n+2}$  for some $y_i,z_i \in R.$  Since $H$ is coprime to $p$ it follows that $y_i$ is divisible by $p.$ We conclude that $x_i \in (H, p^{n+1})$ and  hence $j(f_{n+1})\in F_{n+1}.$ This confirms the claim and ends the proof.
\qed
\begin{Lem} \label{keigemcomp}
Let $j_1,j_2$ be special homomorphisms of nonnegative valuation which are linearly independent. Then the divisors $\mathcal{Z}(j_1)$ and $\mathcal{Z}(j_2)$ do not have any common components. 
\end{Lem}
{\em Proof.} Let $x \in \mathcal{Z}(j_1)(\F)\cap \mathcal{Z}(j_2)(\F).$ If $x$ is superspecial, then we can choose a special homomorphism $j$ of valuation $0$ which is linearly independent of $j_1,j_2$ and for which $x \in \mathcal{Z}(j).$ Then by \cite{KR1}, Theorem 5.2, the intersection $\mathcal{Z}(j)\cap \mathcal{Z}(j_1)\cap \mathcal{Z}(j_2)$ has dimension $0$. Hence the dimension of $\mathcal{Z}(j_1)\cap \mathcal{Z}(j_2)$ in $x$ is $1$, hence there is no common component in $x$. 
If $x$ is not superspecial and there is a  common component $C$ of $\mathcal{Z}(j_1)$ and $\mathcal{Z}(j_2)$ passing through $x$, then it follows that there is a common component $S$  of $\mathcal{Z}(j_1)_p$ and  $\mathcal{Z}(j_2)_p$ which also belongs to $C$. But  then it follows that $S$ is of the form $\mathcal{V}(\Lambda)$  and that, for any superspecial point $z$ in $S$, we have $z\in C$ which is not possible as we just saw. \qed

\begin{Pro}\label{reg}
Let $j$ be a special homomorphism with $h(j,j)\neq 0$. Then the divisor $\mathcal{D}(j)$ is regular.
\end{Pro}
{\em Proof.} We check for any $x\in \mathcal{Z}(j)(\F)$ that $\mathcal{D}(j)$ is regular in $x$. 
First suppose that $x \in \mathcal{V}(\Lambda)(\F)\subset \mathcal{Z}(j)(\F)$ and that in the case that $x$ is superspecial there is no other curve $ \mathcal{V}(\Lambda^{'})$ which contains $x$ and belongs to $\mathcal{Z}(j)$.
 Futher suppose that the multiplicity of $ \mathcal{V}(\Lambda)$ in $\mathcal{Z}(j)_p$ is $1$. Thus locally around $x$ we have $\mathcal{D}(j)=\mathcal{Z}(j).$ Then by Theorem \ref{spezf}, the special  fiber of $\mathcal{D}(j)$ is regular in $x$. Since the dimension of   $\mathcal{D}(j)_p$ is smaller than the dimension of  $\mathcal{D}(j)$, it follows that $\mathcal{D}(j)$ is regular in $x$.

Next we suppose again that $x \in \mathcal{V}(\Lambda)(\F)\subset \mathcal{Z}(j)(\F)$,  but that the multiplicity of $ \mathcal{V}(\Lambda)$ in $\mathcal{Z}(j)_p$ is greater than $1.$ It is easy to see that there is a special homomorphism $y$ which is linearly independent of $j$ such that 
$ \mathcal{V}(\Lambda)\subset \mathcal{Z}(y)$ and such that the multiplicity of  $ \mathcal{V}(\Lambda)$ in $\mathcal{Z}(y)_p$ is $1$ and such that  (in the case that $x$ is superspecial)  there is no other curve $ \mathcal{V}(\Lambda^{'})$ which contains $x$ and belongs to $\mathcal{Z}(y).$ We already know that $\mathcal{Z}(y)$ is regular in $x$. By Lemma \ref{keigemcomp}, there is no common component of $\mathcal{Z}(j)$ and $\mathcal{Z}(y)$, in particular the equation of $\mathcal{Z}(j) \cap \mathcal{Z}(y)$ in $\mathcal{O}_{\mathcal{Z}(y),x}$ is nontrivial. Now Lemma \ref{pmultlem} shows that locally around  $x$ the intersection $\mathcal{D}(j) \cap \mathcal{Z}(y)$ equals $\mathcal{Z}(y)_p$ which is $\mathcal{V}(\Lambda)$. Hence $\mathcal{D}(j)$ is regular in $x$.

Next we suppose that $x$ does not satisfy the conditions of one of the preceding cases. Then either $j$ is of  valuation $r=0$ or $r=2$, and $x$ is the superspecial point of $\mathcal{Z}(j/p^r)$, or  the valuation of $j$ is $1$  and $x$ is the intersection of two (and hence $p+1$) curves (of the form  $\mathcal{V}(\Lambda)$) which belong to $\mathcal{Z}(j)$. First suppose the latter case, i.e. $j$ has valuation $1$. We find special homomorphisms $y_1,y_2$ which are both of valuation $0$ and such that $y_1,y_2$ and $j$ are pairwise perpendicular to each other. Then the intersection $\mathcal{Z}(y_1)\cap \mathcal{Z}(y_2)\cap \mathcal{Z}(j)$ has by \cite{KR1}, Theorem 5.1 length $1$. From this it follows that $\mathcal{Z}(j)$ is regular in $x.$ 

Finally, we suppose that  $j$ is of  valuation $r=0$ or $r=2$ and $x$ is the superspecial point of $\mathcal{Z}(j/p^r)$. If $r=0$, we choose a special homomorphism $y$ which has valuation $0$ and such that  $y \perp j$. Then the length of $\mathcal{O}_{\mathcal{Z}(j)_p\cap \mathcal{Z}(y)_p,x}$ is $1$, hence   $\mathcal{Z}(j)_p$ is regular in $x$, hence $\mathcal{Z}(j)$ is also regular in $x.$
In the case $r=2$ we consider the divisor $\mathcal{Z}(j/p)$ which passes through $x$. Suppose it is given in $\mathcal{O}_{\mathcal{N},x}$ by the equation $f=0.$ We consider the divisor $D$ in $\Spf \mathcal{O}_{\mathcal{N},x}$ given by  $f+p=0.$ Its special fiber is regular since it is the same as the special fiber of $\mathcal{Z}(j/p),$ hence $D$ is regular in $x$. Now Lemma \ref{pmultlem} shows that $D \cap \mathcal{D}(j)$ equals $D_p,$ hence it is regular. Thus $\mathcal{D}(j)$ is also regular in $x.$ \qed

\section{Intersections of special cycles}
\subsection{$GL_3$-invariance and Multilinearity}
\begin{Lem} \label{Lweg}
Let $j_1,j_2$ be linearly independent special homomorphisms.  Then $$\mathcal{O}_{\mathcal{Z}(j_1)}\otimes^{\mathbb{L}}\mathcal{O}_{\mathcal{Z}(j_2)}= \mathcal{O}_{\mathcal{Z}(j_1)}\otimes \mathcal{O}_{\mathcal{Z}(j_2).}$$ More precisely, the object on the left hand side is represented in the derived category by the object on the right hand side.
The same formula holds if $\mathcal{Z}(j_1)$ or $\mathcal{Z}(j_2)$ or both are replaced by $\mathcal{D}(j_1)$ resp. $\mathcal{D}(j_2)$.
\end{Lem}
\emph{Proof.}
This is a simplified version of the proof of Lemma 4.1 in \cite{T}.
 Let $x\in \mathcal{N}(\F)$ and let $R=\mathcal{O}_{ \mathcal{N}, x}$. Let the ideals of $\mathcal{Z}(j_1)$ and $\mathcal{Z}(j_2)$ in $R$ be generated by  $f_1$ and $f_2$. We consider the  exact sequence
\begin{equation*}
\begin{CD}
0 @ >>> R @>f_1\cdot>>R @>>>
R/(f_1)@>>> 0.
\end{CD}
\end{equation*}
Tensoring this with $R/(f_2)$, we see that $$
\mathcal{TOR}_1(\mathcal{O}_{\mathcal{Z}(j_1)}, \mathcal{O}_{\mathcal{Z}(j_2)})_x=\ker(R/(f_2)\stackrel{f_1 \cdot}{\longrightarrow}R/(f_2)).
$$ To show that this vanishes
 we have to show that $f_1$ and $f_2$ have no common divisor in the regular ring $R$.  This follows from the fact that $\mathcal{Z}(j_1)$ and $\mathcal{Z}(j_2)$ have no common component, see Lemma \ref{keigemcomp}. 
The proof in the case that $\mathcal{Z}(j_1)$ or $\mathcal{Z}(j_2)$ or both are replaced by $\mathcal{D}(j_1)$ resp. $\mathcal{D}(j_2)$ is the same.
 \qed

\begin{Pro} \label{diag}
Let $j_1,j_2,j_3$ be linearly independent special homomorphisms. Then the derived tensor product $ \mathcal{O}_{\mathcal{Z}(j_1)}\otimes^{\mathbb{L}}\mathcal{O}_{\mathcal{Z}(j_2)}\otimes^{\mathbb{L}}\mathcal{O}_{\mathcal{Z}(j_3)}$ depends only on the $\Zpq$-span $\bf{j}$ of $j_1,j_2,j_3$ in $\mathbb{V}$.
\end{Pro}

\emph{Proof.}
Again this is a simplified version of the proof of Proposition 4.2 in \cite{T}.
 First we observe that for any basis $y_1, y_2, y_3$ of $\bf{j}$ the derived tensor product $ \mathcal{O}_{\mathcal{Z}(y_1)}\otimes^{\mathbb{L}}\mathcal{O}_{\mathcal{Z}(y_2)}\otimes^{\mathbb{L}}\mathcal{O}_{\mathcal{Z}(y_3)}$ is invariant under any permutation of the $y_i$.
\smallskip

{\bf Claim} $ \mathcal{O}_{\mathcal{Z}(y_1)}\otimes^{\mathbb{L}}\mathcal{O}_{\mathcal{Z}(y_2)}\otimes^{\mathbb{L}}\mathcal{O}_{\mathcal{Z}(y_3)}$  \emph{ is invariant if one $y_i$ is replaced by $\varepsilon y_i + zy_l$, where $\varepsilon \in \Zpq^{\times}, \ z \in \Zpq$ and $l \neq i$. }
\smallskip

By the preceding lemma we have $\mathcal{O}_{\mathcal{Z}(y_i)}\otimes^{\mathbb{L}}\mathcal{O}_{\mathcal{Z}(y_l)}=\mathcal{O}_{\mathcal{Z}(y_i)}\otimes \mathcal{O}_{\mathcal{Z}(y_l)}$ $=\mathcal{O}_{\mathcal{Z}(y_i) \cap \mathcal{Z}(y_l)}$ which only depends on the $\Zpq$-span of $y_i, y_l$ in $\mathbb{V}$.
From this the claim follows.

Since we can transform the basis $y_1,y_2,y_3$ by a suitable sequence of  permutations and operations as in the claim into any other basis of $\bf{j}$, the claim of the proposition follows.
\qed

\begin{Pro}\label{multlin}
Let $j_1,j_2,j_3$ be  special homomorphisms  such that the intersection $\mathcal{Z}(j_1) \cap \mathcal{Z}(j_2) \cap \mathcal{Z}(j_3) $ is non-empty and the corresponding fundamental matrix is non-singular. 
Then 
$$
(\mathcal{Z}(j_1),\mathcal{Z}(j_2),\mathcal{Z}(j_3))=\sum_{l_1,l_2,l_3}(\mathcal{D}(j_1/p^{l_1}),\mathcal{D}(j_2/p^{l_2}),\mathcal{D}(j_3/p^{l_3})),
$$
where the sum is taken over all possible triples $(l_1,l_2,l_3)$ (i.e. setting $a_i=\nu_p(j_i)$, we have $l_1 \leq [a_1/2],\  l_2 \leq [a_2/2], \ l_3 \leq [a_3/2]$, where $[ \ ]$ denotes Gauss brackets).
\end{Pro}

We use here that $(\mathcal{Z}(j_1),\mathcal{Z}(j_2),\mathcal{Z}(j_3))$ is finite. This follows from the fact that  $\mathcal{Z}(j_1) \cap \mathcal{Z}(j_2) \cap \mathcal{Z}(j_3) $ has support in the supersingular locus (cf. \cite{KR1}, section 4) and that this support is proper over $\F$ (follows from the results of section 2).

\emph{Proof.} 
Again, this is a simplified version of the proof of Proposition 4.3 in \cite{T}.
Let $r=[\nu_p(j_1)/2].$ Denote by $\mathcal{I}$ the ideal sheaf of $D:=\mathcal{D}(j_1/p^r)(=\mathcal{Z}(j_1/p^r))$ in $\mathcal{O}_{\mathcal{N}}$. Denote by $\mathcal{J}$ the ideal sheaf of $\Delta:=\mathcal{Z}(j_1)-\mathcal{D}(j_1/p^r)$ (the notation in the latter expression is meant in the same sense as the notation $\mathcal{D}(j)=\mathcal{Z}(j)-\mathcal{Z}(j/p)$). Our first aim is to show that  $$(\mathcal{Z}(j_1),\mathcal{Z}(j_2),\mathcal{Z}(j_3))=(\Delta,\mathcal{Z}(j_2), \mathcal{Z}(j_3))+(D,\mathcal{Z}(j_2),\mathcal{Z}(j_3)).$$ 

We consider the canonical short exact sequences 
$$
\begin{CD}
0 @ >>> \mathcal{J}/(\mathcal{J}\cdot \mathcal{I}) @>>> \mathcal{O}_{\mathcal{N}}/(\mathcal{J}\cdot \mathcal{I}) @>>>\mathcal{O}_{\mathcal{N}}/\mathcal{J}@>>> 0
\end{CD}
$$
and 
$$
\begin{CD}
0 @ >>> (\mathcal{J}+\mathcal{I})/ \mathcal{I}  @>>> \mathcal{O}_{\mathcal{N}}/ \mathcal{I} @>>>\mathcal{O}_{\mathcal{N}}/(\mathcal{J}+\mathcal{I})@>>> 0.
\end{CD}
$$
\begin{Lem}
The inclusion $\mathcal{J}\cdot \mathcal{I} \hookrightarrow \mathcal{J}\cap \mathcal{I}$ is an equality.
\end{Lem}
\emph{Proof.} We can check this locally. Let $x \in \mathcal{N} (\F)$ and let $R= \mathcal{O}_{ \mathcal{N}, x}$. For any $l\leq r$, let $(f_l)$ be the ideal of $\mathcal{D}(j_1/p^l)$ in $R$. By Proposition \ref{reg}, each $f_l$ is a prime element in $R$. It follows from Theorem \ref{spezf}  that the $f_l$ are pairwise coprime (more precisely, Theorem \ref{spezf} shows   that even the special fibres of the several $\Spf R/(f_l)$ are pairwise distinct). Since $\mathcal{I}_x$ equals $(f_r)$ and $\mathcal{J}_x$ equals $(\prod_{l\neq r} f_l)$,  the claim of the lemma follows.
\qed
\newline

By the first exact sequence have $$(\mathcal{Z}(j_1),\mathcal{Z}(j_2),\mathcal{Z}(j_3))=(\Delta,\mathcal{Z}(j_2), \mathcal{Z}(j_3))+ \chi(\mathcal{J}/(\mathcal{J}\cdot \mathcal{I})\otimes^{\mathbb{L}}\mathcal{O}_{\mathcal{Z}(j_2)}\otimes^{\mathbb{L}}\mathcal{O}_{\mathcal{Z}(j_3)}).$$
 
Now using the lemma we see that $ (\mathcal{J}+\mathcal{I})/ \mathcal{I}=\mathcal{J}/(\mathcal{J}\cap \mathcal{I})=\mathcal{J}/(\mathcal{J}\cdot \mathcal{I})$. This shows together with the second exact sequence that
\begin{equation*}
 \begin{split} \chi(\mathcal{J}/(\mathcal{J}\cdot \mathcal{I})\otimes^{\mathbb{L}}\mathcal{O}_{\mathcal{Z}(j_2)}\otimes^{\mathbb{L}}\mathcal{O}_{\mathcal{Z}(j_3)})=  (D,\mathcal{Z}(j_2),\mathcal{Z}(j_3)) 
-(D\cap \Delta,\mathcal{Z}(j_2),\mathcal{Z}(j_3)). 
\end{split} 
\end{equation*}

Thus in order to show that $$(\mathcal{Z}(j_1),\mathcal{Z}(j_2),\mathcal{Z}(j_3))=(\Delta,\mathcal{Z}(j_2), \mathcal{Z}(j_3))+(D,\mathcal{Z}(j_2),\mathcal{Z}(j_3))$$ it remains to show the 
\smallskip  

{\bf Claim} $(D\cap \Delta,\mathcal{Z}(j_2),\mathcal{Z}(j_3))=0$. 
\smallskip

By the proof of the lemma we know  $D$ and $\Delta$ have no common component, hence  $\codim D\cap \Delta =2.$

By Corollary  \ref{cor1} below,  $\Delta \cap D$ is as a divisor in $D$ of the form $rD_p+h$, where 
$h$ is a divisor of the form $\sum_x h_x$, and where the sum runs over a discrete set of $\F$-valued points of $D$,  and $h_x$ is a horizontal divisor (meaning that its equation is coprime to $p$) meeting the underlying reduced subscheme of $D$ only in $x$. Then $(D\cap \Delta,\mathcal{Z}(j_2),\mathcal{Z}(j_3))= (rD_p, \mathcal{Z}(j_2), \mathcal{Z}(j_3))+\sum (h_x,\mathcal{Z}(j_2),\mathcal{Z}(j_3))- (h \cap rD_p,\mathcal{Z}(j_2), \mathcal{Z}(j_3))$. (This is shown by the same reasoning as above.) 
We observe that $(rD_p, \mathcal{Z}(j_2)\cap \mathcal{Z}(j_3))=0$ as  follows from the exact sequence
$$
\begin{CD}
0 @ >>> \mathcal{O}_D @> p^{r}\cdot>> \mathcal{O}_D @>>>\mathcal{O}_D/(p^{r})@>>> 0.
\end{CD}
$$
Thus  $(D\cap \Delta,\mathcal{Z}(j_2),\mathcal{Z}(j_3))=\sum (h_x,\mathcal{Z}(j_2),\mathcal{Z}(j_2)) - (h \cap rD_p,\mathcal{Z}(j_2), \mathcal{Z}(j_3))$. But since   $\codim h_x +  \codim \mathcal{Z}(j_2) + \codim \mathcal{Z}(j_3)=  \codim D\cap \Delta + \codim \mathcal{Z}(j_2) + \codim \mathcal{Z}(j_3) = \dim \mathcal{N}+1$, it follows that each  $(h_x,\mathcal{Z}(j_2),\mathcal{Z}(j_2))$ vanishes, see \cite{SABK}, Corollary 1 on p.15.
In the same way we see that $(h \cap rD_p,\mathcal{Z}(j_2), \mathcal{Z}(j_3))=0$. This confirms the claim.

Now repeating this reasoning, we see that $$(\mathcal{Z}(j_1),\mathcal{Z}(j_2),\mathcal{Z}(j_3))=\sum_{l_1}(\mathcal{D}(j_1/p^{l_1}),\mathcal{Z}(j_2),\mathcal{Z}(j_3)).$$ The remaining multilinearity in the other variables follows in the same way.
\qed 
\begin{Lem}
Let $j_1,j_2,j_3$ be special homomorphisms such that $\mathcal{Z}(j_1)\cap \mathcal{Z}(j_2) \cap \mathcal{Z}(j_3)$ has dimension $0.$ Then $$\mathcal{O}_{\mathcal{Z}(j_1)}\otimes^{\mathbb{L}}\mathcal{O}_{\mathcal{Z}(j_2)}\otimes^{\mathbb{L}} \mathcal{O}_{\mathcal{Z}(j_3)}=\mathcal{O}_{\mathcal{Z}(j_1)}\otimes \mathcal{O}_{\mathcal{Z}(j_2)}\otimes \mathcal{O}_{\mathcal{Z}(j_3)}$$
(meaning that the right hand side represents the left hand side in the derived category). In particular,
$$(\mathcal{Z}(j_1),\mathcal{Z}(j_2),\mathcal{Z}(j_3))=\lg_W(\mathcal{O}_{\mathcal{Z}(j_1)\cap \mathcal{Z}(j_2) \cap \mathcal{Z}(j_3),x}),$$ where $x$ is the unique $\F$-valued point of $\mathcal{Z}(j_1)\cap \mathcal{Z}(j_2) \cap \mathcal{Z}(j_3)$.
\end{Lem}
This can be proved as \cite{T} Proposition  4.2 or \cite{KR2}, Proposition 11.6. \qed

Thus we may apply in the situation of the lemma the results of \cite{KR1}, Theorem 5.2 to determine $(\mathcal{Z}(j_1),\mathcal{Z}(j_2),\mathcal{Z}(j_3)).$

\subsection{Vertical components of $\mathcal{D}(j_1)\cap \mathcal{D}(j_2)$}

Let $\Lambda \in \mathcal{L}^{(3)}$. Let $j_1, j_2 $ be linearly independent special homomorphisms of nonnegative valuations such that  $\mathcal{V}(\Lambda)\subset \mathcal{D}(j_1)\cap \mathcal{D}(j_2)$. The divisor $\mathcal{D}(j_i)_p$ in $\mathcal{N}_p$ contains $ \mathcal{V}(\Lambda)$ with multiplicity $p^{r_i}$ for some even integer $r_i$, see Theorem \ref{spezf}. Now $ \mathcal{V}(\Lambda)$ is a closed irreducible reduced subscheme of codimension $1$ in (the regular) $\mathcal{D}(j_i)$, hence a prime divisor. The subsequent two propositions will tell us the multiplicity of   $\mathcal{V}(\Lambda)$ in the divisor $\mathcal{D}(j_1) \cap \mathcal{D}(j_2)$ in $\mathcal{D}(j_1)$ (or $\mathcal{D}(j_2)$).
\begin{Pro} \label{unglmult}
In the situation just described, suppose further that $r_1 \neq r_2 $. Then the divisor $\mathcal{D}(j_1)\cap \mathcal{D}(j_2)$ in  $\mathcal{D}(j_1)$ contains  $\mathcal{V}(\Lambda)$ with multiplicity $p^{\min\{r_1,r_2\}}$.
\end{Pro}
\emph{Proof.} 
We choose a special homomorphism $j$ of valuation $1$  such that  $\mathcal{V}(\Lambda)\subset \mathcal{Z}(j)$ and $j\perp j_1.$ (It is easy to see such  a special homomorphism $j$ exists.)
We choose $x \in \mathcal{V}(\Lambda)(\F)$ not superspecial such that there is no horizontal component of $\mathcal{Z}(j)\cap \mathcal{Z}(j_1)$ passing through $x$ (comp. Lemma \ref{keigemcomp}).
Let $R= \mathcal{O}_{\mathcal{N}, x}$ and let $(f_1), \ (f_2), \ (f)$ be the ideals of $\mathcal{D}(j_1), \ \mathcal{D}(j_2), \  \mathcal{Z}(j)$ in $R$. 
\smallskip

{\bf Claim 1} \emph{ The ideal of $\mathcal{V}(\Lambda)$ in $R$ is given by $(p,f)=(f,f_1)$.}
\smallskip

The fact that the  ideal of $\mathcal{V}(\Lambda)$ in $R$ is given by $(p,f)$ follows from Theorem \ref{spezf}.
If $r_1>0$, then Lemma \ref{pmultlem} shows that then  ideal of $\mathcal{V}(\Lambda)$ in $R$ is given by $(f,f_1).$ 
Suppose now that $r_1=0.$ We choose a superspecial point $z \in \mathcal{V}(\Lambda)(\F)$  such that all $p+1$ curves   $\mathcal{V}(\Lambda^{'})$ passing through $z$ belong to $\mathcal{Z}(j)$ but no curve $\mathcal{V}(\Lambda^{'})\neq\mathcal{V}(\Lambda)$ passing through $z$ belongs to $\mathcal{Z}(j_1).$ We find special homomorphisms $j_0, \tilde{j_0}$ of valuation $0$ such that  $j_0, \tilde{j_0}, j$ are all perpendicular to each other and $z\in \mathcal{Z}(j_0)(\F)$ and $z\in \mathcal{Z}(\tilde{j_0})(\F).$
Then $j_1=\alpha j_0 + \beta \tilde{j_0}$ for some $\alpha, \beta \in \Zpq^{\times}.$ Hence the intersection multiplicity $(\mathcal{Z}(j_0),\mathcal{Z}(j_1),\mathcal{Z}(j))$ equals  $(\mathcal{Z}(j_0),\mathcal{Z}(\tilde{j_0}),\mathcal{Z}(j))$ which is $1$ by \cite{KR1}, Theorem 5.1. From this it follows that the multiplicity of $ \mathcal{V}(\Lambda)$ in $\mathcal{Z}(j)\cap \mathcal{Z}(j_1)$ as divisor in $\mathcal{Z}(j_1)$ (or $\mathcal{Z}(j)$) is $1$. Since  there is no horizontal component of $\mathcal{Z}(j)\cap \mathcal{Z}(j_1)$ passing through $x$, this shows the claim.

The rest of the proof is more or less the same as  the proof of Proposition 4.6 in \cite{T} but we prove it completely. 
\smallskip

{\bf Claim 2} \emph{We have the following identity of ideals in $R$:} $(f^{p^{r_1}},f_1)=(f^{p^{r_1}},p)=(f_1,p)$. \smallskip

For any $z \in R$, we denote by $\overline{z}$ the image of $z$ in $\overline{R}:=R/(p)$ which is a UFD.
First we observe that $(f^{p^{r_1}},p)=(f_1,p)$, since the equation of $\mathcal{V}(\Lambda)$ in $\overline{R}$ is given by $\overline {f}=0$.
After perhaps multiplying $f_1$ by a unit we can therefore write 
$$f_1=f^{p^{r_1}}+p \varrho, $$ for some $\varrho \in R$. By the descripion of the ideal of $\mathcal{V}(\Lambda)$ in $R$ (claim 1) we can also write
$$f_1=\varepsilon p + \sigma f$$ for some $\varepsilon \in R^{\times}$ and some $\sigma \in R$.
From these equations we get $\overline{f_1}=\overline{f}^{p^{r_1}}= \overline{\sigma}\overline{f}$, hence $\overline{\sigma}=\overline{f}^{p^{r_1}-1}$, hence $\sigma=f^{p^{r_1}-1}+p\sigma^{'}$ for some $\sigma^{'}\in R$.
Hence we have $f_1=\varepsilon p+ f^{p^{r_1}}+pf\sigma^{'}= p(\varepsilon+f\sigma^{'})+f^{p^{r_1}}$. Since $\varepsilon+f\sigma^{'}$ is a unit it follows that $p\in (f_1, f^{p^{r_1}})$.
Hence $(f_1,f^{p^{r_1}})=(f_1,f^{p^{r_1}},p)=(f_1,p)$. This confirms the claim.

Now we distinguish the cases $r_2 < r_1$ and $r_1 < r_2$.

{\bf First case} $r_2 < r_1$.
Since the ideal of $\mathcal{V}(\Lambda)$ in the local ring of $\mathcal{D}(j_1)$ in $x$ equals $(f)$, it is enough is enough to show that  $(f_1,f_2) \subset (f^{p^{r_2}},f_1)$ and that  $(f_1,f_2)\not \subset (f^{p^{r_2}+1},f_1)$. This is the content of claims 3 and 4. 
\smallskip

{\bf Claim 3}  $(f_1,f_2) \subset (f^{p^{r_2}},f_1)$.
\smallskip

We have $(f_2)\subset (f^{p^{r_2}},p)$, hence $(f_1,f_2)\subset (f^{p^{r_2}},f_1,p)$. The latter ideal equals by  claim 2 the ideal $(f^{p^{r_2}},f^{p^{r_1}},f_1)= (f^{p^{r_2}},f_1)$ since $r_2 < r_1$. Hence $ (f_1,f_2)\subset (f^{p^{r_2}},f_1)$ as claimed.
\smallskip

{\bf Claim 4} $(f_1,f_2)\not \subset (f^{p^{r_2}+1},f_1)$.
\smallskip

Suppose $(f_1,f_2) \subset (f^{p^{r_2}+1},f_1)$. Then $(f_1,f_2,p) \subset (f^{p^{r_2}+1},f_1,p)$. Since $r_2 < r_1$ we have $(f_1,f_2,p)=(f^{p^{r_2}},p)$ and $(f^{p^{r_2}+1},f_1,p)=  (f^{p^{r_2}+1},p)$. Thus we get 
$(f^{p^{r_2}},p)\subset (f^{p^{r_2}+1},p)$, a contradiction which confirms the claim.

Combining claims 3 and 4 ends the proof in the case $r_2 < r_1$.
\smallskip

{\bf Second case} $r_2 > r_1$. We may assume that no horizontal component of $\mathcal{D}(j_1) \cap \mathcal{D}(j_2)$ passes through $x.$ Thus, by the first case, the ideal of  $\mathcal{D}(j_1) \cap \mathcal{D}(j_2)$ is 
of the form $(f_1,f_2)=(f_2,g^{p^{r_1}})$ where $g$ is an element of $R$ for which the ideal of $\mathcal{V}(\Lambda)$ in $R$ is given by $(p,g)=(g,f_2)$ and for which  $(g^{p^{r_2}},f_2)=(g^{p^{r_2}},p)=(f_2,p).$ (We can obtain such an element $g$ in the same way as $f$.) 
Thus $p \in (f_1,f_2),$ hence $(f_1,f_2)=(p,f_2,g^{p^{r_1}})=(p,g^{p^{r_1}})=(p,f_1)=(f^{p^{r_1}},f_1)$ by claim 2. 
\qed
\newline

Let $j$ be  a special homomorphism of nonnegative valuation. Suppose that $\mathcal{D}(j)$ contains  $\mathcal{V}(\Lambda)$ with multiplicity $p^r$ for some $r \geq 0.$ Let $l\geq 1.$
Then the same proof shows the following
\begin{Cor} \label{cor1} 
The divisor $\mathcal{D}(j)\cap \mathcal{D}(p^lj)$  in $\mathcal{D}(j)$ contains the curve $\mathcal{V}(\Lambda)$ also with multiplicity $p^r$. \qed
\end{Cor} 
From claims 1 and 2 in the proof we also obtain the following
\begin{Cor}\label{coruglmult}
In the case $r_1=r_2$  the divisor $\mathcal{D}(j_1)\cap \mathcal{D}(j_2)$ in  $\mathcal{D}(j_1)$ contains  $\mathcal{V}(\Lambda)$ with multiplicity at least $p^{\min\{r_1,r_2\}}$. \qed
\end{Cor}

\begin{Pro}\label{glmult}
In the situation described before Proposition \ref{unglmult}, suppose further  that $j_1 \perp j_2$ and $r_1=r_2=:r$. Let $a_1$ resp. $a_2$ be the valuations of $j_1$ resp. $j_2.$ Suppose further that  $a_1$ and $a_2$ are not both even and $a_1 < a_2$.
 Then $a_1$ is odd and the divisor $\mathcal{D}(j_1)\cap \mathcal{D}(j_2)$ in $\mathcal{D}(j_1)$ contains $\mathcal{V}(\Lambda)$ with multiplicity $\frac{a_1+1-r}{2}\cdot p^r$. If $a_2$ is also odd, then  the divisor $\mathcal{D}(j_1)\cap \mathcal{D}(j_2)$ in $\mathcal{D}(j_2)$ also contains $\mathcal{V}(\Lambda)$with multiplicity $\frac{a_1+1-r}{2}\cdot p^r$. 
\end{Pro}

\emph{Proof.}
Using Proposition \ref{redschnitt} it is easy to see that $a_1$ is odd.

 To prove the claim on the multiplicity of $\mathcal{V}(\Lambda)$ in $\mathcal{D}(j_1)\cap \mathcal{D}(j_2)$, we use induction on $r$, starting with the case $r=0$.
We choose a superspecial point $x\in  \mathcal{V}(\Lambda)(\F)$ such that there is no other curve   $\mathcal{V}(\Lambda^{'})$ passing $x$ which belongs to $\mathcal{Z}(j_1)$ or $\mathcal{Z}(j_2)$. We choose a special homomorphism $j_0$  with $h(j_0,j_0)=1$ and $h(j_0,j_1)\neq 0$ and  such that $x \in \mathcal{Z}(j_0)(\F).$ One easily sees that we may write $j_1=\alpha j_0 + \beta \tilde{j_0}$ for some  special homomorphism $\tilde{j_0}\perp j_0$  with $h(\tilde{j_0},\tilde{j_0})=1$ (so that $x \in \mathcal{Z}(\tilde{j_0})(\F)$) and some $\alpha, \beta \in \Zpq^{\times}.$ By perhaps multiplying $j_1$ by a unit we may assume that $\alpha=1.$
We find a special homomorphism $y$ of valuation $1$ such that $y \perp j_0, \tilde{j_0}.$ Then all curves $\mathcal{V}(\Lambda^{'})$ passing through $x$ belong to $\mathcal{Z}(y).$ Since $j_1 \perp j_2$ (and after  perhaps multiplying $j_2$ by a unit), we may write 
$j_2= j_0 - \bar{\beta}^{-1}\tilde{j_0}+\gamma y.$ It follows that $\gamma \equiv 0 \mod p^{(a_1-1)/2}$ but $\gamma \not\equiv 0 \mod p^{(a_1+1)/2}$
Thus the intersection multiplicity $(\mathcal{Z}(j_0),\mathcal{Z}(j_1),\mathcal{Z}(j_2))$ equals  the intersection multiplicity $(\mathcal{Z}(j_0),\mathcal{Z}(\tilde{j_0}),\mathcal{Z}(p^{(a_1-1)/2}y)$ which is $\frac{a_1+1}{2}$ by \cite{KR1}, Theorem 5.1. On the other hand, we have  $(\mathcal{Z}(j_0),\mathcal{Z}(j_1),\mathcal{Z}(y))=(\mathcal{Z}(j_0),\mathcal{Z}(\tilde{j_0}),\mathcal{Z}(y))=1$ which shows that the intersection
  multiplicity of $  \mathcal{V}(\Lambda)$ and   $\mathcal{D}(j_1)\cap \mathcal{Z}(j_0)$ (as divisors in $\mathcal{D}(j_1)$)
is $1$. We may interpret   $(\mathcal{Z}(j_0),\mathcal{Z}(j_1),\mathcal{Z}(j_2))$ as the intersection multiplicity of $\mathcal{D}(j_1) \cap \mathcal{Z}(j_0)$ with $\mathcal{D}(j_1)\cap \mathcal{D}(j_2)$ as divisors in $\mathcal{D}(j_1).$ Hence we see that the multiplicity of  
  $\mathcal{V}(\Lambda)$ in $\mathcal{D}(j_1)\cap \mathcal{D}(j_2)$ is $\leq \frac{a_1+1}{2}$, and it is equal to $ \frac{a_1+1}{2}$ if and only if there is no horizontal component of $\mathcal{D}(j_1)\cap \mathcal{D}(j_2)$ passing through $x$.

Thus it remains to show that there is no such horizontal component.
Assume there is a horizontal component $h$  of $\mathcal{Z}(j_1)\cap \mathcal{Z}(j_2)$ passing through $x.$
Since $j_1$ has valuation $a_1$ we have $\nu_p(1+\beta \bar{\beta})=a_1,$ hence  $\nu_p(\beta + \bar{\beta}^{-1})=a_1.$
We have $j_1-j_2= (\beta+ \bar{\beta}^{-1})\tilde{j_0}-\gamma y.$ This is divisble by $p^{(a_1-1)/2}$ (as a homomorphism from $\overline{\mathbb{Y}}$ to the $p$-divisible group belonging to $x$) and has valuation $a_1.$ Now $\mathcal{Z}(\frac{j_1-j_2}{p^{(a_1-1)/2}})\cap \mathcal{Z}(j_1)$ contains $ \mathcal{V}(\Lambda)$, and Lemma \ref{pmultlem} shows that $\mathcal{Z}(\frac{j_1-j_2}{p^{(a_1-1)/2}})\cap \mathcal{Z}(j_1)$ also contains the horizontal component $h$. But the $\Zpq$-span of $\frac{j_1-j_2}{p^{(a_1-1)/2}}, j_0,j_1$ in $\mathbb{V}$ is the same as the 
 $\Zpq$-span of $j_0,\tilde{j_0},y.$ Hence by \cite{KR1}, Theorem 5.1, the intersection multiplicity  $(\mathcal{Z}(\frac{j_1-j_2}{p^{(a_1-1)/2}}),\mathcal{Z}(j_0),\mathcal{Z}(j_1))$  is $1$. We can regard it as the intersection multiplicity of $\mathcal{Z}(\frac{j_1-j_2}{p^{(a_1-1)/2}}) \cap \mathcal{D}(j_1)$ and $\mathcal{D}(j_1)\cap \mathcal{D}(j_0)$ in $\mathcal{D}(j_1).$ But since $h \subset \mathcal{Z}(\frac{j_1-j_2}{p^{(a_1-1)/2}}) \cap \mathcal{D}(j_1)$ and $ \mathcal{V}(\Lambda) \subset \mathcal{Z}(\frac{j_1-j_2}{p^{(a_1-1)/2}}) \cap \mathcal{D}(j_1)$ the assumption that $h \neq 0$ implies that    $(\mathcal{Z}(\frac{j_1-j_2}{p^{(a_1-1)/2}}),\mathcal{Z}(j_0),\mathcal{Z}(j_1)) \geq 2.$ Thus there is  no horizontal component of $\mathcal{D}(j_1)\cap \mathcal{D}(j_2)$ passing through $x$. This confirms the claim of the proposition in the case $r=0.$

Now we come to the induction step from $r-2$ to $r$.

If  $j_1$ and $j_2$ were replaced by  $j_1/p^{r/2}$ and $j_2/p^{r/2},$ then we were in the situation of the induction start. Now
 we 
choose an $\F$-valued point $x$ of $\mathcal{V}(\Lambda)$ for $j_1/p^{r/2}$ and $j_2/p^{r/2}$ as we did in the induction start.
 We also choose $j_0, \tilde{j_0}$ and $y$ for   $j_1/p^{r/2}$ and $j_2/p^{r/2}$ as before. 
The $\Zpq$-span of $ j_0,j_1,j_2$ in $\mathbb{V}$ is the same as the 
 $\Zpq$-span of $j_0,p^{r/2}\tilde{j_0},p^{(a_1-1)/2}y.$ Using \cite{KR1}, Theorem 5.1, one checks that the intersection multiplicity $(\mathcal{Z}(j_0),\mathcal{D}(j_1),\mathcal{D}(j_2))$ equals $\frac{a_1+1-r}{2}p^r+(\frac{a_1+1-r}{2}+1)\cdot p^{r-1}.$
 By \cite{KR1}, Theorem 5.1, the intersection multiplicity $(\mathcal{Z}(j_0),\mathcal{D}(j_1),\mathcal{Z}(y))$ is $p+1.$  This intersection multiplicity may be considered as the intersection multiplicity of $\mathcal{D}(j_1) \cap \mathcal{Z}(j_0)$ with $\mathcal{D}(j_1)\cap \mathcal{Z}(y)$ as divisors in $\mathcal{D}(j_1).$ Thus the  intersection multiplicity  (in $\mathcal{D}(j_1)$) of $\mathcal{D}(j_1) \cap \mathcal{Z}(j_0)$ with any of the $p+1$ curves $\mathcal{V}(\Lambda^{'})$ passing through $x$ is $1.$ 
 By the induction hypothesis, any of the $p$  curves $\mathcal{V}(\Lambda^{'}) \neq \mathcal{V}(\Lambda)$ passing through $x$ is contained in  $\mathcal{D}(j_1)\cap \mathcal{D}(j_2)$ with multiplicity $\frac{a_1+1-r+2}{2}p^{r-2}.$ Thus the multiplicity of $\mathcal{V}(\Lambda)$ in $\mathcal{D}(j_1)\cap \mathcal{D}(j_2)$ is at most  $\frac{a_1+1-r}{2}p^r+(\frac{a_1+1-r}{2}+1)\cdot p^{r-1}-p\cdot \frac{a_1+1-r+2}{2}p^{r-2}=\frac{a_1+1-r}{2}p^r$, and this is the precise multiplicity if and only if there is no horizontal component of $\mathcal{D}(j_1)\cap \mathcal{D}(j_2)$ passing through $x$.
Thus it remains again to show that there is no such horizontal component.

Assume there is a horizontal component in $\mathcal{D}(j_1)\cap \mathcal{Z}(j_2)$ passing through $x.$ Then it follows as above from Lemma \ref{pmultlem} that this horizontal component also belongs to $\mathcal{Z}(p^{r/2}(j_1-j_2)/p^{(a_1-1)/2})\cap \mathcal{D}(j_1).$ The $\Zpq$-span of $p^{r/2}\frac{j_1-j_2}{p^{(a_1-1)/2}}, j_0,j_1$ in $\mathbb{V}$ is the same as the 
 $\Zpq$-span of $j_0,p^{r/2}\tilde{j_0},p^{r/2}y$. Thus by  \cite{KR1}, Theorem 5.1, the intersection multiplicity  $(\mathcal{Z}(p^{r/2}\frac{j_1-j_2}{p^{(a_1-1)/2}}), \mathcal{Z}(j_0),\mathcal{D}(j_1))$ is $(1+p+...+p^r)+p^{r-1}.$ By  Corollary \ref{coruglmult}, the intersection $\mathcal{Z}(p^{r/2}\frac{j_1-j_2}{p^{(a_1-1)/2}})\cap \mathcal{D}(j_1)$ contains  any of the $p$  curves $\mathcal{V}(\Lambda^{'}) \neq \mathcal{V}(\Lambda)$ passing through $x$ with multiplicity at least $(1+p^2+...+p^{r-2})+p^{r-2}$.  Further, it contains $\mathcal{V}(\Lambda)$ with multiplicity at least $1+p^2+...+p^r$. Thus the intersection multiplicity of the part of  $\mathcal{Z}(p^{r/2}\frac{j_1-j_2}{p^{(a_1-1)/2}})\cap \mathcal{D}(j_1)$ which has support in the special fiber with $\mathcal{Z}(j_0)\cap \mathcal{D}(j_1)$ is at least $p\cdot(1+p^2+...+p^{r-2})+p\cdot p^{r-2}+(1+p^2+...+p^{r})=(1+p+...+p^r)+p^{r-1}$ which is already  the intersection multiplicity  $(\mathcal{Z}(p^{r/2}\frac{j_1-j_2}{p^{(a_1-1)/2}}), \mathcal{Z}(j_0),\mathcal{D}(j_1)).$  Thus there is no  horizontal component in $\mathcal{D}(j_1)\cap \mathcal{Z}(j_2)$ passing through $x.$ 
 The last assertion is proved in the same way.
 \qed

\subsection{Horizontal components of $\mathcal{D}(j_1)\cap \mathcal{D}(j_2)$}

Let $j_1$ be a special homomorphism of odd valuation $a_1,$ let $j_2$ be a  special homomorphism of even  valuation $a_2$ which is perpendicular to $j_1.$  Suppose that $a_1 > a_2.$ Let $y_1=j_1/p^{(a_1-1)/2}$. It has valuation $1.$ Let $y_2=j_2/p^{a_2/2}$. It has valuation $0$. Let $x$ be the unique (superspecial) point in $\mathcal{Z}(y_2)(\F).$ We denote by $\mathcal{V}(\Lambda_1),...,\mathcal{V}(\Lambda_{p+1})$ the $p+1$ irreducible supersingular  curves passing through $x$. They all belong to $\mathcal{Z}(y_1).$ We find a special homomorphism $y_0 \perp y_1,y_2$ of valuation $0$ so that $x \in \mathcal{Z}(y_0)(\F).$ It is easy to see that we find elements $\varepsilon_1,...,\varepsilon_{p+1}\in \Zpq^{\times}$ such that, for each $i \in \{1,...,p+1 \}$, the special cycle  $\mathcal{Z}(y_0+ \varepsilon_i y_2)$ contains $\mathcal{V}(\Lambda_i)$ but does not contain  $\mathcal{V}(\Lambda_j)$ for $j\neq i.$
Let $R= \mathcal{O}_{\mathcal{N},x}$ and let $f_1$ resp. $f_2$ be generators of the ideals of $\mathcal{D}(j_1)$ resp. $\mathcal{D}(j_2)$ in $R$. Further, let $g_i$ for $i \in \{1,...,p+1 \}$ be a generator of the ideal of  $\mathcal{Z}(y_0+ \varepsilon_i y_2)$ in $R.$

{\bf Claim } {\em The ideal of  $\mathcal{V}(\Lambda_i)$ in $R$ is given by $(f_1,g_i).$}  

Since $\mathcal{V}(\Lambda_i) \subset \mathcal{D}(j_1)\cap \mathcal{Z}(y_0+ \varepsilon_i y_2)$ it is enough to show that multiplicity of  $\mathcal{V}(\Lambda_i)$ in  the divisor $\mathcal{D}(j_1)\cap \mathcal{Z}(y_0+ \varepsilon_i y_2)$ in $\mathcal{D}(j_1)$ is $\leq 1$ and that there are no other components passing through $x$. This follows from the fact that  the intersection multiplicity $(\mathcal{Z}(y_0), \mathcal{Z}(y_0+ \varepsilon_i y_2), \mathcal{D}(j_1))$ is $1$. (To see this, observe that 
the $\Zpq$-span of $y_0, y_0+ \varepsilon_i y_2, j_1$ in $\mathbb{V}$  is the $\Zpq$-span of $y_0,y_2, j_1$ and use \cite{KR1}, Theorem 5.1.)

Using Proposition \ref{unglmult}, it follows that, for $a_2 \geq 2$, we may write $$f_2= (g_1\cdot...\cdot g_{p+1})^{p^{a_2-2}}\cdot h+\rho f_1$$ for some $\rho,h \in R$ such that the image of $h$ in $R/(f_1)$ is comprime to $p$. Similarly, for $a_2=0$, we write $f_2= h+\rho f_1.$ 
Let $t \in R$ be an element for which its image $\overline{t}$ in $\overline{R}=R/(p)$ is a generator of the ideal of $\mathcal{Z}(y_2)_p$ in $\overline{R}$ (comp. Theorem \ref{spezf}). Let $\beta = p^{a_2-2}(p^2-1)$ if $a_2 \geq 2$, and let $\beta = 1$ if $a_2=0.$
\begin{Lem}\label{horiz1}
For any $i \in \{1,...,p+1 \}$,  the ideals $I_1=(f_1,g_i,h)$ and $I_2=(t^{\beta}, g_i, p)$ in $R$ are equal. The length (over $W$) of $R/I_1$ is $\beta.$
\end{Lem}
{\em Proof.}
We show that the ideals $I^{'}_1=(f_1, g_1\cdot...\cdot g_{p+1}, h)$ and $I^{'}_2=(t^{\beta}, g_1\cdot...\cdot g_{p+1}, p)$ are equal. Then it follows that $I_1=I_2.$ 

By the claim above we have $p \in  (f_1,g_1\cdot...\cdot g_{p+1})$ and by Theorem \ref{spezf} $f_1 \in (g_1\cdot...\cdot g_{p+1},p)$. Thus  $I^{'}_1=(f_1, g_1\cdot...\cdot g_{p+1}, h,p)$ and $I^{'}_2=(f_1, t^{\beta}, g_1\cdot...\cdot g_{p+1}, p)$.
Thus it is enough to show that the images $\overline{I^{'}_1}$ of $I^{'}_1$ and $\overline{I^{'}_2}$ of $I^{'}_2$ in $\overline{R}$ are equal. For any $r \in R$, we denote by $\overline{r}$ its image in $\overline{R}.$
First assume that  $a_2 \geq 2.$ 
 It follows from Theorem \ref{spezf} that  we may write 
$$f_2=(g_1\cdot...\cdot g_{p+1})^{p^{a_2-2}}\cdot t^{\beta}+ \sigma p$$ for some $\sigma \in R.$
Thus we get 
$$
\overline{f_2}=(\overline{g_1}\cdot...\cdot \overline{g_{p+1}})^{p^{a_2-2}}\cdot \overline{h} + \overline{\rho}\overline{f_1}= (\overline{g_1}\cdot...\cdot \overline{g_{p+1}})^{p^{a_2-2}}\cdot \overline{t}^{\beta}
$$
Further, by Theorem \ref{spezf}, $$\overline{f_1}=(\overline{g_1}\cdot...\cdot \overline{g_{p+1}})^{p^{a_1-1}}.$$
Combining both equations (and using $a_1>a_2$) we get 
$$
 \overline{t}^{\beta}=\overline{h}+\overline{\rho}(\overline{g_1}\cdot...\cdot \overline{g_{p+1}})^{p^{a_1-1}-p^{a_2-2}}.
$$
Since $p^{a_1-1}-p^{a_2-2}>0$, it follows that $\overline{t}^{\beta} \in \overline{I^{'}_1}$ and  $\overline{h} \in \overline{I^{'}_2}.$ Thus indeed  $\overline{I^{'}_2} = \overline{I^{'}_2},$ hence  $I^{'}_1= I^{'}_2$  and hence $I_1=I_2.$
It follows from Theorem \ref{spezf} that the length of $R/I_2$ and hence also the length of $R/I_1$ is $\beta.$ 
The proof in the case $a_2=0$  is the same, one just has to replace the expression  
$(g_1\cdot...\cdot g_{p+1})^{p^{a_2-2}}$ in the proof above by $1.$ 
\qed
\begin{Cor}\label{Corhoriz1}
In the situation described above there is a horizontal component of $\mathcal{D}(j_1)\cap \mathcal{D}(j_2)$ passing through $x$. Its intersection multiplicity (in $\mathcal{D}(j_1)$) with each of the curves   $\mathcal{V}(\Lambda_i)$  is $\beta.$
 \qed
\end{Cor}

\begin{Pro}\label{horiz2}
Let $j_1,j_2$ be perpendicular special homomorphisms of valuations $a_1$ resp. $a_2.$ Suppose that $a_1$ and $a_2$ are not both even and that $1 \leq a_1 < a_2.$
Then, provided that $a_1$ is odd, there is no horizontal component in $\mathcal{D}(j_1) \cap \mathcal{D}(j_2).$ If $a_1$ is even, then the horizontal component described in Corollary \ref{Corhoriz1} is the only horizontal component in $\mathcal{D}(j_1)\cap \mathcal{D}(j_2)$.
\end{Pro}
We will prove this in the next section.

\section{Calculation of intersection multiplicities}

We note that, by Lemma \ref{Lweg}, for any triple of linearly independent special homomorphisms $j_1,j_2,j_3$, we can identify  $\mathcal{O}_{\mathcal{Z}(j_1)}\otimes^{\mathbb{L}}_{\mathcal{O}_{\mathcal{N}}}\mathcal{O}_{\mathcal{Z}(j_2)}\otimes^{\mathbb{L}}_{\mathcal{O}_{\mathcal{N}}}\mathcal{O}_{\mathcal{Z}(j_3)}$ with  $(\mathcal{O}_{\mathcal{Z}(j_1)}\otimes_{\mathcal{O}_{\mathcal{N}}}\mathcal{O}_{\mathcal{Z}(j_2)})\otimes^{\mathbb{L}}_{\mathcal{O}_{\mathcal{Z}(j_2)}}(\mathcal{O}_{\mathcal{Z}(j_3)}\otimes_{\mathcal{O}_{\mathcal{N}}}\mathcal{O}_{\mathcal{Z}(j_2)})=\mathcal{O}_{\mathcal{Z}(j_1) \cap \mathcal{Z}(j_2)}\otimes^{\mathbb{L}}_{\mathcal{O}_{\mathcal{Z}(j_2)}}\mathcal{O}_{\mathcal{Z}(j_3)\cap \mathcal{Z}(j_2)}.$ Here, of course, one can replace the role of $\mathcal{Z}(j_2)$ by $\mathcal{Z}(j_1)$ or $\mathcal{Z}(j_3)$, and one can replace one or several $\mathcal{Z}(j_i)$ by $\mathcal{D}(j_i).$ If $j$ is odd, then $\mathcal{D}(j)$ is a regular formal scheme whose 
 special fibre is a scheme and its underlying reduced subscheme is a union of copies of the Fermat curve. The same reasoning as in \cite{D}  now shows that the intersection number of two divisors $E=E_1+E_2$ and $F$ in $\mathcal{D}(j)$ (defined as usual as the Euler-Poincar\'e characteristic of their structure sheaves) is bilinear i.e. satisfies $(E,F)=(E_1,F)+(E_2,F)$, provided that the support of  $E \cap F$ is contained in $\mathcal{D}(j)_p$ and is proper over $\F$.

Next we compute the intersection multiplicity of two curves $\mathcal{V}(\Lambda)$ and $\mathcal{V}(\Lambda^{'})$ in $\mathcal{D}(j)_p$.

\begin{Lem}\label{selsch}
Let $j \in \mathbb{V}$ be of   odd valuation $u$ and
let $\mathcal{V}(\Lambda),\mathcal{V}(\Lambda^{'})\subset \mathcal{D}(j)_p$. Then their intersection multiplicity in $\mathcal{D}(j)$ has the value:
$$
(\mathcal{V}(\Lambda),\mathcal{V}(\Lambda^{'}))=
\begin{cases}
0, &   \ \text{if }\mathcal{V}(\Lambda) \text{ and }\mathcal{V}(\Lambda{'}) \text{ do not intersect} \\
1, &   \ \text{if } \mathcal{V}(\Lambda) \text{ and }\mathcal{V}(\Lambda{'}) \text{ intersect in precisely one point} \\
-p^2-p, & \ \text{if } \Lambda= \Lambda^{'} \text{ and } u=1 \\ 
-2p^2-p+1, & \ \text{if } \Lambda= \Lambda^{'} \text{ and } u>1 \text{ and }  \mathcal{V}(\Lambda)\subset \mathcal{D}(j/p)_p \\
-p^2-p+1, &  \ \text{if }\Lambda= \Lambda^{'} \text{ and } u>1 \text{ and }  \mathcal{V}(\Lambda)\not \subset \mathcal{D}(j/p)_p \\
\end{cases}
$$
\end{Lem}
\emph{Proof.}
The first point is clear.
Suppose we are in the situation of point 2 and denote the intersection point of  $\mathcal{V}(\Lambda)$ and $\mathcal{V}(\Lambda^{'})$ by $x.$ We choose perpendicular  special homomorphisms $j_0, j_1$ of valuation $0$ which are perpendicular to $j$ and for which $x \in \mathcal{Z}(j_0)(\F)$ and  $x \in \mathcal{Z}(j_1)(\F)$. We find elements $\varepsilon, \varepsilon^{'}$ such that  $\mathcal{V}(\Lambda) \subset \mathcal{Z}(j_0 +\varepsilon j_1)$
 and   $\mathcal{V}(\Lambda^{'}) \subset \mathcal{Z}(j_0 +\varepsilon^{'} j_1)$, and such that these are the only curves in $ \mathcal{Z}(j_0 +\varepsilon j_1)$ resp $ \mathcal{Z}(j_0 +\varepsilon^{'} j_1)$ passing through $x$. Then the intersection multiplicity $ (\mathcal{Z}(j_0 +\varepsilon j_1),  \mathcal{Z}(j_0 +\varepsilon^{'} j_1), \mathcal{D}(j))$ is the same as the intersection multiplicity $(\mathcal{Z}(j_0),\mathcal{Z}(j_1),\mathcal{D}(j))$, which is $1$ by \cite{KR1}, Theorem 5.1. From this claim follows.
To compute the intersection multiplicity in the remaining cases, we use the fact that $(\mathcal{V}(\Lambda), \mathcal{D}(j)_p)=0$. This follows as in the proof of Proposition \ref{multlin} from the exact sequence
$$
\begin{CD}
0 @ >>> \mathcal{O}_{\mathcal{D}(j)} @> p \cdot>> \mathcal{O}_{\mathcal{D}(j)} @>>>\mathcal{O}_{\mathcal{D}(j)}/(p)@>>> 0.
\end{CD}
$$

Suppose we are in the third case. 
Then $\mathcal{V}(\Lambda)$ has multiplicity $1$ in $\mathcal{D}(j)_p$ and there are $p\cdot (p+1)$ further supersingular curves which intersect  $\mathcal{V}(\Lambda)$. All have multiplicity $1$ in $\mathcal{D}(j)_p.$ Claim 2 in the proof of Proposition \ref{unglmult} shows that the multiplicity of $\mathcal{V}(\Lambda)$ in the divisor  $\mathcal{D}(j)_p$ in  $\mathcal{D}(j)$ is the same as the multiplicity of  $\mathcal{V}(\Lambda)$ in the divisor  $\mathcal{D}(j)_p$ in $\mathcal{N}_p$. Thus we can rewrite  $(\mathcal{V}(\Lambda), \mathcal{D}(j)_p)=0$ as 
$(\mathcal{V}(\Lambda),\mathcal{V}(\Lambda))+p(p+1)=0.$ This shows the claim.

Suppose we are in the fourth case and suppose that $\mathcal{V}(\Lambda)$ belongs to $\mathcal{Z}(j/p^{(u-1)/2}).$
Then $\mathcal{V}(\Lambda)$ has multiplicity $p^{u-1}$ in $\mathcal{D}(j)_p$ and there are $p\cdot (p+1)$ further  supersingular curves which intersect  $\mathcal{V}(\Lambda)$ and have  multiplicity $p^{u-1}$ in $\mathcal{D}(j)_p.$ 
There are $p\cdot (p^3-p)$ further  supersingular curves which intersect $\mathcal{V}(\Lambda)$ and have  multiplicity $p^{u-3}$ in $\mathcal{D}(j)_p.$ 
Thus we can rewrite  $(\mathcal{V}(\Lambda), \mathcal{D}(j)_p)=0$ as 
$p^{u-1}(\mathcal{V}(\Lambda),\mathcal{V}(\Lambda))+p^{u-1}\cdot p(p+1)+p(p^3-p)p^{u-3}=0$. Thus $(\mathcal{V}(\Lambda),\mathcal{V}(\Lambda))=-2p^2-p+1.$

Suppose we are in the fourth case and suppose that $\mathcal{V}(\Lambda)$ does not belong to $\mathcal{Z}(j/p^{(u-1)/2}).$ Suppose the multiplicity of $\mathcal{V}(\Lambda)$ in $\mathcal{D}(j)_p$ is $p^{2r}.$ Then there are $p-1$  further  supersingular curves which intersect  $\mathcal{V}(\Lambda)$ and have  multiplicity $p^{2r}$ in $\mathcal{D}(j)_p.$ There are $p\cdot p^3$  further  supersingular curves which intersect  $\mathcal{V}(\Lambda)$ and have  multiplicity $p^{2r-2}$. There is one   further  supersingular curve which intersects $\mathcal{V}(\Lambda)$ and has  multiplicity $p^{2r+2}$ in $\mathcal{D}(j)_p.$ Thus we can rewrite  $(\mathcal{V}(\Lambda), \mathcal{D}(j)_p)=0$ as 
$p^{2r}(\mathcal{V}(\Lambda),\mathcal{V}(\Lambda))+p^{2r}\cdot (p-1)+p\cdot p^3\cdot p^{2r-2}+p^{2r+2}=0$. Thus $(\mathcal{V}(\Lambda),\mathcal{V}(\Lambda))=-2p^2-p+1.$

In the last case it follows that  the multiplicity of $\mathcal{V}(\Lambda)$ in $\mathcal{D}(j)_p$ is $1.$  There are $p-1$  further  supersingular curves which intersect  $\mathcal{V}(\Lambda)$ and have  multiplicity $1$ in $\mathcal{D}(j)_p.$
There is one   further  supersingular curve which intersects $\mathcal{V}(\Lambda)$ and has  multiplicity $p^{2}$ in $\mathcal{D}(j)_p.$
It follows as before that $(\mathcal{V}(\Lambda),\mathcal{V}(\Lambda))=-p^2-p+1.$ \qed
\newline

We introduce the following terminologies. 
\begin{Def}
Let $j$ be a special homomorphism of valuation $a\geq 0$. If $a \geq 2$, then the boundary of $\mathcal{D}(j)$ (and of $\mathcal{Z}(j)$) is the set $B(j)$ which consists of all $\mathcal{V}(\Lambda)$ (where $\Lambda \in \mathcal{L}^{(3)}$) for which 
$\mathcal{V}(\Lambda) \subset \mathcal{Z}(j)$ but $\mathcal{V}(\Lambda) \not\subset \mathcal{Z}(j/p).$ If $a \in \{0,1\}$, then $B(j)$ is defined to be empty.

Write $j=p^r y$, where $y$ is a special homomorphism of valuation $0$ or $1$. Then the center $C(j)$ is defined as follows. If $j$ is even, then it simply consists of the unique $\F$-valued point of $\mathcal{Z}(y).$ If $j$ is odd, then it consists of all $\F$-valued points of $\mathcal{Z}(y)$ which are the intersection point of two (and hence $p+1$) curves  $\mathcal{V}(\Lambda), \mathcal{V}(\Lambda^{'}) \subset \mathcal{Z}(y).$

Let $j^{'}$ be a further  special homomorphism which is linearly independent of $j$. Then  the part of $(\mathcal{D}(j)\cap \mathcal{D}(j^{'}))$ (as divisor in $\mathcal{D}(j)$) which has support in the special fiber is denoted by $(\mathcal{D}(j)\cap \mathcal{D}(j^{'}))_v$. The horizontal part is denoted by $(\mathcal{D}(j)\cap \mathcal{D}(j^{'}))_h.$
\end{Def}

\begin{Lem}\label{nullsch}
Let $j_1,j_2$ be perpendicular special homomorphisms of  valuations $a_1 \neq a_2$. Suppose that $j_1$  is odd. Suppose the curve $\mathcal{V}(\Lambda)$ lies in $\mathcal{D}(j_1)\cap \mathcal{D}(j_2)$ but does not belong to $B(j_2)$ and also does not contain any point in $ C(j_2) $. Then the intersection multiplicity $(\mathcal{V}(\Lambda),(\mathcal{D}(j_1)\cap \mathcal{D}(j_2))_v )$ in $\mathcal{D}(j_1)$ is zero.
\end{Lem}
{\em Proof.} This follows from Propositions \ref{unglmult}  and \ref{glmult} together with Lemma \ref{selsch}.
Suppose, for example, that $\mathcal{V}(\Lambda)$ has multiplicities $p^{r_1}$ resp. $p^{r_2}$ in $\mathcal{D}(j_1)$ resp. $\mathcal{D}(j_2)$ where $r_1 > r_2>0$. Suppose further  that $\mathcal{V}(\Lambda)$ does not contain any points in $C(j_1).$
It follows that  $(\mathcal{V}(\Lambda),(\mathcal{D}(j_1)\cap \mathcal{D}(j_2))_v )= (-2p^2-p+1)p^{r_2}+p^3\cdot p \cdot p^{r_2-2}+(p-1)p^{r_2}+p^{r_2+2}=0.$
Here the first summand comes from the self intersection of $\mathcal{V}(\Lambda)$, the second comes from the intersections in the $p^3$ supersingular points on $\mathcal{V}(\Lambda)$ for which all $p$ further  supersingular curves passing through the corresponding point have multiplicity $p^{r_2-2}$ in $\mathcal{D}(j_2)_p.$ The third summand comes from the single supersingular curve of multiplicity $p^{r_2+2}$ which intersects $\mathcal{V}(\Lambda)$ and the last summand comes from the remaining $p-1$ supersingular curves of multiplicity $p^{r_2}$ intersecting $\mathcal{V}(\Lambda)$ (in the same point as the single  supersingular curve of multiplicity $p^{r_2+2}$).
 The other cases are checked analogously. \qed
\newline

{\em Proof of Proposition \ref{horiz2} }

Suppose first that $a_1$ and $a_2$ are odd. We find a third special homomorphism $j_3$ which is perpendicular to $j_1,j_2$ and has odd valuation $a_3$ so large that, for any    $\mathcal{V}(\Lambda) \subset (\mathcal{Z}(j_1)\cap \mathcal{Z}(j_2))_{\red}$, the multiplicity of  $\mathcal{V}(\Lambda)$ in $\mathcal{D}(j_3)_p$ is greater than the  multiplicity of  $\mathcal{V}(\Lambda)$ in $\mathcal{D}(j_1)_p$ and $\mathcal{D}(j_2)_p.$ In particular,  $(\mathcal{Z}(j_1)\cap \mathcal{Z}(j_2))_{\red} \subset \mathcal{Z}(j_3)_{\red}$ (compare Proposition \ref{redschnitt}).  Lemma \ref{pmultlem}  implies that there is no horizontal component in $\mathcal{D}(j_1)\cap \mathcal{D}(j_3)$ and in  $\mathcal{D}(j_2)\cap \mathcal{D}(j_3).$ By Proposition  \ref{unglmult} we know the vertical parts of  $\mathcal{D}(j_1)\cap \mathcal{D}(j_3)$ and  $\mathcal{D}(j_2)\cap \mathcal{D}(j_3).$ Thus we can compute the intersection multiplicity $(\mathcal{D}(j_1),\mathcal{D}(j_2),\mathcal{D}(j_3))= ((\mathcal{D}(j_1)\cap \mathcal{D}(j_3)),(\mathcal{D}(j_2)\cap \mathcal{D}(j_3)))$ where the latter is meant as intersection multiplicity in $\mathcal{D}(j_3).$ By Propositions \ref{unglmult}  and \ref{glmult} we also know $(\mathcal{D}(j_1)\cap \mathcal{D}(j_2))_v$. 
Thus we can also compute $((\mathcal{D}(j_1)\cap \mathcal{D}(j_2))_v,(\mathcal{D}(j_1)\cap \mathcal{D}(j_3)))$ (intersection multiplicity in $\mathcal{D}(j_1)$). Now the claim that  $(\mathcal{D}(j_1)\cap \mathcal{D}(j_2))_v= (\mathcal{D}(j_1)\cap \mathcal{D}(j_2))$ follows from the following 
\smallskip

{\bf Claim} $ ((\mathcal{D}(j_1)\cap \mathcal{D}(j_3)),(\mathcal{D}(j_2)\cap \mathcal{D}(j_3)))=((\mathcal{D}(j_1)\cap \mathcal{D}(j_2))_v,(\mathcal{D}(j_1)\cap \mathcal{D}(j_3))).$
\smallskip

We show that both intersection multiplicities are $0$.

It follows from the choice of $j_3$  and Lemma \ref{pmultlem} that   $((\mathcal{D}(j_1)\cap \mathcal{D}(j_2))_v,(\mathcal{D}(j_1)\cap \mathcal{D}(j_3)))=((\mathcal{D}(j_1)\cap \mathcal{D}(j_2))_v,\mathcal{D}(j_1)_p)$ which is $0$  since $((\mathcal{D}(j_1)\cap \mathcal{D}(j_2))_v$ has support in the special fiber and this support is proper over $\F$.

Next we compute  $ ((\mathcal{D}(j_1)\cap \mathcal{D}(j_3)),(\mathcal{D}(j_2)\cap \mathcal{D}(j_3))).$ For any  $\mathcal{V}(\Lambda)$ in $\mathcal{D}(j_1) \cap \mathcal{D}(j_3)$, we compute the intersection multiplicity   $ (\mathcal{V}(\Lambda),(\mathcal{D}(j_2)\cap \mathcal{D}(j_3)))$ and then  add all these with the corresponding multiplicity of $\mathcal{V}(\Lambda)$ in  $\mathcal{D}(j_1) \cap \mathcal{D}(j_3)$.
We assume first that $a_1>1.$
Then, for any $\mathcal{V}(\Lambda)$ in $\mathcal{D}(j_2) \cap \mathcal{D}(j_3)$ that does not belong to $B(j_2)$, we have    $ (\mathcal{V}(\Lambda),(\mathcal{D}(j_2)\cap \mathcal{D}(j_3)))=0.$ (This is a slightly stronger claim than the claim of Lemma \ref{nullsch} in this case but it is proved in the same way.)  
Thus we only need to consider those  $\mathcal{V}(\Lambda)$ in $\mathcal{D}(j_1) \cap \mathcal{D}(j_3)$ which  belong to $B(j_2)$ and those  $\mathcal{V}(\Lambda)$ in  $\mathcal{D}(j_1)\cap \mathcal{D}(j_3)$ which do not belong to $\mathcal{D}(j_2)$ but intersect one curve in $\mathcal{D}(j_2)\cap \mathcal{D}(j_3).$
Let   $\mathcal{V}(\Lambda) \subset \mathcal{D}(j_1) \cap \mathcal{D}(j_3)$  belong to $B(j_2)$ and suppose that the multiplicity of $\mathcal{V}(\Lambda)$ in $\mathcal{D}(j_1)_p$ is $p^r$, where $1<r<a_1-1.$ Then $(p^r\mathcal{V}(\Lambda), \mathcal{D}(j_2)\cap \mathcal{D}(j_3))=p^r(-2p^2-p+1+(p-1)+p^2)=-p^{r+2}$. Now there are $p^3\cdot p$ further  curves in $\mathcal{D}(j_1)\cap \mathcal{D}(j_3)$ which do not belong to $\mathcal{D}(j_2)$ but intersect $\mathcal{V}(\Lambda)$. Each has multiplicity $p^{r-2}$ in $\mathcal{D}(j_1)_p.$ The sum of their contribution to the intersection multiplicity is $p^3\cdot p^{r-2}=p^{r+2}$, which cancels with  $(p^r\mathcal{V}(\Lambda), \mathcal{D}(j_2)\cap \mathcal{D}(j_3))=p^r(-2p^2-p+1+(p-1)+p^2)=-p^{r+2}$.

Next we consider the curves    $\mathcal{V}(\Lambda) \subset \mathcal{D}(j_1) \cap \mathcal{D}(j_3)$ which  belong to $B(j_2)$ and whose multiplicity in $\mathcal{D}(j_1)_p$ is $1.$ Each contributes with  $(\mathcal{V}(\Lambda), \mathcal{D}(j_2)\cap \mathcal{D}(j_3))=-2p^2-p+1+p-1+p^2=-p^2$ to the intersection multiplicity, and one easily checks that there are $(p\cdot p)^{(a_2-a_1)/2}(p^3-p)(p^3\cdot p)^{(a_1-1)/2-1}$ such curves. The total contribution of all is $-p^{a_2-a_1-2}(p^2-1).$ 

Next we consider the curves    $\mathcal{V}(\Lambda) \subset \mathcal{D}(j_1) \cap \mathcal{D}(j_3)$ which  belong to $B(j_2)$ and whose multiplicity in $\mathcal{D}(j_1)_p$ is $p^{a_1-1}.$ There are $(p\cdot p)^{(a_2-1)/2}$ such curves . We have  again $p^{(a_1-1)/2}(\mathcal{V}(\Lambda), \mathcal{D}(j_2)\cap \mathcal{D}(j_3))$ $=p^{(a_1-1)}(-p^2).$ Each such curve is intersected by  $(p^3-p)\cdot p$ further  curves in $\mathcal{D}(j_1)\cap \mathcal{D}(j_3)$ which do not belong to $\mathcal{D}(j_2)$ but intersect $\mathcal{V}(\Lambda)$ and have multiplicity $p^{a_1-3}$ in $\mathcal{D}(j_1)_p.$ There are  $(p\cdot p)$ further  curves in $\mathcal{D}(j_1)\cap \mathcal{D}(j_3)$ which do not belong to $\mathcal{D}(j_2)$ but intersect $\mathcal{V}(\Lambda)$ and have multiplicity $p^{a_1-1}$ in $\mathcal{D}(j_1)_p.$ Adding everything up we get a total contribution of $(p\cdot p)^{(a_2-1)/2}(p^{a_1-1}(-p^2)+(p^3-p)\cdot p\cdot p^{a_1-3}+p^{a_1-1}p^2)=p^{a_2-a_1-2}(p^2-1).$
This contribution cancels the contribution coming  from the curves    $\mathcal{V}(\Lambda) \subset \mathcal{D}(j_1) \cap \mathcal{D}(j_3)$ which  belong to $B(j_2)$ and whose multiplicity in $\mathcal{D}(j_1)_p$ is $1.$

Thus alltogether we get  $ ((\mathcal{D}(j_1)\cap \mathcal{D}(j_3)),(\mathcal{D}(j_2)\cap \mathcal{D}(j_3)))=0.$ The case $a_1=1$ is computed analogously.
\newline

Now we consider the case that $a_1$ is even and hence $a_2$ is odd. 
In this case the claim follows from Lemma \ref{horiz1} and Lemma \ref{pmultlem}.
\newline

Now we consider the case that $a_1$ is odd and $a_2$ is even. 
 We find a third special homomorphism $\tilde{j_3}$ which is perpendicular to $j_1,j_2$ and has even valuation $\tilde{a_3}$ so large that, for any   $\mathcal{V}(\Lambda) \subset (\mathcal{Z}(j_1)\cap \mathcal{Z}(j_2))_{\red}$, the multiplicity of  $\mathcal{V}(\Lambda)$ in $\mathcal{D}(\tilde{j_3})_p$ is greater than the  multiplicity of  $\mathcal{V}(\Lambda)$ in $\mathcal{D}(j_1)_p$ and $\mathcal{D}(j_2)_p.$ In particular  $(\mathcal{Z}(j_1)\cap \mathcal{Z}(j_2))_{\red} \subset \mathcal{Z}(\tilde{j_3})_{\red}$. We find a unit $\varepsilon \in \Zpq^{\times}$ such that $j_3:= \tilde{j_3}+\varepsilon p^{(\tilde{a_3}-a_2)/2}j_2$ has valuation $a_3=\tilde{a_3}+1.$ Then still, for any   $\mathcal{V}(\Lambda) \subset (\mathcal{Z}(j_1)\cap \mathcal{Z}(j_2))_{\red}$, the multiplicity of  $\mathcal{V}(\Lambda)$ in $\mathcal{D}(j_3)_p$ is greater than the  multiplicity of  $\mathcal{V}(\Lambda)$in $\mathcal{D}(j_1)_p$ and $\mathcal{D}(j_2)_p$ and $(\mathcal{Z}(j_1)\cap \mathcal{Z}(j_2))_{\red} $ are contained in $ \mathcal{Z}(j_3)_{\red}$.
There is one curve $\mathcal{V}(\Lambda)$ which contains the unique point $x \in C(j_2)$ and has multiplicity $p^{a_3-1}$ in $\mathcal{D}(j_3)$. The other $p$  supersingular curves which contain $x$ have multiplicity $p^{a_3-3}$ in $\mathcal{D}(j_3)$.
Note that $(\mathcal{D}(j_1),\mathcal{D}(j_2),\mathcal{D}(\tilde{j_3}))=(\mathcal{D}(j_1),\mathcal{D}(j_2),\mathcal{D}(j_3)) $ and that $\mathcal{D}(j_2)\cap \mathcal{D}(\tilde{j_3})=\mathcal{D}(j_2)\cap \mathcal{D}(j_3).$ By Lemma \ref{pmultlem},  the latter contains no horizontal component passing through any $\F$-valued point $\neq x$.

\smallskip
{\bf Claim} {\em There is a horizontal component $h_{23}$ in $\mathcal{D}(j_2)\cap \mathcal{D}(j_3).$ Its intersection multiplicity (in $\mathcal{D}(j_3)$) with each supersingular curve passing through $x$ is $p^{a_2-2}(p^2-1).$}
\smallskip

This is proved in the same way as Lemma \ref{horiz1}.

We already know that $\mathcal{D}(j_1)\cap \mathcal{D}(j_3)$ contains no horizontal components.
As in the first case  $(\mathcal{D}(j_1),\mathcal{D}(j_2),\mathcal{D}(j_3)) = ((\mathcal{D}(j_1)\cap \mathcal{D}(j_3)),(\mathcal{D}(j_2)\cap \mathcal{D}(j_3)))=((\mathcal{D}(j_1)\cap \mathcal{D}(j_2)),(\mathcal{D}(j_1)\cap \mathcal{D}(j_3)))$, where the second term is meant as intersection multiplicity in $\mathcal{D}(j_3)$ and the third  term is meant as intersection multiplicity in $\mathcal{D}(j_1)$. The claim that there is no horizontal component in $\mathcal{D}(j_1)\cap \mathcal{D}(j_2)$ follows from the following
\smallskip

{\bf Claim} $ ((\mathcal{D}(j_1)\cap \mathcal{D}(j_3)),(\mathcal{D}(j_2)\cap \mathcal{D}(j_3)))=((\mathcal{D}(j_1)\cap \mathcal{D}(j_2))_v,(\mathcal{D}(j_1)\cap \mathcal{D}(j_3))).$
\smallskip

We show that both expressions are equal to $0.$

First we calculate  $ ((\mathcal{D}(j_1)\cap \mathcal{D}(j_3)),(\mathcal{D}(j_2)\cap \mathcal{D}(j_3))).$
 For any  $\mathcal{V}(\Lambda)$ in $\mathcal{D}(j_2) \cap \mathcal{D}(j_3)$, we compute the intersection multiplicity   $ (\mathcal{V}(\Lambda),(\mathcal{D}(j_1)\cap \mathcal{D}(j_3)))$ and then add all these with the corresponding multiplicity of $\mathcal{V}(\Lambda)$ in  $\mathcal{D}(j_2) \cap \mathcal{D}(j_3)$.
We assume first that $a_1>1.$
Then as above, for any $\mathcal{V}(\Lambda)$ in $\mathcal{D}(j_1) \cap \mathcal{D}(j_3)$ that does not belong to $B(j_1)$, we have    $ (\mathcal{V}(\Lambda),(\mathcal{D}(j_1)\cap \mathcal{D}(j_3)))=0.$ 
Thus we only need to consider those  $\mathcal{V}(\Lambda)$ in $\mathcal{D}(j_2) \cap \mathcal{D}(j_3)$ which  belong to $B(j_1)$ and those  $\mathcal{V}(\Lambda)$ in  $\mathcal{D}(j_2)\cap \mathcal{D}(j_3)$ which do not belong to $\mathcal{D}(j_1)$ but intersect one curve in $\mathcal{D}(j_1)\cap \mathcal{D}(j_3).$ 
Suppose that  $\mathcal{V}(\Lambda) \subset \mathcal{D}(j_2) \cap \mathcal{D}(j_3)$  belongs to $B(j_1)$. Then one checks that in the case that its multiplicity in $\mathcal{D}(j_2)_p$ is greater than $1$, its contribution to the intersection multiplicity cancels with all contributions coming from those supersingular curves which intersect   $\mathcal{V}(\Lambda)$ but do not belong to $\mathcal{D}(j_1)$.
Suppose now that  $\mathcal{V}(\Lambda) \subset \mathcal{D}(j_2) \cap \mathcal{D}(j_3)$  belongs to $B(j_1)$ and to $B(j_2)$. Then its contribution to the intersection multiplicity is $-p^2$ (as above.) There are $(p+1)(p\cdot p)^{(a_2-a_1-1)/2}(p^3-p)p\cdot (p^3 \cdot p)^{(a_1-1)/2-1}$ such curves. Together their conribution is $-(p+1)(p^2-1)p^{a_1+a_2-3}.$

Finally, the contribution of the contribution coming from the horizontal component in $\mathcal{D}(j_2)\cap \mathcal{D}(j_3)$ is  $((\mathcal{D}(j_1)\cap \mathcal{D}(j_3)),h_{23})= (p+1)\cdot p^{a_1-1}\cdot p^{a_2-2}(p^2-1).$ Thus  $((\mathcal{D}(j_1)\cap \mathcal{D}(j_3)),(\mathcal{D}(j_2)\cap \mathcal{D}(j_3)))=0$
The case $a_1=1$ is computed similarly.

Finally,  it is easy to see that  $((\mathcal{D}(j_1)\cap \mathcal{D}(j_2))_v,(\mathcal{D}(j_1)\cap \mathcal{D}(j_3)))=0$, since there is no curve  $\mathcal{V}(\Lambda)$ in $\mathcal{D}(j_1) \cap \mathcal{D}(j_2)$ which belongs to $B(j_3).$ \qed

\begin{Pro} \label{alleung} Let $j_1,j_2,j_3$ be odd special homomorphisms of valuations $1\leq a_1\leq a_2 \leq a_3$ which are perpendicular to each other.
\begin{enumerate}
\item If $a_1<a_2<a_3$,  then $(\mathcal{D}(j_1),\mathcal{D}(j_2),\mathcal{D}(j_3))=0.$
\item If $a_1<a_2=a_3$, then $(\mathcal{D}(j_1),\mathcal{D}(j_2),\mathcal{D}(j_3))=-(\frac{a_1+1}{2}p^2-\frac{a_1-1}{2})(p+1)p^{a_1+a_2-3}.$
\item If $a_1=a_2<a_3$,  then $(\mathcal{D}(j_1),\mathcal{D}(j_2),\mathcal{D}(j_3))=(p+1)p^{2a_1-2}.$
\item If $a_1=a_2=a_3=:a$, then $$(\mathcal{D}(j_1),\mathcal{Z}(j_2),\mathcal{Z}(j_3))-(\mathcal{D}(j_1),\mathcal{Z}(j_2/p),\mathcal{Z}(j_3/p))=-(\frac{a+1}{2}p-\frac{a+3}{2})(p+1)p^{2a-2}.$$
\end{enumerate}
\end{Pro}

In the proof of this and in the proof of the following proposition there occur sums which might be empty depending on $a_1,a_2,a_3$. In such cases these sums are meant to be zero. There may also be expressions of the form $p^{d},$ where $d$ is of the form $d=\frac{a_i-l}{k}$ for some integers $l,k$. If $l >a_i$ these expressions must be replaced by zero.

{\em Proof.} In the first case was proved in the proof of Proposition \ref{horiz2}

In the second case we compute $((\mathcal{D}(j_1)\cap \mathcal{D}(j_2))(\mathcal{D}(j_1)\cap \mathcal{D}(j_3)))$ as intersection in $\mathcal{D}(j_1).$  For any  $\mathcal{V}(\Lambda)$ in $\mathcal{D}(j_1) \cap \mathcal{D}(j_2)$, we compute the intersection multiplicity   $ (\mathcal{V}(\Lambda),(\mathcal{D}(j_1)\cap \mathcal{D}(j_3)))$ and then add all these with the corresponding multiplicity of $\mathcal{V}(\Lambda)$ in  $\mathcal{D}(j_1) \cap \mathcal{D}(j_2)$.
We assume first that $a_1>1.$ 
We only need to consider $\mathcal{V}(\Lambda)\in B(j_3).$  Suppose that the multiplicity of   $\mathcal{V}(\Lambda)\in B(j_3)$   in $\mathcal{D}(j_1)$ is $1$. Then the intersection multiplicity $ (\mathcal{V}(\Lambda),(\mathcal{D}(j_1)\cap \mathcal{D}(j_3)))$
is $\frac{a_1+1}{2}(-p^2-p+1)+p^2\frac{a_1-1}{2}+(p-1)\frac{a_1+1}{2}=-p^2.$ Further, the multiplicity of  $\mathcal{V}(\Lambda)$ in  $\mathcal{D}(j_1) \cap \mathcal{D}(j_2)$ is $\frac{a_1+1}{2}.$ There are $(p+1)p(p\cdot p)^{(a_2-a_1)/2-1}\cdot (p^3-p)p(p^3\cdot p)^{(a_1-3)/2}$ such  $\mathcal{V}(\Lambda)$. Thus the total contribution of these curves is $-(p+1)p^{a_1+a_2-3}(\frac{a_1+1}{2}p^2-\frac{a_1+1}{2}).$
 Consider now  the curves   $\mathcal{V}(\Lambda)\in B(j_3)$  whose multiplicity in $\mathcal{D}(j_1)$ is $p^{2r}$, where $0<2r<(a_1-1).$ Then one checks similarly that  the contribution coming from one such curve is $-p^2$ and
that the total contribution coming from all such curves is $-(p+1)p^{a_1+a_2-3-2r}(p^2-1).$ 
 Consider now  the curves   $\mathcal{V}(\Lambda)\in B(j_3)$  whose multiplicity in $\mathcal{D}(j_1)$ is $p^{a_1-1}.$ Then one checks similarly that  the contribution coming from one such curve is $-p^2$ and
that the total contribution coming from all such curves is $-(p+1)p^{a_2}.$ Adding all these contributions we get the claimed intersection multiplicity. The case $a_1=1$ is computed similarly.

Let us consider the third case. We  compute $((\mathcal{D}(j_1)\cap \mathcal{D}(j_3))(\mathcal{D}(j_2)\cap \mathcal{D}(j_3)))$ as intersection in $\mathcal{D}(j_3).$  Again we assume $a_1>1$. For any  $\mathcal{V}(\Lambda)$ in $\mathcal{D}(j_1) \cap \mathcal{D}(j_3)$, we compute the intersection multiplicity   $ (\mathcal{V}(\Lambda),(\mathcal{D}(j_2)\cap \mathcal{D}(j_3)))$ and then add all these with the corresponding multiplicity of $\mathcal{V}(\Lambda)$ in  $\mathcal{D}(j_1) \cap \mathcal{D}(j_3)$. We only need to consider such curves which are in $B(j_2)$ or which do not belong to $\mathcal{D}(j_2) \cap \mathcal{D}(j_3)$ but intersect one curve lying in $B(j_2).$ First, we consider such curves which are in $B(j_1)\cap B(j_2).$ Each contributes with $-p^2$ to intersection multiplicity and there are $(p^3+1-2(p+1))p\cdot (p^3 \cdot p)^{(a_1-3)/2}.$ Suppose now that the multiplicity $p^r$ of $ \mathcal{V}(\Lambda)$ in $\mathcal{D}(j_1)$ is greater than $1$ but less than $p^{a_1-1}$ such curves. Suppose further  that it is not equal to the multiplicity of $ \mathcal{V}(\Lambda)$ in $\mathcal{D}(j_3)_p$ (which is $p^{a_3-a_1}$.) Then one easily checks that the contribution of  $ \mathcal{V}(\Lambda)$ cancels with the contributions of the supersingular curves  which do not belong to $\mathcal{D}(j_2) \cap \mathcal{D}(j_3)$ but intersect $ \mathcal{V}(\Lambda)$. In the case that $a_3<2a_1-1$ the same is true if $r=a_1-1.$ Suppose that $r=a_1-1$ and  $a_3 \geq 2a_1-1.$  Then one easily checks that the contribution coming from $ \mathcal{V}(\Lambda)$ together with the contributions coming from the curves  which do not belong to $\mathcal{D}(j_2) \cap \mathcal{D}(j_3)$ but intersect $\mathcal{V}(\Lambda)$ is $p^{a_1-1}(p^2-1).$ Now there are $(p+1)p^{(a_1-3)/2}$ such curves. Further, the case the case  $r=a_3-a_1$ can only happen if $r=a_1-1$ which is the case just discussed. Adding everything we get the desired contribution in the case  $a_3 \geq 2a_1-1.$
Suppose now that     $a_3 < 2a_1-1.$ Then $r=a_1-1$ cannot occur. Suppose that  $r=a_3-a_1.$ Then the contribution coming from  $ \mathcal{V}(\Lambda)$ is $\frac{a_1+1-r}{2}p^{r}(-p^2)$.  Each curve  which do not belong to $\mathcal{D}(j_2) \cap \mathcal{D}(j_3)$ but intersects $\mathcal{V}(\Lambda)$ contributes with $\frac{a_1+1-r+2}{2}p^{r-2}$ and there are $p^3\cdot p$ such curves. Adding all their contributions to the contribution of $\mathcal{V}(\Lambda)$, we get $p^{r+2}.$ There are $(p+1)p(p^2)^{r/2-1}(p^3-p)p\cdot (p^3 \cdot p)^{(a_1-1-r)/2}$ such curves $\mathcal{V}(\Lambda)$. Their total contribution is $(p+1)p(p^2)^{r/2-1}(p^3-p)p\cdot (p^3 \cdot p)^{(a_1-3-r)/2}p^{r+2}.$ Adding everything gives again the desired result. The case $a_1=1$ is obvious.

Let us consider the fourth case. Here we are faced with the problem that the propositions in section 3 do not give us enough information about the structure of $\mathcal{D}(j_i) \cap \mathcal{D}(j_k).$ We proceed as follows. We find elements $\varepsilon, \eta \in \Zpq^{\times}$ such that $y_2:= j_2+\varepsilon j_3$ and  $y_2:= j_2+\eta j_3$ both have valuation $a+1$ and such that the $\Zpq$-span of $j_2,j_3$ in $\mathbb{V}$ is the same as the $\Zpq$-span of $y_2,y_3$ in $\mathbb{V}.$ Thus $(\mathcal{D}(j_1),\mathcal{Z}(j_2),\mathcal{Z}(j_3))=(\mathcal{D}(j_1),\mathcal{Z}(y_2),\mathcal{Z}(y_3))$ and  $(\mathcal{D}(j_1),\mathcal{Z}(j_2/p),\mathcal{Z}(j_3/p))=(\mathcal{D}(j_1),\mathcal{Z}(y_2/p),\mathcal{Z}(y_3/p))$. Using Proposition \ref{diag}, we obtain the relation 
$(\mathcal{D}(j_1),\mathcal{Z}(j_2),\mathcal{Z}(j_3))-(\mathcal{D}(j_1),\mathcal{Z}(j_2/p),\mathcal{Z}(j_3/p)=(\mathcal{D}(j_1),\mathcal{D}(y_2),\mathcal{Z}(y_3/p))+(\mathcal{D}(j_1),\mathcal{Z}(y_2/p),\mathcal{D}(y_3))+(\mathcal{D}(j_1),\mathcal{D}(y_2),\mathcal{D}(y_3)).$
Note that the point in $C(y_2)$ and the point in $C(y_3)$ both lie on the curve which belongs to $\mathcal{Z}(j_1/p^{(a-1)/2}) \cap \mathcal{Z}(j_1/p^{(a-1)/2}) \cap \mathcal{Z}(j_1/p^{(a-1)/2}).$ We assume that $a \geq 5$ for a consistent notation, but the cases $a=1$ and $a=3$ are easily checked in the same way.
We calculate $(\mathcal{D}(j_1),\mathcal{D}(y_2),\mathcal{Z}(y_3/p))=((\mathcal{D}(j_1) \cap \mathcal{D}(y_2))(\mathcal{D}(j_1)\cap \mathcal{Z}(y_3/p)))$ as intersection multiplicity in $\mathcal{D}(j_1).$ We write $(\mathcal{D}(j_1)\cap \mathcal{Z}(y_3/p)))=(\mathcal{D}(j_1)\cap \mathcal{Z}(y_3/p))_v+(\mathcal{D}(j_1)\cap \mathcal{Z}(y_3/p))_h$ for its decomposition in its vertical and horizontal part. Note that, by Proposition \ref{horiz2}, there is no horizontal part in $\mathcal{D}(j_1)\cap \mathcal{D}(y_2)$ and that any curve $\mathcal{V}(\Lambda)$ in $(\mathcal{D}(j_1)\cap \mathcal{Z}(y_3/p))_v$ also belongs to  $\mathcal{D}(j_1)\cap \mathcal{D}(y_2)$. To calculate 
$((\mathcal{D}(j_1 \cap \mathcal{D}(y_2)),(\mathcal{D}(j_1)\cap \mathcal{Z}(y_3/p))_v)$  we compute, for any  $\mathcal{V}(\Lambda)$ in $\mathcal{D}(j_1) \cap \mathcal{Z}(y_3/p)$,  the intersection multiplicity   $ (\mathcal{V}(\Lambda),(\mathcal{D}(j_1)\cap \mathcal{D}(y_2)))$ and then add all these with the corresponding multiplicity of $\mathcal{V}(\Lambda)$ in  $\mathcal{D}(j_1) \cap \mathcal{D}(y_3/p)$. We only need to consider such curves which are in $B(j_2).$ 
One easily sees that there are $p\cdot (p^3\cdot p)^{(a-3)/2}$ such curves and each contributes with the value $-p^2$ (as above). Thus $((\mathcal{D}(j_1 \cap \mathcal{D}(y_2)),(\mathcal{D}(j_1)\cap \mathcal{Z}(y_3/p))_v)=-p^{2a-3}.$
Next it follows from Proposition \ref{horiz2} that  $((\mathcal{D}(j_1) \cap \mathcal{D}(y_2)),(\mathcal{D}(j_1)\cap \mathcal{Z}(y_3/p))_h)=(p^{a-1}+p\cdot p^{a-3})p^{a-1}.$ Altogether we get $((\mathcal{D}(j_1) \cap \mathcal{D}(y_2))(\mathcal{D}(j_1)\cap \mathcal{Z}(y_3/p)))=p^{2a-2}.$ 
Of course the intersection multiplicity $(\mathcal{D}(j_1),\mathcal{Z}(y_2/p),\mathcal{D}(y_3))$ is calculated in the same way and has the same value.

Now we compute $(\mathcal{D}(j_1),\mathcal{D}(y_2),\mathcal{D}(y_3))=((\mathcal{D}(j_1) \cap \mathcal{D}(y_2))(\mathcal{D}(j_1)\cap \mathcal{D}(y_3)))$ as intersection multiplicity in $\mathcal{D}(j_1).$ By Proposition \ref{horiz2}, both factors have support in the special fiber. For any  $\mathcal{V}(\Lambda)$ in $\mathcal{D}(j_1) \cap \mathcal{D}(y_3)$, we compute the intersection multiplicity   $ (\mathcal{V}(\Lambda),(\mathcal{D}(j_1)\cap \mathcal{D}(y_2)))$ and then add all these with the corresponding multiplicity of $\mathcal{V}(\Lambda)$ in  $\mathcal{D}(j_1) \cap \mathcal{D}(y_3)$. We only need to consider such curves which are in $B(y_2)$ or which do not belong to $\mathcal{D}(j_1) \cap \mathcal{D}(y_2)$ but intersect one curve lying in $B(y_2).$ First we consider such curves which are in $B(y_2)\cap B(y_3) \cap B(j_1).$ One easily checks as above  that each contributes with the value $-p^2$  and has multiplicity $\frac{a+1}{2}$ in $\mathcal{D}(j_1)\cap \mathcal{D}(y_3).$  There are $(p^3-p)\cdot p (p^3\cdot p)^{(a-3)/2}$ such curves. Their total contribution thus is  $\frac{a+1}{2}(p^3-p)\cdot p \cdot (p^4)^{(a-3)/2}(-p^2).$

Now we consider such curves which are in $B(y_2)\cap B(y_3)$ but not in $B(j_1).$ Again each contributes with $(-p^2)$ and there are $(p-1)p\cdot (p^3\cdot p)^{(a-3)/2}$ such curves. Their total contribution thus is  $(p^2-p)\cdot  (p^4)^{(a-3)/2}(-p^2).$

Suppose now that   $\mathcal{V}(\Lambda)$ belongs to $B(y_2)$ but has multiplicity $\geq p^2$ in $\mathcal{D}(j_1)_p$ and $\mathcal{D}(y_3)_p.$ It follows that  the multiplicity of  $\mathcal{V}(\Lambda)$  in $\mathcal{D}(y_3)_p$ is precisely $p^2.$ Then one easily checks that in the case that the multiplicity of  $\mathcal{V}(\Lambda)$  in $\mathcal{D}(j_1)$ is greater than $p^2$ the contribution of $\mathcal{V}(\Lambda)$  cancels with the contributions coming from the curves which do not belong to $\mathcal{D}(y_2)$ but intersect $\mathcal{V}(\Lambda).$   In the case that  the multiplicity of  $\mathcal{V}(\Lambda)$  in $\mathcal{D}(j_1)_p$ is precisely $p^2$ it follows that the contribution coming from  $\mathcal{V}(\Lambda)$  is again $-p^2.$ Further, the multiplicity of   $\mathcal{V}(\Lambda)$ in $\mathcal{D}(y_3)\cap \mathcal{D}(y_1)$ is $\frac{a-1}{2}p^2.$ The total contribution coming from the curves which do not belong to $\mathcal{D}(y_2)$ but intersect $\mathcal{V}(\Lambda)$  is $\frac{a+1}{2}(p^3\cdot p).$ There are $p\cdot (p^3-p)p(p^4)^{(a-5)/2}$ such curves. 

Adding everything belonging to $((\mathcal{D}(j_1) \cap \mathcal{D}(y_2))(\mathcal{D}(j_1)\cap \mathcal{D}(y_3)))$ we get a value of $p^{2a-2}(-p^2\frac{a+1}{2}+\frac{a-1}{2}+p).$ 
Adding now  $(\mathcal{D}(j_1),\mathcal{D}(y_2),\mathcal{Z}(y_3/p))+(\mathcal{D}(j_1),\mathcal{Z}(y_2/p),\mathcal{D}(y_3))=2p^{2a-2}$ we get the desired result. \qed

\begin{Pro} \label{einerung} Let $j_1,j_2,j_3$ be special homomorphisms of valuations $1\leq a_1\leq a_2 \leq a_3$ which are perpendicular to each other and such that precisely one of them is odd.
\begin{enumerate}
\item If $a_2$ is odd, then  $(\mathcal{D}(j_1),\mathcal{D}(j_2),\mathcal{D}(j_3))=(p^2-1)(p+1)p^{a_1+a_2-3}.$
\item  If $a_1=a_2$, then $(\mathcal{Z}(j_1),\mathcal{D}(j_2),\mathcal{D}(j_3))=p^{2a_1-2}(p^2-1).$
\item If $a_1<a_2$ and $a_3$ is odd, then   $(\mathcal{D}(j_1),\mathcal{D}(j_2),\mathcal{D}(j_3))=(p^2-1)(p+1)p^{a_1+a_2-3}.$
\item If $a_1$ is odd and $a_2<a_3$, then  $(\mathcal{D}(j_1),\mathcal{D}(j_2),\mathcal{D}(j_3))=0.$
\item If $a_2=a_3$, then  $(\mathcal{D}(j_1),\mathcal{D}(j_2),\mathcal{D}(j_3))=-(\frac{a_1+1}{2}p^2-\frac{a_1-1}{2})(p+1)p^{a_1+a_2-3}.$
\end{enumerate}
\end{Pro}
{\em Proof.} 
For the first point, we write $(\mathcal{D}(j_1),\mathcal{D}(j_2),\mathcal{D}(j_3))=((\mathcal{D}(j_2)\cap \mathcal{D}(j_1)),(\mathcal{D}(j_2) \cap \mathcal{D}(j_3)))$ as intersection multiplicity in $\mathcal{D}(j_2)$.  With the help of Lemma \ref{nullsch} it is easy to see that $((\mathcal{D}(j_2)\cap \mathcal{D}(j_1))_v,(\mathcal{D}(j_2) \cap \mathcal{D}(j_3))_v)=0.$
By Proposition \ref{horiz2}, there is no horizontal component in $(\mathcal{D}(j_2) \cap \mathcal{D}(j_3))$ and the horizontal part $h$ of $(\mathcal{D}(j_2) \cap \mathcal{D}(j_1))$ is given by Corollary \ref{Corhoriz1}. Its intersection multiplicity with $(\mathcal{D}(j_2) \cap \mathcal{D}(j_3))$ is $(p+1)p^{a_2-1}(p^2-1)p^{a_1-2}$ which is the claimed intersection multiplicity.

Let us consider the second point. We would like to compute  $(\mathcal{Z}(j_1),\mathcal{D}(j_2),\mathcal{D}(j_3))=((\mathcal{Z}(j_1)\cap \mathcal{D}(j_3)),(\mathcal{D}(j_2) \cap \mathcal{D}(j_3)))$. But there are horizontal components in both factors $(\mathcal{Z}(j_1)\cap \mathcal{D}(j_3))$ and $(\mathcal{D}(j_2) \cap \mathcal{D}(j_3))$ which pass through the same $\F$-valued point in $\mathcal{D}(j_3)$ and we do not know what the intersection multiplicity of these two horizontal components is. To avoid this problem, we proceed as follows. We find a unit $\varepsilon \in \Zpq^{\times}$ such that $y_1:= j_1 + \varepsilon j_2$ has valuation $a_2+2.$ The points in $C(y_1)$ and in $C(j_2)$ lie on one supersingular curve.
 Then  $(\mathcal{Z}(j_1),\mathcal{D}(j_2),\mathcal{D}(j_3))=(\mathcal{Z}(y_1),\mathcal{D}(j_2),\mathcal{D}(j_3))=((\mathcal{Z}(y_1)\cap \mathcal{D}(j_3)),(\mathcal{D}(j_2) \cap \mathcal{D}(j_3)))=((\mathcal{Z}(y_1/p^2)\cap \mathcal{D}(j_3)),(\mathcal{D}(j_2) \cap \mathcal{D}(j_3)))+((\mathcal{D}(y_1/p)\cap \mathcal{D}(j_3)),(\mathcal{D}(j_2) \cap \mathcal{D}(j_3)))+((\mathcal{D}(y_1)\cap \mathcal{D}(j_3)),(\mathcal{D}(j_2) \cap \mathcal{D}(j_3)))$. 
Let us assume that $a_3 \geq a_2+3.$  
The remaining cases are proved in the same way, but to keep the notation simpler we make this assumption.
First, we compute $((\mathcal{Z}(y_1/p^2)\cap \mathcal{D}(j_3)),(\mathcal{D}(j_2) \cap \mathcal{D}(j_3)))$. We write $(\mathcal{Z}(y_1/p^2)\cap \mathcal{D}(j_3))=(\mathcal{Z}(y_1/p^2)\cap \mathcal{D}(j_3))_v+(\mathcal{Z}(y_1/p^2)\cap \mathcal{D}(j_3))_h$ for the vertical and horizontal component and analogously for
$(\mathcal{D}(j_2) \cap \mathcal{D}(j_3))$. Then  $((\mathcal{Z}(y_1/p^2)\cap \mathcal{D}(j_3))_h,(\mathcal{D}(j_2) \cap \mathcal{D}(j_3)))_v= p^{a_1-2}(p^{a_1-2}+p\cdot p^{a_1-4}).$ Further,  $((\mathcal{Z}(y_1/p^2)\cap \mathcal{D}(j_3))_v,(\mathcal{D}(j_2) \cap \mathcal{D}(j_3)))_h=(p^{a_1-2})(p^2-1)(1+p^2+...+p^{a_1-4}+p\cdot (1+p^2+...+p^{a_1-6}))$.   Obviously,  $((\mathcal{Z}(y_1/p^2)\cap \mathcal{D}(j_3))_h,(\mathcal{D}(j_2) \cap \mathcal{D}(j_3)))_h=0.$
To  compute  $((\mathcal{Z}(y_1/p^2)\cap \mathcal{D}(j_3))_v,(\mathcal{D}(j_2) \cap \mathcal{D}(j_3)))_v)$, we compute again, for any  $\mathcal{V}(\Lambda)$ in $(\mathcal{Z}(y_1/p^2)\cap \mathcal{D}(j_3)$,  the intersection multiplicity   $ (\mathcal{V}(\Lambda),(\mathcal{D}(j_2)\cap \mathcal{D}(j_3)))$ and then add all these with the corresponding multiplicity of $\mathcal{V}(\Lambda)$ in  $(\mathcal{Z}(y_1/p^2)\cap \mathcal{D}(j_3)$. We only need to consider such curves which are in $B(j_2)$, or which do not belong to $\mathcal{D}(j_2) \cap \mathcal{D}(j_3)$ but intersect one curve lying in $B(j_2)$  or such curves which contain the point in $C(j_2).$ The contribution coming from the curves containing the point in  $C(j_2)$ is $(p^{a_1-2}-p^{a_1})(1+p^2+...+p^{a_1-4}+p\cdot (1+p^2+...+p^{a_1-6}))$. 
The contribution coming from the curves in $B(j_2)$ is $p(p^3\cdot p)^{(a_1-4)/2}(-p^2).$

Similarly adding the summands coming from horizontal resp. vertical parts, one checks that  $((\mathcal{D}(y_1/p)\cap \mathcal{D}(j_3)),(\mathcal{D}(j_2) \cap \mathcal{D}(j_3)))=   p^{a_1-2}(p^2-1)(p^{a_1-2}+p\cdot p^{a_1-4})+p^{a_1-2}(p^2-1)(p^{a_1-2}+p\cdot p^{a_1-4})+p(p^3-1)(p^3\cdot p)^{(a_1-4)/2}(-p^2)+(p^{a_1-2}-p^{a_1})(p^{a_1-2}+p\cdot p^{a_1-4}).$
Similarly one computes $((\mathcal{D}(y_1)\cap \mathcal{D}(j_3)),(\mathcal{D}(j_2) \cap \mathcal{D}(j_3)))=  p^{a_1}(p^2-1)(p^{a_1-2}+p^{a_1-2}(p^2-1)(p^{a_1}+p\cdot p^{a_1-2})+(p^{a_1-2}-p^{a_1})(p^{a_1}-p\cdot p^{a_1-2})+p\cdot (p^3\cdot p)^{(a_1-2)/2}(-p^2).$

Adding everything we get the desired result.

We consider the third point. By the same reasons as in the second case we choose a  unit $\varepsilon \in \Zpq^{\times}$ such that $y_2:= j_1 + \varepsilon j_2$ has valuation $a_2+2$. Then the points in $C(y_1)$ and in $C(j_2)$ lie on one supersingular curve and we have  $(\mathcal{D}(j_1),\mathcal{D}(j_2),\mathcal{D}(j_3)) = (\mathcal{D}(j_1),\mathcal{D}(y_2),\mathcal{D}(j_3))= ((\mathcal{D}(j_1)\cap \mathcal{D}(j_3)),(\mathcal{D}(y_2) \cap \mathcal{D}(j_3))).$ We write $((\mathcal{D}(j_1)\cap \mathcal{D}(j_3))=((\mathcal{D}(j_1)\cap \mathcal{D}(j_3))_v+((\mathcal{D}(j_1)\cap \mathcal{D}(j_3))_h$ and 
$  ((\mathcal{D}(y_2)\cap \mathcal{D}(j_3))=((Dy_2)\cap \mathcal{D}(j_3))_v+((\mathcal{D}(y_2)\cap \mathcal{D}(j_3))_h.$ 
Let us assume that $a_3<a_2+1.$
We have by Proposition \ref{horiz2} $((\mathcal{D}(j_1)\cap \mathcal{D}(j_3))_h,(\mathcal{D}(y_2) \cap \mathcal{D}(j_3))_v)=p^{a_1-2}(p^2-1)(p\cdot p^{a_2-2}+p^{a_2}).$ Obviously,  $((\mathcal{D}(j_1)\cap \mathcal{D}(j_3))_h,(\mathcal{D}(y_2) \cap \mathcal{D}(j_3))_h)=0.$
Further,  $((\mathcal{D}(j_1)\cap \mathcal{D}(j_3))_v,(\mathcal{D}(y_2) \cap \mathcal{D}(j_3))_h)=p^{a_2}(p^2-1)(p\cdot p^{a_1-4}+p^{a_1-2})$. Finally,  one computes the intersection  $((\mathcal{D}(j_1)\cap \mathcal{D}(j_3))_v,(\mathcal{D}(y_2) \cap \mathcal{D}(j_3))_v)$ as usual using Lemma \ref{nullsch}. The result is 
$(p^{a_2-2}-p^{a_2})(p\cdot p^{a_1-2}+p^{a_1})$. Adding everything gives the desired result. 
The proof in the case $a_3=a_2+1$ is similar, in this case again  $((\mathcal{D}(j_1)\cap \mathcal{D}(j_3))_h,(\mathcal{D}(y_2) \cap \mathcal{D}(j_3))_v)=p^{a_1-2}(p^2-1)(p\cdot p^{a_2-2}+p^{a_2}),$ but the other terms are all zero.

The fourth point follows using Lemma \ref{nullsch} and the fact that $a_1<a_2<a_3$ together with the fact that by Proposition \ref{horiz2} there are no horizontal components in $\mathcal{D}(j_1)\cap \mathcal{D}(j_2)$ and $\mathcal{D}(j_1)\cap \mathcal{D}(j_3)$.

The fifth point is proved as the second point in Proposition \ref{alleung}. \qed
\newline

Recall that the {\em fundamental matrix} of three special cycles $j_1,j_2,j_3$ has, by definition,  the entry $h(j_k,j_l)$ at $k,l.$

\begin{The}\label{intmult}
Let $j_1,j_2,j_3$ be special homomorphisms such that $\mathcal{Z}(j_1) \cap \mathcal{Z}(j_2) \cap \mathcal{Z}(j_3)$ is non-empty and such that the fundamental matrix $T$ of $j_1,j_2,j_3$ is non-singular. Suppose that $T$ is $GL_3(\Zpq)$-equivalent to the diagonal matrix $diag(p^{a_1},p^{a_2},p^{a_3})$, where $0\leq a_1 \leq a_2 \leq a_3$. Then the intersection multiplicity $(\mathcal{Z}(j_1),\mathcal{Z}(j_2),\mathcal{Z}(j_3))$ is finite and is given by the formula
$$
(\mathcal{Z}(j_1),\mathcal{Z}(j_2),\mathcal{Z}(j_3))= -\frac{1}{2}\sum_{k=0}^{a_1}\sum_{l=0}^{a_1+a_2-2k}(-1)^k((k+l)p^{2k+l}-(k+l+a_3+1)p^{a_1+a_2-l}).
$$
\end{The}
{\em Proof.} By Proposition \ref{diag}  we may assume that $j_1,j_2,j_3$ are perpendicular to each other and have valuations $a_1,a_2,a_3.$ Now the theorem follows via induction on $a_1+a_2+a_3$ from Theorem 5.1 in \cite{KR1} (which treats the case $a_1=0$) and Propositions \ref{alleung} and \ref{einerung} 
(note that $a_1+a_2+a_3$ is odd). 

\qed

\section{The connection to hermitian representation densities }
For the moment, let $n$ be arbitrary again. For non-singular hermitian matrices $S\in \Herm_m(\Zpq)$ and $T\in \Herm_n(\Zpq),$ where $m \geq n$, the representation density $\alpha_p(S,T)$ is defined as $$\alpha_p(S,T)= \operatorname*{lim}_{l\rightarrow\infty} p^{-ln(2m-n)} \mid \{x \in M_{m,n}(\Zpq/p^l \Zpq);
\ S[x] \equiv T \mod p^l \}\mid,$$
where $S[x]= {^t} xS\sigma(x).$
Note that  $\alpha_p(S,T)$ depends only on the $GL_m(\Zpq)$ resp. $GL_n(\Zpq)$ -equivalence classes of  $S$ and $T$.
Further, let $$\beta_p(S,T)=\operatorname*{lim}_{l\rightarrow\infty} p^{-ln(2m-n)} \mid \{x \in M_{m,n}(\Zpq/p^l \Zpq), \ {\operatorname*{rank}}_{\Zpq / p \Zpq}(x)=n ;
\ S[x] \equiv T \mod p^l \}\mid.$$
For $r \geq 0$, let $S_r= \diag(S, 1_r).$ Then one can show that there is a polynomial $F_p(S,T;X)\in \Q[X]$ such that $\alpha_p(S_r, T)=F_p(S,T; (-p)^{-r}).$ The derivative   $\alpha^{'}_p(S,T)$ is defined to be
$$\alpha^{'}_p(S,T):= -\frac{\partial}{\partial X}F_p(S,T;X)\mid_{X=1}. $$
We denote the unit matrix of size $s$ by $1_s.$
In particular, for any  $T\in \Herm_n(\Zpq),$ there is a polynomial $F_p(X;T)$ such that $\alpha_p(1_s,T)=F_p((-p)^{-s};T).$
\begin{The}\label{repdens}
Suppose that the matrix $T\in \Herm_n(\Zpq)$ is $GL_n(\Zpq)$-equivalent to the diagonal matrix $diag(p^{a_1},...,p^{a_n})$, where  $0 \leq a_1 \leq ... \leq a_n.$ 
\begin{enumerate}
\item If $n=1$, then 
$$F_p(X;T)=(1-X)\sum_{l=0}^{a_1}(pX)^l $$
\item If  $n\geq 2$, let $T^{-}= diag(p^{a_1},...,p^{a_{n-1}}).$ Further, let $H_p(X)=(1-X)(1+pX)$ and let $$A_p(X;T)=(1-(-p)^{n-1}X)(1-(-p)^{-n}X^{-1})(-p)^{a_1+...+a_{n-1}}(-1)^{(n+1)a_n}(p^nX)^{a_n+2}.$$ Then there is the following inductive formula, 
$$
F_p(X;T)=\frac{H_p(X)F_p(p^2X;T^{-})-A_p(X;T)F_p(X;T^{-})}{1-p^{2n}X^2}.
$$
\end{enumerate}
\end{The}
{\em Proof.} 
The first point follows from \cite{N}, Example 2, p. 239.

Let us prove the second point. We define the matrix $T_+=diag(p^{a_1},...,p^{a_{n-1}},p^{a_n+2}).$ 
\smallskip

{\bf Claim} {\em For even $s$ we have} $$\alpha_p(1_s,T_+)=p^{2n-2s}\alpha_p(1_s,T)+\beta_p(1_s,(p^{a_n+2}))\alpha_p(1_{s-2},T^{-}).$$
\smallskip

This claim can be proved in the same way as Theorem 2.6 (1) in \cite{Ka}.

We note that it follows from point 1 and the formula $$\alpha_p(S,(p^a))=\beta_p(S,(p^a))+p^{3-2m}\alpha_p(S,(p^{a-2}))$$ (which can be proved as Proposition 2.1 in \cite{Ka}) that we have $\beta_p(1_s,(p^{a_n+2}))=(1-p^{-s})(1+p^{1-s}).$
Inserting this into the formula given by the claim we obtain the equation
\[
F_p(X;T_+)=p^{2n}X^2F_p(X;T)+H_p(X)F_p(p^2X;T^{-}). \tag{$*$}
\]
Now let $t_n(X)=\prod_{i=0}^{n-1}(1-(-p)^iX)^{-1}$. 
 We obtain from \cite{I}, Corollary 3.2  the following functional equation. 
\[
F_p(X;T)t_n(X)=((-1)^{n+1}(p^nX))^{a_1+...+a_n}F_p(p^{-2n}X^{-1};T)t_n(p^{-2n}X^{-1}) \tag{$**$}.
\]
(Of course the corresponding formulas for $T^{-}$ and $T_+$ also hold).
Next one plugs in $p^{-2n}X^{-1}$ for $X$ in formula $(*)$ and replaces in the resulting equation the expressions of the form  $F_p(... ; ... )$ by the corresponding expressions obtained from formula $(**)$ resp. the corresponding formulas for $T^- $ and $T_+$.  One thus obtains a new equation which expresses $F_p(X;T_+)$ as a new linear combination (where the coefficients are rational functions) of $F_p(X;T)$ and $F_p(p^2X;T^{-})$. Now we equate this linear combination with the right hand side of $(*).$ From this the result follows easily. \qed

\begin{The}\label{Hauptsatz}
Let $j_1,j_2,j_3$ be special homomorphisms such that $\mathcal{Z}(j_1) \cap \mathcal{Z}(j_2) \cap \mathcal{Z}(j_3)$ is non-empty and such that the fundamental matrix $T$ of $j_1,j_2,j_3$ is non-singular. 
 Let $S=1_3$. Then there is the following identity between intersection multiplicities of special cycles and representation densities, 
 $$
 (\mathcal{Z}(j_1),\mathcal{Z}(j_2),\mathcal{Z}(j_3))=\frac{\alpha^{'}_p(S,T)}{\alpha_p(S,S)}.
 $$  
\end{The}
{\em Proof.} 
By our assumptions on $j_1,j_2,j_3$, the matrix $T$ is $GL_3(\Zpq)$-equivalent to a  diagonal matrix of the form $diag(p^{a_1},p^{a_2},p^{a_3})$ where $0\leq a_1 \leq a_2 \leq a_3$.
Using the preceding theorem, one easily checks that 
$$
 F_p(S,T;X)=\frac{(1+\frac{X}{p^3})(1-\frac{X}{p^2})(1+\frac{X}{p})}{1+X}\sum_{k=0}^{a_1}\sum_{l=0}^{a_1+a_2-2k}(-1)^k X^{k+l}(p^{2k+l}-p^{a_1+a_2-l}X^{a_3+1}).
$$ Here we use that $a_1+a_2+a_3$ is odd.
Further, it follows from the theorem that  $\alpha_p(S,S)=(1+\frac{1}{p})(1-\frac{1}{p^2})(1+\frac{1}{p^3})$ (comp. also \cite{KR1}, p. 44). Using this, one easily computes the expression $\frac{\alpha^{'}_p(S,T)}{\alpha_p(S,S)}$ and checks that it is the same as the expression for $ (\mathcal{Z}(j_1),\mathcal{Z}(j_2),\mathcal{Z}(j_3))$ given in Theorem \ref{intmult}. \qed
\newline

\bigskip
\bigskip
\centerline{\bf\Large Part II:  The Global Result }
\bigskip

In this part we will apply the local result in Theorem \ref{Hauptsatz} to relate the intersection multiplicities of global special cycles to Fourier coefficients of derivatives of certain  Eisenstein series. 
\section{The global intersection problem}
We briefly recall the notion of special cycles in the global sense and the corresponding intersection problem as introduced in \cite{KR2}. 
The whole section summarizes contents of \cite{KR2} to which we refer  for more details. 

Let $\bk$ be an imaginary quadratic field and denote by $O_{\bk}$ its ring of integers. Let $n\geq 1$ and let $0 \leq r \leq n$. We consider the Deligne-Mumford-stack $\M(n-r,r)$ over  $O_{\bk}$ whose functor of points associates to any  locally noetherian $O_{\bk}$-scheme $S$ the groupoid 
of tuples $(A,\iota, \lambda),$ where $A$ is an abelian scheme over $S$, where $\iota:  O_{\bk} \rightarrow \End_S(A)$ is an  $O_{\bk}$-action on $A$ satisfying the signature condition $(n-r,r)$, i.e.  the characteristic polynomial of $\iota(a)$ is of the form
$$
\charpol (\iota(a), \Lie A)(T)=(T-a)^{n-r}(T-a^{\sigma})^r \in \mathcal{O}_S [ T ], 
$$
 and where $\lambda: A \rightarrow A^{\vee}$ is a principal polarization for which the Rosati involution satisfies $\iota(a)^*=\iota(a^{\sigma}).$ Further, we require  that that the $O_{\bk}$-action on $\Lie A$ satisfies  $\wedge^{n-r+1}(\iota(a)-a)=0$ and $\wedge^{r+1}(\iota(a)-a^{\sigma})=0.$
Let $$\M=\M(n-r,r)\times_{\Spec O_{\bk}} \M(1,0).$$

As in \cite{KR2}, we denote by  $\mathcal{R}_{(n-r,r)}(\bk)$ the set of relevant hermitian spaces, i.e. the set of isomorphism classes of hermitian spaces $V$ over $\bk$ of dimension $n$ and signature $(n-r,r)$ which contain a self-dual $O_{\bk}$-lattice. We denote by   $\mathcal{R}_{(n-r,r)}(\bk)^{\sharp}$ the set isomorphism classes of pairs $(V, [[L]])$, where $V \in \mathcal{R}_{(n-r,r)}(\bk)$ and $[[L]]$ is a $U(V)$-genus of self-dual hermitian lattices in $V$ (see \cite{KR2}).   Note that, for odd $n$, there is only one such genus, whereas, for $n$ even, the number of genera is $2$ or $1$ depending on whether or not the conditions of \cite{KR2}, Proposition 2.14 are satisfied.
We write $\M[\frac{1}{2}]=\M\times_{\Spec O_{\bk}}{\Spec O_{\bk}[\frac{1}{2}]}$. 
As shown in \cite{KR2}, Propositions 2.12 and 2.19, there are natural decompositions
$$
\M=\coprod_{V \in \mathcal{R}_{(n-r,r)}(\bk)}\M^V \ \ \text{ and } \ \ \M[\frac{1}{2}]=\coprod_{V^{\sharp} \in \mathcal{R}_{(n-r,r)}(\bk)^{\sharp}}\M[\frac{1}{2}]^{V^{\sharp}}.
$$ 
For  $V^{\sharp} = (V,L)$, one can identify $\M[\frac{1}{2}]^{V^{\sharp}}\times_{{\Spec O_{\bk}[\frac{1}{2}]}} \Spec \bk$ with (the canonical model of) the Shimura variety for $GU(V)$ with respect to the stabilizer $K$ of $L$ in $GU(V)(\A_f)$, if this stabilizer is small enough. See \cite{KR2}, section 4, in particular Proposition 4.4, for the precise statement. 

Now let $p>2$ be a prime which is inert in $\bk.$
We denote by $\widehat{\M}(n-r,r)^{ss}$ the completion of $\M(n-r,r)\times_{O_{\bk}}W(\overline{\F}_p)$ along its  supersingular locus.
For $V^{\sharp} = (V,L) \in \mathcal{R}_{(n-r,r)}(\bk)^{\sharp}$, we denote by  $\widehat{\M}(n-r,r)^{V^{\sharp}, ss}$ the open and closed sublocus where the rational Tate module $T^p(A)^0$  is isomorphic to $V \otimes \A_f^p$ and where the type of the hermitian lattice $T_2(A)$ coincides with the type of  the U(V)-genus $[[L]].$ (See \cite{KR2}.) Let $G^V=GU(V)$, and  let $K^p$ be the stabilizer of $L$ in $G^V(\A_f^p)^0,$ where $$G^V(\A_f^p)^0=\{g\in G^V(\A_f^p)\  |\  \nu(g)\in (\widehat{\Z}^p)^{\times}\}.$$ Let $\mathcal{N}$ be as in the introduction  and  in part I but with signature condition $(n-r,r)$. Then there is an isomorphism of formal algebraic stacks over $W$,
$$
I^V(\Q)\setminus (\mathcal{N}\times G^V(\A_f^p)^0/K^p) \cong \widehat{\M}(n-r,r)^{V^{\sharp}, ss},
$$ where $I^V(\Q)$ denotes the group of quasi-isogenies of some (fixed) $(A^o, \iota^o, \lambda^o)\in \widehat{\M}(n-r,r)^{V^{\sharp}, ss}(\overline{\F}_p)$ that respect $\iota^o$ and $\lambda^o$, see \cite{KR2}, Theorem 5.5.
Combining the uniformizations for $\M(n-1,1)$ and $\M(1,0)$, we obtain a uniformization for $\M$,
$$ 
\widehat{\M}^{{\tilde{V}}^{\sharp},ss}= \coprod_{(V^{\sharp},V_0)}\widehat{\M}^{({V^{\sharp}},V_0),ss},
$$
where $\widehat{\M}^{{\tilde{V}}^{\sharp},ss}$ is the completion  of $\M[\frac{1}{2}]^{\tilde{V}^{\sharp}}\times_{{\Spec O_{\bk}[\frac{1}{2}]}}W(\overline{\F}_p)$ along the supersingular locus, further   
$$
\widehat{\M}^{({V^{\sharp}},V_0),ss}= \widehat{\M}(n-r,r)^{V^{\sharp}, ss}\times \widehat{\M_0}^{V_0,ss}, 
$$
and where the paris $({V^{\sharp}},V_0)$ run over $\mathcal{R}_{(n-r,r)}(\bk)^{\sharp}\times \mathcal{R}_{(1,0)}(\bk) $ such that $\Hom(V_0,V)\cong\tilde{V}.$

Given a point in $\M(S)$, i.e., a pair $(A,\iota, \lambda) \in \M(n-r,r)(S)$ and $(E,\iota_0, \lambda_0) \in \M(1,0)(S)$, we consider the free $O_{\bk}$-module of finite rank 
$$V^{'}(E,A)=\Hom_{O_{\bk}}(E,A).$$
It is equipped with an $O_{\bk}$-valued  hermitian form $h^{'}$ given by $$h^{'}(x,y)=\iota_0^{-1}(\lambda_0^{-1}\circ y^{\vee}\circ \lambda \circ x),$$
where $y^{\vee}$ is the dual of $y.$ Given elements $x_1,...,x_m \in V^{'}(E,A)$, we define the {\em fundamental matrix} of the collection $x_1,...,x_m$ to be the hermitian $m \times m$ matrix over $ O_{\bk}$ with entries $h^{'}(x_i,x_j).$
\newline

Let now $T$ be an hermitian $m \times m$ matrix with entries in $O_{\bk}$. The  (global) {\em special cycle} $Z(T)$ is the stack of collections $(A,\iota, \lambda, E,\iota_0, \lambda_0, x_1,...,x_m)$, where $(A,\iota, \lambda) \in \M(n-r,r)(S)$, where $(E,\iota_0, \lambda_0) \in \M(1,0)(S)$, and  where $x_1,...,x_m \in V^{'}(E,A)$ such that the fundamental matrix of the $x_i$ is $T.$
Then $Z(T)$ is representable by a Deligne-Mumford-stack which is finite and unramified over $\M.$
\newline

We assume from  now on that $r=1$.

Given $T\in \Herm_n(O_{\bk})_{>0}$, one denotes by $\Diff_0(T)$ the set of primes which are inert in $\bk$ and for which $\ord_p ( \det T)$ is odd. Let $V_T=\bk^n$ with hermitian form given by $T$. In \cite{KR2}, Proposition 2.22, the following statements are shown.
\begin{itemize}
\item
 If the cardinality of  $\Diff_0(T)$ is at least $2$, then $Z(T)$ is empty.  

\item If $\Diff_0(T)=\{p\}$, then $Z(T)$ has support in  $\M_p^{\tilde{V},ss}$, where $\tilde{V}$ is the unique relevant hermitian space in $\mathcal{R}_{(n-1,1)}$ that is isomorphic to $V_T$ at all primes other than $p$ and $\infty$.
\item  If $\Diff_0(T)$ is empty, then $Z(T)$ has support in the union (of the supersingular locus) of the fibers at the ramified primes.
\end{itemize}
Let now $n=3$ and suppose   that $\Diff_0(T)=\{p\}$ for some $p>2$. 
We want to identify in this situation $\widehat{\deg}\ Z(T)$ (see introduction)   with the $T$-th Fourier coefficient of the derivative of a certain Eisenstein series (up to a constant independent of $T$). This was done in the case that $Z(T)$ has dimension $0$ in \cite{KR2}. Since in our situation  $Z(T)$ has support in  $\M_p^{\tilde{V},ss}$, we can use the above uniformizations to compute $\widehat{\deg}\ Z(T)$ using the local result of Theorem \ref{Hauptsatz}.

\section{The connection to Eisenstein series}

We first recall from \cite{KR2} some notions and facts on the relevant Eisenstein series. 

Let $V$ be any non-degenerate hermitian space over $\bk$ of dimension $m$. Let $G_1^V=U(V)$ and let $H=U(n,n).$ We fix a character $\eta$ of $\bk_{\A}^{\times}$ whose restriction to $\Q_{\A}^{\times}$ is $\chi^m$, where $\chi$ is the global quadratic character attached to $\bk$. One obtains a Weil representation of $G_1^V(\A)\times H(\A)$ on the Schwartz space $S(V(\A)^n)$. To $\varphi \in S(V(\A)^n)$ one attaches  the theta function 
$$
\theta(g,h; \varphi)= \sum_{x \in V(\Q)^n}\omega(h)\varphi(g^{-1}x).
$$
Let now $m=n$ and assume that $V$ has signature $(n,0).$ 
The theta integral is then $$
I(h;\varphi)=\int_{G_1^V(\Q)\setminus G_1^V(\A)}\theta(g,h;\varphi)dg,
$$ where the Haar measure is such that $\vol(G_1^V(\Q)\setminus G_1^V(\A))=1.$

Let $T$ be  non-singular and fix $x \in \Omega_T=\{x \in V^n \ | \ (x,x)=T \}.$ The $T$-th Fourier coefficient $I_T(h;\varphi)$ is of the form
$$
I_T(h,\varphi)=C \int_{ G_1^V(\A)}\omega(h)\varphi(g^{-1}x)dg,
$$
for some constant $C$ (independent of $T$ and $h$). Note that the stabilizer of $x$ in $G_1^V$ is trivial. 
\newline

Assume now that $\tilde{V}$ has signature $(1,n-1)$ and assume again that $\tilde{V}$ contains a self-dual lattice $\tilde{L}$. For non-singular $T\in \Herm_n({\bk})$, let $$\Diff(T,\tilde{V})=\{p < \infty \  | \  \chi_p(\det(T))=-\chi_p(\det(\tilde{V})) \}.$$ Let now $T\in \Herm_n(O_{\bk})_{>0}$ with $\Diff(T,\tilde{V})=\{p\}$ for an odd inert prime $p.$ Let $\tilde{V}^{'}=V_T$ (i.e. $\bk^n$ with hermitian form given by $T$). Then $\tilde{V}(\A_f^p)=\tilde{V}^{'}(\A_f^p).$ Let $\varphi^{p}\in S(\tilde{V}(\A_f^p)^n)$ be  the characteristic function of $(\tilde{L}\otimes \widehat{\Z}^p)^n$, let $\varphi_p$ be the characteristic function of $(\tilde{L}\otimes \Zp)^n$, and let $\varphi^{'}_p$ be the characteristic function of $(\Lambda^{\vee})^n$, where $\Lambda \subset \tilde{V}^{'}\otimes \Zp$ is a lattice of type $1$ as in section 1. Further, let  $\varphi_{\infty}(x)=e^{-2\pi \tr(x,x)}$ be the Gaussian associated to an hermitian space of signature $(n,0).$ To  $\varphi =\varphi_{\infty}\otimes \varphi_f^p\otimes \varphi_p$ one associates a Siegel-Weil section $\Phi$,  and to this one associates  an (incoherent) Eisenstein series $E(h,s,\tilde{V}^{\sharp}):=E(h,s, \Phi)$, where we  write as usual $\tilde{V}^{\sharp}=(\tilde{V},[[\tilde{L}]])$.  (See \cite{KR2} for more details on the definition of the Eisenstein series).  Let   $\varphi^{'} =\varphi_{\infty}\otimes \varphi_f^p\otimes \varphi^{'}_p$. 

For any $z\in M_n(\C)$ such that $v(z):=(2i)^{-1}(z-{}^t\bar{z})>0$, let  $u(z)=\frac{1}{2}(z+{}^t\bar{z})$ and let $$
h_z=
\begin{pmatrix}
1_n&u(z)\\
&1_n\\
 \end{pmatrix} \cdot \begin{pmatrix}
a&\\
&{}^t\bar{a}^{-1}\\
 \end{pmatrix}\in H(\R),
$$ where $a\in GL_n(\C)$ with $v(z)=a{}^t\bar{a}$. In the sequel we assume that  $h=h_{\infty}=h_z$ for some $z$ as above. 
As in \cite{KR2}, section 9, one has for the derivative in $s=0,$
$$
E^{'}_T(h_z,0,\tilde{V}^{\sharp})= \frac{W^{'}_{T,p}(0,\Phi_p)}{W_{T,p}(0,\Phi^{'}_p)}\cdot 2I_T(h_z,\varphi^{'}),
$$ where $\Phi^{'}_p$ is the Siegel-Weil section  associated to $\varphi^{'}_p$ (here we use that $W_{T,p}(0,\Phi^{'}_p) \neq 0$). 
We consider the Eisenstein series
$$
E(z,s,\Phi)=\eta_{\infty}(\det(a)^{-1})\det(v(z))^{-\frac{n}{2}}E(h_z,s,\Phi).
$$
Let $$E(z,s,\tilde{V})=\sum_{\tilde{V}^{\sharp}=(\tilde{V},[[L]])}E(z,s,\tilde{V}^{\sharp})$$ and let 
$$
E(z,s)=\sum_{V\in\mathcal{R}_{(n-r,r)}(\bk)}E(z,s,V).
$$
  From now on we assume again that $n=3.$

Let $h_{\bk}$ be the class number of $\bk$  and let $w_{\bk}=|O^{\times}_{\bk}|$. Given $\tilde{V}$ and $\tilde{V}^{'}$ as above, we define the constant $C_1=\frac{h_k}{w_k}\vol(U(\tilde{V}^{'})(\R))\vol(K_1)$, where $K_1$ is the stabilizer of a self-dual $O_{\bk}$-lattice in $\tilde{V}$ in the group $U(\tilde{V})(\A_f)$, and where the measures are   the Tamagawa measures with respect to a gauge form as in \cite{KR2}, section 8. 
(Note that $C_1$ does not depend on $T$.)

We are now ready to extend  Theorem 11.9 of \cite{KR2}  to all $T\in \Herm_3(O_{\bk})_{>0}$.

\begin{The}\label{GlobalerHS}
Let $T\in \Herm_3(O_{\bk})_{>0}$ with $Z(T)\neq\emptyset$ and let $p>2$ be a prime which is inert in $\bk$ and for which $\ord_p ( \det T)$ is odd. 
 Then
$$
E^{'}_T(z,0)=E^{'}_T(z,0,\tilde{V})=C_1\cdot \widehat{\deg}\ Z(T) \cdot q^T,  
$$ for any $z$ as above.
Here $\tilde{V}$ is the unique relevant hermitian space in $\mathcal{R}_{(2,1)}$ that is isomorphic to $V_T$ at all primes other than $p$ and $\infty$ and $q^T=e^{2\pi i \tr(Tz)}$ . \end{The}
\begin{Lem}\label{GloLweg}  In the situation of the theorem suppose that $T$ has diagonal blocks $T_1$ and $T_2$, where one of these blocks has size $2\times 2$ and hence the other has size $1\times 1$. Then  $\widehat{\deg}(T)=\chi(Z(T)_p, \mathcal{O}_{Z(T_1)}\otimes^{\mathbb{L}}\mathcal{O}_{Z(T_2)})\cdot \log p.$
\end{Lem}
{\em Proof.} 
Suppose, for example, that $T_1=(t_1)$ and that $T_2$ has diagonal entries $t_2,t_3$. Then it follows from Lemma \ref{Lweg} and \cite{KR2},  Proposition 6.3 that locally around any point of $Z(t_2)\times_{\M} Z(t_3)$  lying over the supersingular locus of $\M_p$ and having fundamental matrix $T_2$ we have $\mathcal{O}_{Z(t_2)}\otimes^{\mathbb{L}}\mathcal{O}_{Z(t_3)}=\mathcal{O}_{Z(t_2)}\otimes \mathcal{O}_{Z(t_3)}= \mathcal{O}_{Z(T_2)}$. From this the claim follows. \qed

\begin{Cor}
In the situation of the theorem suppose that $T$ has diagonal blocks $T_1,...,T_r$, where the size of $T_i$ is $m_i\times m_i$ and $\sum m_i=3$. Suppose further  that each cycle $Z(T_i)$ has pure dimension $n-m_i$. Let $\langle Z(T_1),...,Z(T_r)\rangle_T=\chi(Z(T)_p, \mathcal{O}_{Z(T_1)}\otimes^{\mathbb{L}}...\otimes^{\mathbb{L}}\mathcal{O}_{Z(T_r)})\cdot \log p.$  Then $$
E^{'}_T(z,0)=E^{'}_T(z,0,\tilde{V})=C_1\cdot \langle Z(T_1),...,Z(T_r)\rangle_T \cdot q^T,  
$$ for any $z$ as above.
\end{Cor}
{\em Proof.} This follows from Theorem \ref{GlobalerHS}, Lemma \ref{GloLweg} and (in the case $Z(T)$ has dimension $0$) from \cite{KR2}, Theorem 11.9 \qed
\newline

{\em Proof of the theorem.} As in \cite{KR2} we see that $E^{'}_T(z,0)=E^{'}_T(z,0,\tilde{V})$. It follows from our assumption on $p$  that $\supp(Z(T))\subset \M_p^{\tilde{V},ss}$,  see \cite{KR2}, Proposition 2.22. The completion of $Z(T)$ along the supersingular locus is a disjoint sum of the $\widehat{Z(T)}^{({V^{\sharp}},V_0),ss}:=Z(T)\times_{\M}\widehat{\M}^{({V^{\sharp}},V_0),ss},$ where the paris $({V^{\sharp}},V_0)$ run over $\mathcal{R}_{(2,1)}(\bk)^{\sharp}\times \mathcal{R}_{(1,0)}(\bk) $ such that $\Hom(V_0,V)\cong\tilde{V}.$ Let now  $({V^{\sharp}},V_0)$ be such a pair and write  as usual $V^{\sharp}=(V,L)$. We fix a base point $(A^o,\iota^o, \lambda^o, E^o,\iota^o_0, \lambda^o_0, x_1,x_2,x_3)$ in $\widehat{Z(T)}^{({V^{\sharp}},V_0),ss}(\F)$. Let $j_1,j_2,j_3$ be the special homomorphisms (in the sense of section 1) induced by $x_1,x_2,x_3$ between the $p$-divisble groups of $E^o$ and $A^o$. We fix an isomorphism $\eta^o: T^p(A^o)^0 \rightarrow V(\A_f^p)$ which is an isometry (w.r.t. the Weil Pairing on $T^p(A^o)^0$) and such that $\eta^o( T^p(A^o))=L\otimes \widehat{\Z}^p$. Likewise we choose $\eta_0^o$ for $E^o$ and a self-dual lattice $L_0$ in $V_0$.
Let $x^0=(x_1,x_2,x_3)$ and denote also by $x^0$ the induced element in $\Hom_{O_{\bk}\otimes \widehat{\Z}^p} (T^p(E^o),T^p(A^o))^3$. Then it follows from \cite{KR2}, Proposition 6.3 that the contribution to the intersection multiplicity $\widehat{\deg}\ Z(T)$ coming from $\widehat{Z}^{({{V}}^{\sharp},V_0),ss}(T) $ is of the form 
\[
\begin{split}
\chi_{\widehat{Z}^{({{V}}^{\sharp},V_0),ss}(T)} = \  &
(\mathcal{Z}(j_1),\mathcal{Z}(j_2),\mathcal{Z}(j_3))\times \\ &  |I^{V_0}(\Q)\setminus \{\   (gK^p,g_0K^p_0) \  | \ g^{-1}{x}^0g_0\in \Hom_{O_{\bk}\otimes \widehat{\Z}^p} (T^p(E^o),T^p(A^o))^3 \  \}|\log p.
\end{split}
\]
Here $g$ resp. $g_0$ are elements of $ G^V(\A_f^p)^0$ resp. of  $G^{V_0}(\A_f^p)^0$ and the actions on $x^0$ are through $\eta^o$ resp. $\eta^o_0$, i.e. $g^{-1}{x}^0g_0$ is an abbreviation for $({\eta^o}^{-1} g^{-1}\eta^o) {x}^0 ({\eta_0^o}^{-1} g_0\eta_0^o)$.  The group $I^{V_0}$ is the group  of quasi-isogenies of $E^o$ (respecting $\iota^o_0$ and $ \lambda^o$), further   
  $K^p \subset G^V(\A_f^p)^0$ is  the stabilizer of $L\otimes \widehat{\Z}^p$, and $K^p_0 $ is defined analogously.  We consider the lattice $ \tilde{L}:=\Hom_{O_{\bk}}(L_0,L)$ in $\tilde{V}.$ Note that $\Hom_{{\bk}\otimes \A_f^p} (T^p(E^o)^0,T^p(A^o)^0) \cong \tilde{V}(\A_f^p)$ and $   \Hom_{O_{\bk}\otimes \widehat{\Z}^p} (T^p(E^o),T^p(A^o))\cong \tilde{L}\otimes {\widehat{\Z}}^p$. We can identify the groups $G^V(\A_f^p)^0$ and  $G^{\tilde{V}^{}}(\A_f^p)^0$, and we can also identify $K^p$ with the stabilizer of $\tilde{L}\otimes \widehat{\Z}^p$ in  $G^{\tilde{V}^{}}(\A_f^p)^0$ (we denote it also by $K^p$). Let $K^p_1$ be the stabilizer of $\tilde{L}\otimes \widehat{\Z}^p $ in $G_1^{\tilde{V}^{}}(\A_f^p)=U(\tilde{V}^{})(\A_f^p)$. Let now as above $\varphi_f$  be the characteristic function of $(\tilde{L}\otimes \widehat{\Z})^3$, and let $\varphi_p^{'}, \varphi_{\infty}, \varphi$ and $\varphi^{'}$ also be as above.
Since on the one hand,  by Theorem \ref{Hauptsatz},  $ (\mathcal{Z}(j_1),\mathcal{Z}(j_2),\mathcal{Z}(j_3))=C_2\cdot W^{'}_{T,p}(0,\Phi_p)$ for some constant $C_2$,   
and  on the other hand  $G^{{V}^{}}(\A_f^p)^0/K^p=G^{\tilde{V}^{}}(\A_f^p)^0/K^p=G_1^{\tilde{V}^{}}(\A_f^p)/K_1^p$,  the above expression is of the form $$ \chi_{\widehat{Z}^{({{V}}^{\sharp},V_0),ss}(T)} =C_3\cdot W^{'}_{T,p}(0,\Phi_p) \cdot I_T(h_z,\varphi_f^p),$$ for some constant $C_3$.  Now 
$$E^{'}_T(z,0,\tilde{V})= \frac{W^{'}_{T,p}(0,\Phi_p)}{W_{T,p}(0,\Phi^{'}_p)}\cdot 2I_T(h_z,\varphi^{'})= C_4\cdot W^{'}_{T,p}(0,\Phi_p) \cdot \frac{I_T(h_z,\varphi_p^{'})}{W_{T,p}(0,\Phi^{'}_p)}\cdot I_T(h,\varphi_f^p)q^T,$$
for some constant $C_4$ (see also \cite{KR2}, (7.4)). The expression  $\frac{I_T(h_z,\varphi_p^{'})}{W_{T,p}(0,\Phi^{'}_p)}$ is a constant by the same reasoning as in \cite{KRY}, Proposition 5.3.3. Adding the contributions from all $(V^{\sharp},V_0)$
it follows that $$ E^{'}_T(z,0)=E^{'}_T(z,0,\tilde{V})=\tilde{C_1}\cdot \widehat{\deg}\ Z(T) \cdot q^T,$$  for some constant $\tilde{C_1}$. (Here we use that all pairs  $(V^{\sharp},V_0)$ give the same contribution which follows from the fact that by the reasoning above the contribution from  $(V^{\sharp},V_0)$ depends only on  $ I_T(h_z,\varphi_f^p)$  and hence only on   $\tilde{V}$. Note here that, since $n=3$ is odd, there is only one $U(\tilde{V})$ genus of self dual lattices.)  We have to show that $\tilde{C_1}=C_1.$ But by Theorem 11.9 in \cite{KR2}, this is true in the case that $T$ is not divisible by $p$ (as a matrix over $O_{\bk}$). Since the claim is independent of $T$, it follows that it is true. \qed
\begin{Rem}
{\em As already mentioned in the introduction, this proof shows that, in order to show the corresponding statement for any $n$, it is enough to show the corresponding local statement of Theorem \ref{Hauptsatz}. (In the case that  $n$ is even and that there are two $U(\tilde{V})$-genera of self-dual lattices,  one proceeds for each genus separately as above.) To conclude  the analogue of the corollary, we also proceed similarly as above:  Let $n$ be arbitrary, let  $Z(T)$ be a special cycle, where $T\in \Herm_m(O_{\bk})_{>0}$,  and suppose that $Z(T)$ has pure dimension $n-m$. Let $t_1,...,t_m$ be the diagonal entries of $T$. Then we see similarly as in the proof of the lemma that locally around any point of $Z(t_1)\times_{\M}...\times_{\M} Z(t_m)$  lying over the supersingular locus of $\M_p$ and having fundamental matrix $T$ there is the relation $\mathcal{O}_{Z(t_1)}\otimes^{\mathbb{L}} ...\otimes^{\mathbb{L}}\mathcal{O}_{Z(t_m)}=\mathcal{O}_{Z(t_1)}\otimes ...\otimes \mathcal{O}_{Z(t_m)}= \mathcal{O}_{Z(T)}$. Compare also the proof of \cite{KR2}, Proposition 11.6. Using this we can conclude the analogue of the corollary  if we assume the analogue of the theorem.}
\end{Rem}

\end{document}